\input amstex
\input amsppt.sty

\def\stydate{May 10, 2002}

\chardef\tempcat\catcode`\@ \ifx\undefined\amstexloaded\input amstex
\else\catcode`\@\tempcat\fi \expandafter\ifx\csname
amsppt.sty\endcsname\relax\input amsppt.sty \fi
\let\tempcat\undefined

\immediate\write16{This is LABEL.DEF by A.Degtyarev <\stydate>}
\expandafter\ifx\csname label.def\endcsname\relax\else
  \message{[already loaded]}\endinput\fi
\expandafter\edef\csname label.def\endcsname{%
  \catcode`\noexpand\@\the\catcode`\@\edef\noexpand\styname{LABEL.DEF}%
  \def\expandafter\noexpand\csname label.def\endcsname{\stydate}%
    \toks0{}\toks2{}}
\catcode`\@11
\def\labelmesg@ {LABEL.DEF: }
{\edef\temp{\the\everyjob\W@{\labelmesg@<\stydate>}}
\global\everyjob\expandafter{\temp}}

\def\@car#1#2\@nil{#1}
\def\@cdr#1#2\@nil{#2}
\def\eat@bs{\expandafter\eat@\string}
\def\eat@ii#1#2{}
\def\eat@iii#1#2#3{}
\def\eat@iv#1#2#3#4{}
\def\@DO#1#2\@{\expandafter#1\csname\eat@bs#2\endcsname}
\def\@N#1\@{\csname\eat@bs#1\endcsname}
\def\@Nx{\@DO\noexpand}
\def\@Name#1\@{\if\@undefined#1\@\else\@N#1\@\fi}
\def\@Ndef{\@DO\def}
\def\@Ngdef{\global\@Ndef}
\def\@Nedef{\@DO\edef}
\def\@Nxdef{\global\@Nedef}
\def\@Nlet{\@DO\let}
\def\@undefined#1\@{\@DO\ifx#1\@\relax\true@\else\false@\fi}
\def\@@addto#1#2{{\toks@\expandafter{#1#2}\xdef#1{\the\toks@}}}
\def\@@addparm#1#2{{\toks@\expandafter{#1{##1}#2}%
    \edef#1{\gdef\noexpand#1####1{\the\toks@}}#1}}
\def\make@letter{\edef\t@mpcat{\catcode`\@\the\catcode`\@}\catcode`\@11 }
\def\donext@{\expandafter\egroup\next@}
\def\x@notempty#1{\expandafter\notempty\expandafter{#1}}
\def\lc@def#1#2{\edef#1{#2}%
    \lowercase\expandafter{\expandafter\edef\expandafter#1\expandafter{#1}}}
\newif\iffound@
\def\find@#1\in#2{\found@false
    \DNii@{\ifx\next\@nil\let\next\eat@\else\let\next\nextiv@\fi\next}%
    \edef\nextiii@{#1}\def\nextiv@##1,{%
    \edef\next{##1}\ifx\nextiii@\next\found@true\fi\FN@\nextii@}%
    \expandafter\nextiv@#2,\@nil}
{\let\head\relax\let\specialhead\relax\let\subhead\relax
\let\subsubhead\relax\let\proclaim\relax
\gdef\let@relax{\let\head\relax\let\specialhead\relax\let\subhead\relax
    \let\subsubhead\relax\let\proclaim\relax}}
\newskip\@savsk
\let\@ignorespaces\ignorespaces
\def\@ignorespacesp{\ifhmode
  \ifdim\lastskip>\z@\else\penalty\@M\hskip-1sp%
        \penalty\@M\hskip1sp \fi\fi\@ignorespaces}
\def\ignorespaces{\protect\@ignorespacesp}
\def\@bsphack{\relax\ifmmode\else\@savsk\lastskip
  \ifhmode\edef\@sf{\spacefactor\the\spacefactor}\fi\fi}
\def\@esphack{\relax
  \ifx\penalty@\penalty\else\penalty\@M\fi   
  \ifmmode\else\ifhmode\@sf{}\ifdim\@savsk>\z@\@ignorespacesp\fi\fi\fi}
\let\@frills@\identity@
\let\@txtopt@\identyty@
\newif\if@star
\newif\if@write\@writetrue
\def\@numopt@{\if@star\expandafter\eat@\fi}
\def\checkstar@#1{\DN@{\@writetrue
  \ifx\next*\DN@####1{\@startrue\checkstar@@{#1}}%
      \else\DN@{\@starfalse#1}\fi\next@}\FN@\next@}
\def\checkstar@@#1{\DN@{%
  \ifx\next*\DN@####1{\@writefalse#1}%
      \else\DN@{\@writetrue#1}\fi\next@}\FN@\next@}
\def\checkfrills@#1{\DN@{%
  \ifx\next\nofrills\DN@####1{#1}\def\@frills@####1{####1\nofrills}%
      \else\DN@{#1}\let\@frills@\identity@\fi\next@}\FN@\next@}
\def\checkbrack@#1{\DN@{%
    \ifx\next[\DN@[####1]{\def\@txtopt@########1{####1}#1}%
    \else\DN@{\let\@txtopt@\identity@#1}\fi\next@}\FN@\next@}
\def\check@therstyle#1#2{\bgroup\DN@{#1}\ifx\@txtopt@\identity@\else
        \DNii@##1\@therstyle{}\def\@therstyle{\DN@{#2}\nextii@}%
    \expandafter\expandafter\expandafter\nextii@\@txtopt@\@therstyle.\@therstyle
    \fi\donext@}

\newread\@inputcheck
\def\@input#1{\openin\@inputcheck #1 \ifeof\@inputcheck \W@
  {No file `#1'.}\else\closein\@inputcheck \relax\input #1 \fi}

\def\loadstyle#1{\edef\next{#1}%
    \DN@##1.##2\@nil{\if\notempty{##2}\else\def\next{##1.sty}\fi}%
    \expandafter\next@\next.\@nil\lc@def\next@\next
    \expandafter\ifx\csname\next@\endcsname\relax\input\next\fi}

\let\pagebody@\pagebody
\let\pagetop@\empty
\let\pagebot@\empty
\let\@Xend\empty
\def\pagebody{\pagetop@\pagebody@\pagebot@\@Xend}
\let\@Xclose\empty

\newwrite\@Xmain
\newwrite\@Xsub
\def\W@X{\write\@Xout}
\def\make@Xmain{\global\let\@Xout\@Xmain\global\let\end\endmain@
  \xdef\@Xname{\jobname}\xdef\@inputname{\jobname}}
\begingroup
\catcode`\(\the\catcode`\{\catcode`\{12
\catcode`\)\the\catcode`\}\catcode`\}12
\gdef\W@count#1((\lc@def\@tempa(#1)%
    \def\\##1(\W@X(\global##1\the##1))%
    \edef\@tempa(\W@X(%
        \string\expandafter\gdef\string\csname\space\@tempa\string\endcsname{)%
        \\\pageno\\\cnt@toc\\\cnt@idx\\\cnt@glo\\\footmarkcount@
        \@Xclose\W@X(}))\expandafter)\@tempa)
\endgroup
\def\readaux{\bgroup\checkbrack@\readaux@}
\let\begin\readaux
\def\readaux@{%
    \W@{>>> \labelmesg@ Run this file twice to get x-references right}%
    \global\everypar{}%
    {\def\\{\global\let}%
        \def\/##1##2{\gdef##1{\wrn@command##1##2}}%
        \disablepreambule@cs}%
    \make@Xmain{\make@letter\setboxz@h{\@input{\@txtopt@{\@Xname.aux}}%
            \lc@def\@tempa\jobname\@Name\open@\@tempa\@}}%
  \immediate\openout\@Xout\@Xname.aux%
    \immediate\W@X{\relax}\egroup}
\everypar{\global\everypar{}\readaux} {\toks@\expandafter{\topmatter}
\global\edef\topmatter{\noexpand\readaux\the\toks@}}
\let\@@end@@\end

\def\@Xclose@{{\def\@Xend{\ifnum\insertpenalties=\z@
        \W@count{close@\@Xname}\closeout\@Xout\fi}%
    \vfill\supereject}}
\def\endmain@{\@Xclose@
    \W@{>>> \labelmesg@ Run this file twice to get x-references right}%
    \@@end@@}
\def\disablepreambule@cs{\\\disablepreambule@cs\relax}

\def\include#1{\bgroup
  \ifx\@Xout\@Xsub\DN@{\errmessage
        {\labelmesg@ Only one level of \string\include\space is supported}}%
    \else\edef\@tempb{#1}\clearpage
      \DN@##1 {\if\notempty{##1}\edef\@tempb{##1}\DN@####1\eat@ {}\fi\next@}%
    \DNii@##1.{\edef\@tempa{##1}\DN@####1\eat@.{}\next@}%
        \expandafter\next@\@tempb\eat@{} \eat@{} %
    \expandafter\nextii@\@tempb.\eat@.%
        \relaxnext@
      \if\x@notempty\@tempa
          \edef\nextii@{\write\@Xmain{%
            \noexpand\string\noexpand\@input{\@tempa.aux}}}\nextii@
        \ifx\undefined\@includelist\found@true\else
                    \find@\@tempa\in\@includelist\fi
            \iffound@\ifx\undefined\@noincllist\found@false\else
                    \find@\@tempb\in\@noincllist\fi\else\found@true\fi
            \iffound@\lc@def\@tempa\@tempa
                \if\@undefined\close@\@tempa\@\else\edef\next@{\@Nx\close@\@tempa\@}\fi
            \else\xdef\@Xname{\@tempa}\xdef\@inputname{\@tempb}%
                \W@count{open@\@Xname}\global\let\@Xout\@Xsub
            \openout\@Xout\@tempa.aux \W@X{\relax}%
            \DN@{\let\end\endinput\@input\@inputname
                    \@Xclose@\make@Xmain}\fi\fi\fi
  \donext@}
\def\includeonly#1{\edef\@includelist{#1}}
\def\noinclude#1{\edef\@noincllist{#1}}

\def\arabicnum#1{\number#1}

\def\Romannum#1{\expandafter\uppercase\expandafter{\romannumeral#1}}
\def\alphnum#1{\ifcase#1\or a\or b\or c\or d\else\@ialph{#1}\fi}
\def\@ialph#1{\ifcase#1\or \or \or \or \or e\or f\or g\or h\or i\or j\or
    k\or l\or m\or n\or o\or p\or q\or r\or s\or t\or u\or v\or w\or x\or y\or
    z\else\fi}
\def\Alphnum#1{\ifcase#1\or A\or B\or C\or D\else\@Ialph{#1}\fi}
\def\@Ialph#1{\ifcase#1\or \or \or \or \or E\or F\or G\or H\or I\or J\or
    K\or L\or M\or N\or O\or P\or Q\or R\or S\or T\or U\or V\or W\or X\or Y\or
    Z\else\fi}

\def\ST@P{step}
\def\ST@LE{style}
\def\N@M{no}
\def\F@NT{font@}
\outer\def\newcounter{\checkbrack@{\expandafter\newcounter@\@txtopt@{{}}}}
{\let\newcount\relax
\gdef\newcounter@#1#2#3{{%
    \toks@@\expandafter{\csname\eat@bs#2\N@M\endcsname}%
    \DN@{\alloc@0\count\countdef\insc@unt}%
    \ifx\@txtopt@\identity@\expandafter\next@\the\toks@@
        \else\if\notempty{#1}\global\@Nlet#2\N@M\@#1\fi\fi
    \@Nxdef\the\eat@bs#2\@{\if\@undefined\the\eat@bs#3\@\else
            \@Nx\the\eat@bs#3\@.\fi\noexpand\arabicnum\the\toks@@}%
  \@Nxdef#2\ST@P\@{}%
  \if\@undefined#3\ST@P\@\else
    \edef\next@{\noexpand\@@addto\@Nx#3\ST@P\@{%
             \global\@Nx#2\N@M\@\z@\@Nx#2\ST@P\@}}\next@\fi
    \expandafter\@@addto\expandafter\@Xclose\expandafter
        {\expandafter\\\the\toks@@}}}}
\outer\def\copycounter#1#2{%
    \@Nxdef#1\N@M\@{\@Nx#2\N@M\@}%
    \@Nxdef#1\ST@P\@{\@Nx#2\ST@P\@}%
    \@Nxdef\the\eat@bs#1\@{\@Nx\the\eat@bs#2\@}}
\outer\def\everystep{\checkstar@\everystep@}
\def\everystep@#1{\if@star\let\next@\gdef\else\let\next@\@@addto\fi
    \@DO\next@#1\ST@P\@}
\def\counterstyle#1{\@Ngdef\the\eat@bs#1\@}
\def\advancecounter#1#2{\@N#1\ST@P\@\global\advance\@N#1\N@M\@#2}
\def\setcounter#1#2{\@N#1\ST@P\@\global\@N#1\N@M\@#2}
\def\counter#1{\refstepcounter#1\printcounter#1}
\def\printcounter#1{\@N\the\eat@bs#1\@}
\def\refcounter#1{\xdef\@lastmark{\printcounter#1}}
\def\stepcounter#1{\advancecounter#1\@ne}
\def\refstepcounter#1{\stepcounter#1\refcounter#1}
\def\savecounter#1{\@Nedef#1@sav\@{\global\@N#1\N@M\@\the\@N#1\N@M\@}}
\def\restorecounter#1{\@Name#1@sav\@}

\def\warning#1#2{\W@{Warning: #1 on input line #2}}
\def\warning@#1{\warning{#1}{\the\inputlineno}}
\def\wrn@@Protect#1#2{\warning@{\string\Protect\string#1\space ignored}}
\def\wrn@@label#1#2{\warning{label `#1' multiply defined}{#2}}
\def\wrn@@ref#1#2{\warning@{label `#1' undefined}}
\def\wrn@@cite#1#2{\warning@{citation `#1' undefined}}
\def\wrn@@command#1#2{\warning@{Preamble command \string#1\space ignored}#2}
\def\wrn@@option#1#2{\warning@{Option \string#1\string#2\space is not supported}}
\def\wrn@@reference#1#2{\W@{Reference `#1' on input line \the\inputlineno}}
\def\wrn@@citation#1#2{\W@{Citation `#1' on input line \the\inputlineno}}
\let\wrn@reference\eat@ii
\let\wrn@citation\eat@ii
\def\nowarning#1{\if\@undefined\wrn@\eat@bs#1\@\wrn@option\nowarning#1\else
        \@Nlet\wrn@\eat@bs#1\@\eat@ii\fi}
\def\printwarning#1{\if\@undefined\wrn@@\eat@bs#1\@\wrn@option\printwarning#1\else
        \@Nlet\wrn@\eat@bs#1\expandafter\@\csname wrn@@\eat@bs#1\endcsname\fi}
\printwarning\Protect \printwarning\label \printwarning\ref
\printwarning\cite \printwarning\command \printwarning\option

{\catcode`\#=12\gdef\@lH{#}}
\def\@@HREF#1{}
\def\@HREF#1#2{\@@HREF{a #1}{\let\@@HREF\eat@#2}\@@HREF{/a}}
\def\@@Hf#1{file:#1} \let\@Hf\@@Hf
\def\@@Hl#1{\@lH#1} \let\@Hl\@@Hl
\def\@@Hname#1{\@HREF{name="#1"}{}} \let\@Hname\@@Hname
\def\@@Href#1{\@HREF{href="#1"}} \let\@Href\@@Href
\ifx\undefined\pdfoutput
  \csname newcount\endcsname\pdfoutput
\else
  \def\pdflinkattr{attr{/C [0 0.9 0.9]}}
  \def\@pdfHf#1{file{#1}}
  \def\@pdfHl#1{name{#1}}
  \def\@pdfHname#1{\pdfdest name{#1}xyz\relax}
  \def\@pdfHref#1#2{\pdfstartlink \pdflinkattr goto #1\relax{\let\@Href\eat@
    #2}\pdfendlink}
  \def\@ifpdf#1#2{\ifnum\pdfoutput>\z@\expandafter#1\else\expandafter#2\fi}
  \def\@Hf{\@ifpdf\@pdfHf\@@Hf}
  \def\@Hl{\@ifpdf\@pdfHl\@@Hl}
  \def\@Hname{\@ifpdf\@pdfHname\@@Hname}
  \def\@Href{\@ifpdf\@pdfHref\@@Href}
\fi
\def\@Hr#1#2{\if\notempty{#1}\@Hf{#1}\fi\if\notempty{#2}\@Hl{#2}\fi}
\def\@localHref#1{\@Href{\@Hr{}{#1}}}
\def\@countlast#1{\@N#1last\@}
\def\@@countref#1#2{\global\advance#2\@ne
  \@Nxdef#2last\@{\the#2}\@tocHname{#1\@countlast#2}}
\def\@countref#1{\@DO\@@countref#1@HR\@#1}

\def\Href@@#1{\@N\Href@-#1\@}
\def\Href@#1#2{\@N\Href@-#1\@{\@Hl{@#1-#2}}}
\def\Hname@#1{\@N\Hname@-#1\@}
\def\Hlast@#1{\@N\Hlast@-#1\@}
\def\cntref@#1{\global\@DO\advance\cnt@#1\@\@ne
  \@Nxdef\Hlast@-#1\@{\@DO\the\cnt@#1\@}\Hname@{#1}{@#1-\Hlast@{#1}}}
\def\HyperRefs#1{\global\@Nlet\Hlast@-#1\@\empty
  \global\@Nlet\Hname@-#1\@\@Hname
  \global\@Nlet\Href@-#1\@\@Href}
\def\NoHyperRefs#1{\global\@Nlet\Hlast@-#1\@\empty
  \global\@Nlet\Hname@-#1\@\eat@
  \global\@Nlet\Href@-#1\@\eat@}

\HyperRefs{label} {\catcode`\-11 \gdef\@labelref#1{\Hname@-label{r@-#1}}
\gdef\@xHref#1{\Href@-label{\@Hl{r@-#1}}} }
\HyperRefs{toc}
\def\@HR#1{\if\notempty{#1}\string\@HR{\Hlast@{toc}}{#1}\else{}\fi}



\def\bftext{\ifmmode\fam\bffam\else\bf\fi}
\let\@lastmark\empty
\let\@lastlabel\empty
\def\lastmark{\@lastmark}
\let\lastlabel\empty
\let\everylabel\relax
\let\everylabel@\eat@
\let\everyref\relax
\def\newlabel{\bgroup\everylabel\newlabel@}
\def\newlabel@#1#2#3{\if\@undefined\r@-#1\@\else\wrn@label{#1}{#3}\fi
  {\let\protect\noexpand\@Nxdef\r@-#1\@{#2}}\egroup}
\def\w@ref{\bgroup\everyref\w@@ref}
\def\w@@ref#1#2#3#4{%
  \if\@undefined\r@-#1\@{\bftext??}#2{#1}{}\else%
   \@xHref{#1}{\@DO{\expandafter\expandafter#3}\r@-#1\@\@nil}\fi
  #4{#1}{}\egroup}
\def\@@@xref#1{\w@ref{#1}\wrn@ref\@car\wrn@reference}
\def\@xref#1{\rom{\@@@xref{#1}}}
\let\xref\@xref
\def\pageref#1{\w@ref{#1}\wrn@ref\@cdr\wrn@reference}
\def\thepage{\ifnum\pageno<\z@\romannumeral-\pageno\else\number\pageno\fi}
\def\label@{\@bsphack\bgroup\everylabel\label@@}
\def\label@@#1#2{\everylabel@{{#1}{#2}}%
  \@labelref{#2}%
  \let\thepage\relax
  \def\protect{\noexpand\noexpand\noexpand}%
  \edef\@tempa{\edef\noexpand\@lastlabel{#1}%
    \W@X{\string\newlabel{#2}{{\@lastmark}{\thepage}}{\the\inputlineno}}}%
  \expandafter\egroup\@tempa\@esphack}
\def\label#1{\label@{#1}{#1}}
\def\fn@P@{\relaxnext@
    \DN@{\ifx[\next\DN@[####1]{}\else
        \ifx"\next\DN@"####1"{}\else\DN@{}\fi\fi\next@}%
    \FN@\next@}
\def\eat@fn#1{\ifx#1[\expandafter\eat@br\else
  \ifx#1"\expandafter\expandafter\expandafter\eat@qu\fi\fi}
\def\eat@br#1]#2{}
\def\eat@qu#1"#2{}
{\catcode`\~\active\lccode`\~`\@
\lowercase{\global\let\@@P@~\gdef~{\protect\@@P@}}}
\def\Protect@@#1{\def#1{\protect#1}}
\def\disable@special{\let\W@X@\eat@iii\let\label\eat@
    \def\footnotemark{\protect\fn@P@}%
  \let\footnotetext\eat@fn\let\footnote\eat@fn
    \let\refcounter\eat@\let\savecounter\eat@\let\restorecounter\eat@
    \let\advancecounter\eat@ii\let\setcounter\eat@ii
  \let\ifvmode\iffalse\Protect@@\@@@xref\Protect@@\pageref\Protect@@\nofrills
    \Protect@@\\\Protect@@~}
\let\notoctext\identity@
\def\W@X@#1#2#3{\@bsphack{\disable@special\let\notoctext\eat@
    \def\chapter{\protect\chapter@toc}\let\thepage\relax
    \def\protect{\noexpand\noexpand\noexpand}#1%
  \edef\next@{\if\@undefined#2\@\else\write#2{#3}\fi}\expandafter}\next@
    \@esphack}
\newcount\cnt@toc
\def\writeauxline#1#2#3{\W@X@{\cntref@{toc}\let\tocref\@HR}
  \@Xout{\string\@Xline{#1}{#2}{#3}{\thepage}}}
{\let\newwrite\relax
\gdef\@openin#1{\make@letter\@input{\jobname.#1}\t@mpcat}
\gdef\@openout#1{\global\expandafter\newwrite\csname tf@-#1\endcsname
   \immediate\openout\@N\tf@-#1\@\jobname.#1\relax}}
\def\@@openout#1{\@openout{#1}%
  \@@addto\readaux@{\immediate\closeout\@N\tf@-#1\@}}
\def\auxlinedef#1{\@Ndef\do@-#1\@}
\def\@Xline#1{\if\@undefined\do@-#1\@\expandafter\eat@iii\else
    \@DO\expandafter\do@-#1\@\fi}
\def\beginW@{\bgroup\def\do##1{\catcode`##112 }\dospecials\do\@\do\"
    \catcode`\{\@ne\catcode`\}\tw@\immediate\write\@N}
\def\endW@toc#1#2#3{{\string\tocline{#1}{#2\string\page{#3}}}\egroup}
\def\do@tocline#1{%
    \if\@undefined\tf@-#1\@\expandafter\eat@iii\else
        \beginW@\tf@-#1\@\expandafter\endW@toc\fi
} \auxlinedef{toc}{\do@tocline{toc}}

\let\protect\empty
\def\Protect#1{\if\@undefined#1@P@\@\PROTECT#1\else\wrn@Protect#1\empty\fi}
\def\PROTECT#1{\@Nlet#1@P@\@#1\edef#1{\noexpand\protect\@Nx#1@P@\@}}
\def\pdef#1{\edef#1{\noexpand\protect\@Nx#1@P@\@}\@Ndef#1@P@\@}

\Protect\operatorname \Protect\operatornamewithlimits \Protect\qopname@
\Protect\qopnamewl@ \Protect\text \Protect\topsmash \Protect\botsmash
\Protect\smash \Protect\widetilde \Protect\widehat \Protect\thetag
\Protect\therosteritem
\Protect\Cal \Protect\Bbb \Protect\bold \Protect\slanted \Protect\roman
\Protect\italic \Protect\boldkey \Protect\boldsymbol \Protect\frak
\Protect\goth \Protect\dots
\Protect\cong \Protect\lbrace \let\{\lbrace \Protect\rbrace \let\}\rbrace
\let\root@P@@\root \def\root@P@#1{\root@P@@#1\of}
\def\root#1\of{\protect\root@P@{#1}}

\def\frills{\ignorespaces\@txtopt@}
\def\frillsnotempty#1{\x@notempty{\@txtopt@{#1}}}
\def\numberline{\@numopt@}
\newif\if@theorem
\let\@therstyle\eat@
\def\@headtext@#1#2{{\disable@special\let\protect\noexpand
    \def\chapter{\protect\chapter@rh}%
    \edef\next@{\noexpand\@frills@\noexpand#1{#2}}\expandafter}\next@}
\let\AmSrighthead@\rightheadtext
\def\rightheadtext{\checkfrills@{\@headtext@\AmSrighthead@}}
\let\AmSlefthead@\leftheadtext
\def\leftheadtext{\checkfrills@{\@headtext@\AmSlefthead@}}
\def\@head@@#1#2#3#4#5{\@Name\pre\eat@bs#1\@\if@theorem\else
    \@frills@{\csname\expandafter\eat@iv\string#4\endcsname}\relax
        \ifx\protect\empty\@N#1\F@NT\@\fi\fi
    \@N#1\ST@LE\@{\counter#3}{#5}%
  \if@write\writeauxline{toc}{\eat@bs#1}{#2{\counter#3}\@HR{#5}}\fi
    \if@theorem\else\expandafter#4\fi
    \ifx#4\endhead\ifx\@txtopt@\identity@\else
        \headmark{\@N#1\ST@LE\@{\counter#3}{\frills\empty}}\fi\fi
    \@Name\post\eat@bs#1\@\ignorespaces}
\ifx\undefined\endhead\Invalid@\endhead\fi
\def\@head@#1{\checkstar@{\checkfrills@{\checkbrack@{\@head@@#1}}}}
\def\@thm@@#1#2#3{\@Name\pre\eat@bs#1\@
    \@frills@{\csname\expandafter\eat@iv\string#3\endcsname}
    {\@theoremtrue\check@therstyle{\@N#1\ST@LE\@}\frills
            {\counter#2}\@theoremfalse}%
    \@DO\envir@stack\end\eat@bs#1\@
    \@N#1\F@NT\@\@Name\post\eat@bs#1\@\ignorespaces}
\def\@thm@#1{\checkstar@{\checkfrills@{\checkbrack@{\@thm@@#1}}}}
\def\@capt@@#1#2#3#4#5\endcaption{\bgroup
    \edef\@tempb{\global\footmarkcount@\the\footmarkcount@
    \global\@N#2\N@M\@\the\@N#2\N@M\@}%
    \def\shortcaption##1{\global\def\sh@rtt@xt####1{##1}}\let\sh@rtt@xt\identity@
    \DN@{#4{\@tempb\@N#1\ST@LE\@{\counter#2}}}%
    \if\notempty{#5}\DNii@{\next@\@N#1\F@NT\@}\else\let\nextii@\next@\fi
    \nextii@#5\endcaption
  \if@write\writeauxline{#3}{\eat@bs#1}{{} \@HR{\@N#1\ST@LE\@{\counter#2}%
    \if\notempty{#5}.\enspace\fi\sh@rtt@xt{#5}}}\fi
  \global\let\sh@rtt@xt\undefined\egroup}
\def\@capt@#1{\checkstar@{\checkfrills@{\checkbrack@{\@capt@@#1}}}}
\let\captiontextfont@\empty

\ifx\undefined\subsubheadfont@\def\subsubheadfont@{\it}\fi
\ifx\undefined\proclaimfont\def\proclaimfont{\sl}\fi
\ifx\undefined\proclaimfont@\let\proclaimfont@\proclaimfont\fi
\def\proclaimfont{\proclaimfont@}
\ifx\undefined\definitionfont@\def\AmSdeffont@{\rm}
    \else\let\AmSdeffont@\definitionfont@\fi
\ifx\undefined\remarkfont@\def\remarkfont@{\rm}\fi

\def\newfont@def#1#2{\if\@undefined#1\F@NT\@
    \@Nxdef#1\F@NT\@{\@Nx.\expandafter\eat@iv\string#2\F@NT\@}\fi}
\def\newhead@#1#2#3#4{{%
    \gdef#1{\@therstyle\@therstyle\@head@{#1#2#3#4}}\newfont@def#1#4%
    \if\@undefined#1\ST@LE\@\@Ngdef#1\ST@LE\@{\headstyle}\fi
    \if\@undefined#2\@\gdef#2{\headtocstyle}\fi
  \@@addto\moretocdefs@{\\#1#1#4}}}
\outer\def\newhead#1{\checkbrack@{\expandafter\newhead@\expandafter
    #1\@txtopt@\headtocstyle}}
\outer\def\newtheorem#1#2#3#4{{%
    \gdef#2{\@thm@{#2#3#4}}\newfont@def#2#4%
    \@Nxdef\end\eat@bs#2\@{\noexpand\revert@envir
        \@Nx\end\eat@bs#2\@\noexpand#4}%
  \if\@undefined#2\ST@LE\@\@Ngdef#2\ST@LE\@{\proclaimstyle{#1}}\fi}}%
\outer\def\newcaption#1#2#3#4#5{{\let#2\relax
  \edef\@tempa{\gdef#2####1\@Nx\end\eat@bs#2\@}%
    \@tempa{\@capt@{#2#3{#4}#5}##1\endcaption}\newfont@def#2\endcaptiontext%
  \if\@undefined#2\ST@LE\@\@Ngdef#2\ST@LE\@{\captionstyle{#1}}\fi
  \@@addto\moretocdefs@{\\#2#2\endcaption}\newtoc{#4}}}
{
\outer\gdef\newtoc#1{{%
    \@DO\ifx\do@-#1\@\relax
    \global\auxlinedef{#1}{\do@tocline{#1}}{}%
    \@@addto\tocsections@{\make@toc{#1}{}}\fi}}}

\toks@\expandafter{\itembox@}
\toks@@{\bgroup\let\therosteritem\identity@\let\rm\empty
  \let\@Href\eat@\let\@Hname\eat@
  \edef\next@{\edef\noexpand\@lastmark{\therosteritem@}}\donext@}
\edef\itembox@{\the\toks@@\the\toks@}
\def\firstitem@false{\let\iffirstitem@\iffalse
    \global\let\lastlabel\@lastlabel}

\let\rosteritemrefform\therosteritem
\let\rosteritemrefseparator\empty
\def\rosteritemref#1{\hbox{\rosteritemrefform{\@@@xref{#1}}}}
\def\local#1{\label@\@lastlabel{\lastlabel-i#1}}
\def\loccit#1{\rosteritemref{\lastlabel-i#1}}
\def\xRef@P@{\gdef\lastlabel}
\def\xRef#1{\@xref{#1}\protect\xRef@P@{#1}}

\def\iref@P@{\gdef\lastref}
\def\itemref#1#2{\rosteritemref{#1-i#2}\protect\iref@P@{#1}}
\def\iref#1{\@xref{#1}\rosteritemrefseparator\itemref{#1}}
\def\ditto#1{\rosteritemref{\lastref-i#1}}

\def\eqtag{\tag\counter\equation}
\def\eqref#1{\thetag{\@@@xref{#1}}}
\def\tagform@#1{\ifmmode\hbox{\rm\else\rom{\fi
        (\ignorespaces#1\unskip)\iftrue}\else}\fi}

\let\AmSfnote@\makefootnote@
\def\makefootnote@#1{\bgroup\let\footmarkform@\identity@
  \edef\next@{\edef\noexpand\@lastmark{#1}}\donext@\AmSfnote@{#1}}

\def\clearpage{\ifnum\insertpenalties>0\line{}\fi\vfill\supereject}

\def\proof{\checkfrills@{\checkbrack@{%
    \check@therstyle{\@frills@{\demo}{\frills{Proof}}{}}
        {\frills{}\envir@stack\endremark\envir@stack\enddemo}%
  \envir@stack\endproof\ignorespaces}}}
\def\endproof{\nofrillscheck{\frills@{\qed}\revert@envir\endproof\enddemo}}

\let\AmSref\ref
\let\AmSrefstyle\refstyle
\let\plaincite\cite
\def\citei@#1,{\citeii@#1\eat@,}
\def\citeii@#1\eat@{\w@ref{#1}\wrn@cite\@car\wrn@citation}
\def\mcite@#1;{\plaincite{\citei@#1\eat@,\unskip}\mcite@i}
\def\mcite@i#1;{\DN@{#1}\ifx\next@\endmcite@
  \else, \plaincite{\citei@#1\eat@,\unskip}\expandafter\mcite@i\fi}
\def\endmcite@{\endmcite@}
\def\cite#1{\mcite@#1;\endmcite@;}
\PROTECT\cite
\def\refstyle#1{\AmSrefstyle{#1}\uppercase{%
    \ifx#1A\relax \def\@ref@##1{\AmSref\xdef\@lastmark{##1}\key##1}%
    \else\ifx#1C\relax \def\@ref@{\AmSref\no\counter\refno}%
        \else\def\@ref@{\AmSref}\fi\fi}}
\refstyle A
\newcounter\refno\null
\newif\ifRefs
\gdef\Refs{\checkstar@{\checkbrack@{\csname AmSRefs\endcsname
  \nofrills{\frills{References}%
  \if@write\writeauxline{toc}{vartocline}{\@HR{\frills{References}}}\fi}%
  \def\ref{\@ref@}\Refstrue\ignorespaces}}}
\let\ref\xref

\newif\iftoc
\pdef\tocbreak{\iftoc\hfil\break\fi}
\def\tocsections@{\make@toc{toc}{}}
\let\moretocdefs@\empty
\def\newtocline@#1#2#3{%
  \edef#1{\def\@Nx#2line\@####1{\@Nx.\expandafter\eat@iv
        \string#3\@####1\noexpand#3}}%
  \@Nedef\no\eat@bs#1\@{\let\@Nx#2line\@\noexpand\eat@}%
    \@N\no\eat@bs#1\@}
\def\MakeToc#1{\@@openout{#1}}
\def\newtocline#1#2#3{\Err@{\Invalid@@\string\newtocline}}
\def\make@toc#1#2{\penaltyandskip@{-200}\aboveheadskip
    \if\notempty{#2}
        \centerline{\headfont@\ignorespaces#2\unskip}\nobreak
    \vskip\belowheadskip \fi
    \@openin{#1}\relax
    \vskip\z@}
\def\contents{\readaux\checkfrills@{\checkbrack@{\@contents@}}}
\def\@contents@{\toc@{\frills{Contents}}\envir@stack\endcontents%
    \def\nopagenumbers{\let\page\eat@}\let\newtocline\newtocline@\toctrue
  \def\@HR{\Href@{toc}}%
  \def\tocline##1{\csname##1line\endcsname}
  \edef\caption##1\endcaption{\expandafter\noexpand
    \csname head\endcsname##1\noexpand\endhead}%
    \ifmonograph@\def\vartoclineline{\Chapterline}%
        \else\def\vartoclineline##1{\sectionline{{} ##1}}\fi
  \let\\\newtocline@\moretocdefs@
    \ifx\@frills@\identity@\def\\##1##2##3{##1}\moretocdefs@
        \else\let\tocsections@\relax\fi
    \def\\{\unskip\space\ignorespaces}\let\maketoc\make@toc}
\def\endcontents{\tocsections@\vskip-\lastskip\revert@envir\endcontents
    \endtoc}

\if\@undefined\selectf@nt\@\let\selectf@nt\identity@\fi
\def\Err@math#1{\Err@{Use \string#1\space only in text}}
\def\textonlyfont@#1#2{%
    \def#1{\RIfM@\Err@math#1\else\edef\f@ntsh@pe{\string#1}\selectf@nt#2\fi}%
    \PROTECT#1}
\tenpoint

\def\newshapeswitch#1#2{\gdef#1{\selectsh@pe#1#2}\PROTECT#1}
\def\shapeswitch#1#2#3{\@Ngdef#1\string#2\@{#3}}
\shapeswitch\rm\bf\bf \shapeswitch\rm\tt\tt \shapeswitch\rm\smc\smc
\newshapeswitch\em\it
\shapeswitch\em\it\rm \shapeswitch\em\sl\rm
\def\selectsh@pe#1#2{\relax\if\@undefined#1\f@ntsh@pe\@#2\else
    \@N#1\f@ntsh@pe\@\fi}

\def\@itcorr@{\leavevmode
    \edef\prevskip@{\ifdim\lastskip=\z@ \else\hskip\the\lastskip\relax\fi}\unskip
    \edef\prevpenalty@{\ifnum\lastpenalty=\z@ \else
        \penalty\the\lastpenalty\relax\fi}\unpenalty
    \/\prevpenalty@\prevskip@}
\def\rom@P@#1{\@itcorr@{\selectsh@pe\rm\rm#1}}
\def\rom{\protect\rom@P@}
\def\Rom@P@#1{\@itcorr@{\rm#1}}
\def\Rom{\protect\Rom@P@}
{\catcode`\-11 \HyperRefs{idx} \HyperRefs{glo}
\newcount\cnt@idx \global\cnt@idx=10000
\newcount\cnt@glo \global\cnt@glo=10000
\gdef\writeindex#1{\W@X@{\cntref@{idx}}\tf@-idx
 {\string\indexentry{#1}{\Hlast@{idx}}{\thepage}}}
\gdef\writeglossary#1{\W@X@{\cntref@{glo}}\tf@-glo
 {\string\glossaryentry{#1}{\Hlast@{glo}}{\thepage}}}
}
\def\emph#1{\@itcorr@\bgroup\em\ignorespaces#1\unskip\egroup
  \DN@{\DN@{}\ifx\next.\else\ifx\next,\else\DN@{\/}\fi\fi\next@}\FN@\next@}
\def\makequoteactive{\catcode`\"\active}
{\makequoteactive\gdef"{\FN@\quote@}
\gdef\quote@{\ifx"\next\DN@"##1""{\quoteii{##1}}\else\DN@##1"{\quotei{##1}}\fi\next@}}
\let\quotei\eat@
\let\quoteii\eat@
\def\MakeIndex{\@openout{idx}}
\def\MakeGlossary{\@openout{glo}}

\def\endofpar#1{\ifmmode\ifinner\endofpar@{#1}\else\eqno{#1}\fi
    \else\leavevmode\endofpar@{#1}\fi}
\def\endofpar@#1{\unskip\penalty\z@\null\hfil\hbox{#1}\hfilneg\penalty\@M}

\newdimen\normalparindent\normalparindent\parindent
\def\firstparindent#1{\everypar\expandafter{\the\everypar
  \global\parindent\normalparindent\global\everypar{}}\parindent#1\relax}

\@@addto\disablepreambule@cs{%
    \\\readaux\relax
    \\\begin\relax
    \\\readaux@\relax
    \\\@openout\eat@
    \\\@@openout\eat@
    \/\Monograph\empty
    \/\MakeIndex\empty
    \/\MakeGlossary\empty
    \/\MakeToc\eat@
    \/\HyperRefs\eat@
    \/\NoHyperRefs\eat@
}

\csname label.def\endcsname


\def\punct#1#2{\if\notempty{#2}#1\fi}
\def\sppunct{\punct{.\enspace}}
\def\varpunct#1#2{\if\frillsnotempty{#2}#1\fi}

\def\headstyle#1#2{\numberline{#1\sppunct{#2}}\ignorespaces#2\unskip}
\def\headtocstyle#1#2{\numberline{#1\punct.{#2}}\space #2}

\def\specialtocstyle#1#2{#2}
\newcounter\section\null
\newcounter\subsection\section
\newcounter\subsubsection\subsection
\newhead\specialsection[\specialtocstyle]\null\endspecialhead
\newhead\section\section\endhead
\newhead\subsection\subsection\endsubhead
\newhead\subsubsection\subsubsection\endsubsubhead
\def\firstappendix{\global\sectionno0 %
  \counterstyle\section{\Alphnum\sectionno}%
    \global\let\firstappendix\empty}

\def\appendixtocstyle#1#2{\space\numberline{Appendix #1\sppunct{#2}}#2}
\newhead\appendix[\appendixtocstyle]\section\endhead

\let\endAmSdef\enddefinition
\def\proclaimstyle#1#2{\numberline{#2\varpunct{.\enspace}{#1}}\frills{#1}}
\copycounter\thm\subsubsection
\theorem\thm\endproclaim
\proposition\thm\endproclaim
\lemma\thm\endproclaim
\corollary\thm\endproclaim
\definition\thm\endAmSdef
\example\thm\endAmSdef

\def\captionstyle#1#2{\frills{#1}\numberline{\varpunct{ }{#1}#2}}
\newcounter\figure\null
\newcounter\table\null
\newcaption{Figure}\figure\figure{lof}\botcaption
\newcaption{Table}\table\table{lot}\topcaption

\copycounter\equation\subsubsection



\def\stydate{November 27, 1997}
\immediate\write16{This is DEGT.DEF by A.Degtyarev <\stydate>}
{\edef\temp{\the\everyjob\immediate\write16{DEGT.DEF: <\stydate>}}
\global\everyjob\expandafter{\temp}}

\chardef\tempcat\catcode`\@\catcode`\@=11

\let\ge\geqslant
\let\le\leqslant
\def\C{{\Bbb C}}
\def\R{{\Bbb R}}
\def\Z{{\Bbb Z}}
\def\Q{{\Bbb Q}}
\def\N{{\Bbb N}}

\let\ZZ\Z

\def\Cp#1{\C{\operator@font p}^{#1}}
\def\Rp#1{\R{\operator@font p}^{#1}}
\def\Hom{\qopname@{Hom}}
\def\Ext{\qopname@{Ext}}
\def\Tors{\qopname@{Tors}}

\def\Im{\qopname@{Im}}          
\def\Re{\qopname@{Re}}     %
\def\Ker{\qopname@{Ker}}
\def\Coker{\qopname@{Coker}}
\def\Int{\qopname@{Int}}
\def\Cl{\qopname@{Cl}}
\def\Fr{\qopname@{Fr}}
\def\Fix{\qopname@{Fix}}
\def\tr{\qopname@{tr}}
\def\inj{\qopname@{in}}
\def\id{\qopname@{id}}
\def\pr{\qopname@{pr}}
\def\rel{\qopname@{rel}}
\def\pt{{\operator@font{pt}}}
\def\const{{\operator@font{const}}}
\def\codim{\qopname@{codim}}
\def\cdim{\qopname@{dim_{\C}}}
\def\rdim{\qopname@{dim_{\R}}}
\def\conj{\qopname@{conj}}
\def\rank{\qopname@{rk}}
\def\sign{\qopname@{sign}}
\def\gcd{\qopname@{g.c.d.}}
\let\sminus\smallsetminus
\def\set<#1|#2>{\bigl\{#1\bigm|#2\bigr\}}


\def\preprint#1{\hrule height0pt depth0pt\kern-24pt%
  \hbox to\hsize{#1}\kern24pt}
\def\today{\ifcase\month\or January\or February\or March\or
  April\or May\or June\or July\or August\or September\or October\or
  November\or December\fi \space\number\day, \number\year}

\def\n@te#1#2{\leavevmode\vadjust{%
 {\setbox\z@\hbox to\z@{\strut\eightpoint#1}%
  \setbox\z@\hbox{\raise\dp\strutbox\box\z@}\ht\z@=\z@\dp\z@=\z@%
  #2\box\z@}}}
\def\leftnote#1{\n@te{\hss#1\quad}{}}
\def\rightnote#1{\n@te{\quad\kern-\leftskip#1\hss}{\moveright\hsize}}
\def\?{\FN@\qumark}
\def\qumark{\ifx\next"\DN@"##1"{\leftnote{\rm##1}}\else
 \DN@{\leftnote{\rm??}}\fi{\rm??}\next@}

\def\centerpage{\dimen@=6.5truein \advance\dimen@-\hsize\hoffset.5\dimen@}
\ifnum\mag>1000 \centerpage\fi

\def\nologo{\let\logo@\relax}

\expandafter\ifx\csname eat@\endcsname\relax\def\eat@#1{}\fi
\expandafter\ifx\csname operator@font\endcsname\relax
 \def\operator@font{\roman}\fi
\expandafter\ifx\csname eightpoint\endcsname\relax
 \let\eightpoint\small\fi

\catcode`\@\tempcat\let\tempcat\undefined

\def\stydate{September 17, 1998}
\immediate\write16{This is CD.DEF by A.Degtyarev <\stydate>}
{\edef\temp{\the\everyjob\immediate\write16{CD.DEF: <\stydate>}}
\global\everyjob\expandafter{\temp}}

\expandafter\edef\csname cd.def\endcsname{%
\catcode`\noexpand\@\the\catcode`\@\catcode`\noexpand\&\the\catcode`\&
  \edef\expandafter\noexpand\csname cd.def\endcsname{\stydate}}
\catcode`\@=11 \catcode`\&=11


%
\newcount\mindasharrowwidth \mindasharrowwidth=3
\newif\if&CDdash
\newif\if&CDdouble
\def\&CDlefttrue#1{\let\if&CDleft\iftrue\csname &CD\string#1\endcsname}
\def\&CDleftfalse#1{\let\if&CDleft\iffalse\csname &CD\string#1\endcsname}

\def\toz@{to\z@}

\expandafter\ifx\csname msa@group\endcsname\relax 
    \expandafter\ifx\csname msafam\endcsname\relax\let\next\relax 
  \else
    \edef\next{\hexnumber@\msafam}
    \font@\eightmsa=msam8
    \font@\sixmsa=msam6
    \addto\eightpoint{\textfont\msafam\eightmsa \scriptfont\msafam\sixmsa}
    \addto\tenpoint{\textfont\msafam\tenmsa \scriptfont\msafam\sevenmsa}
  \fi\else
    \edef\next{\ifcase\msa@group0 \or1 \or2 \or3 \or4 \or5 \or6 \or7
        \or8 \or9 \or A \or B \or C \or D \or F \fi}
\fi \ifx\next\relax\let\&CDdash\relax\else
  \mathchardef\&CDrdash"3\next4B
  \mathchardef\&CDldash"3\next4C
  \def\&CDdash{{\ifnum\mindasharrowwidth=0 \else 
    \setboxz@h{$\m@th\dabar@$}\dimen@\wdz@\multiply\dimen@\mindasharrowwidth
    \loop\ifdim\dimen@<\bigaw@\advance\dimen@\wdz@\repeat
    \global\bigaw@\dimen@\fi}}
\fi

\def\&CDbar{\if&CDdash\mathrel\dabar@
  \else\if&CDdouble\Relbar\else\relbar\fi\fi}
\def\&CDleft{\if&CDdouble\Leftarrow\else\leftarrow\fi}
\def\&CDright{\if&CDdouble\Rightarrow\else\rightarrow\fi}
\def\&CDfill#1{\mkern-2mu\mathord#1\mkern-2mu}
\def\&CDlap#1{\mathrel{\hbox\toz@{\hss$\m@th#1$\hss}}}

{\def\&CDdef#1{\expandafter\gdef\csname &CD\string#1\endcsname}
\expandafter\let\csname &CD-\endcsname\relbar \expandafter\let\csname
&CD=\endcsname\Relbar \&CDdef<{\relax\if&CDdash\&CDldash\else\&CDleft\fi}
\&CDdef>{\relax\if&CDdash\&CDrdash\else\&CDright\fi}
\&CDdef\<{\relax\if&CDdash\&CDlap{\&CDldash\joinrel}\&CDldash
  \else\&CDleft\mkern-14mu\&CDleft\fi}
\&CDdef\>{\relax\if&CDdash\&CDrdash\&CDlap{\joinrel\&CDrdash}
  \else\&CDright\mkern-14mu\&CDright\fi}
\&CDdef({\&CDlap{\lhook\joinrel}\&CDbar}
\&CDdef){\&CDbar\&CDlap{\mskip-1.5mu\rhook}}
\&CDdef|{\relax\if&CDleft\mapstochar\&CDbar
  \else\&CDbar\mskip-2.5mu\mapstochar\mskip 2.5mu\fi}
\&CDdef{&-}#1#2-#3-#4{\&CDdoublefalse\&CDline#1{\&CDfill-}#4{#2}{#3}}
\&CDdef{&=}#1#2=#3=#4{\&CDdoubletrue\&CDline#1{\&CDfill=}#4{#2}{#3}}
\&CDdef{&:}#1#2:#3:#4{\&CDdashtrue\&CDline#1\dabar@#4{#2}{#3}}}

\def\&CDline#1#2#3#4#5{%
\ifCD@$\global\bigaw@\minCDaw@
    \else\global\bigaw@\minaw@\mathrel{\fi\&CDls
  \setboxz@h{$\m@th\ssize\;{#4}\;\;$}%
  \setbox@ne\hbox{$\m@th\ssize\;{#5}\;\;$}%
  \setbox\tw@\hbox{$\m@th#5$}%
  \ifdim\wdz@>\bigaw@\global\bigaw@\wdz@\fi
  \ifdim\wd@ne>\bigaw@\global\bigaw@\wd@ne\fi\&CDhs
  \ifCD@\enskip\&CDdash\else\if&CDdash\&CDdash\fi\fi\hskip.5\bigaw@ plus.5fil%
  \hbox\toz@{\hss$\m@th\mathop{\vphantom=}\limits^{#4}%
    \ifdim\wd\tw@>\z@ _{#5}\fi$\hss}%
  \hskip-.5\bigaw@ plus-.5fil\rlap{$\m@th\&CDlefttrue#1$}%
  \if&CDdash\else\global\advance\bigaw@-6\ex@\kern 3\ex@\fi
  \cleaders\hbox{$\m@th#2$}\hskip\bigaw@ plus1fil%
  \if&CDdash\else\kern 3\ex@\fi\llap{$\m@th\&CDleftfalse#3$}%
  \ifCD@\enskip\&CDrs$\else}\fi\&CDdashfalse\ampersand@}
\def\&CDhline#1#2{\ampersand@\ifCD@\multispan\mscount\fi
  \csname &CD&#2\endcsname#1}

\def\CDhline#1{\ifCD@\let\<\relax\fi
  \afterassignment\&CDhline\global\mscount0#1}

\atdef@:{\relax\ifmmode\CDhline\else\leavevmode\null:\fi}
\atdef@>#1>#2>{\CDhline--{#1}-{#2}->}
\atdef@<#1<#2<{\CDhline<-{#1}-{#2}--} \atdef@={\CDhline=====}
\atdef@[#1>#2>{\CDhline(-{#1}-{#2}->}
\atdef@]#1<#2<{\CDhline<-{#1}-{#2}-)}

\def\hexpand(#1,#2){\global\def\&CDls{\dimen@#1\global\let\&CDls\relax
  \ifCD@\kern-\dimen@\def\&CDrs{\dimen@#2\kern-\dimen@}%
    \def\&CDhs{\global\bigaw@\z@\mindasharrowwidth\z@}%
  \else\global\advance\bigaw@\dimen@\advance\minaw@#2\fi}}
\let\&CDls\empty
\let\&CDrs\empty
\let\&CDhs\empty

%
\let\&CDpref\Big
\def\short{\let\&CDpref\big}

\def\&CD&&{\ifx\next$\let\next@\relax
  \else\DN@{\ampersand@\ampersand@}\fi\next@}

\def\&CDmbox#1{\hbox to\wdz@{\hss$\m@th#1$\hss}}
\let\&CDover\empty
\let\&CDunder\empty
\def\&CDvert{\if&CDdouble\mathchar"377\else\mathchar"33F\fi}
\def\&CDuparrow{\if&CDdouble\mathchar"37E\else\mathchar"378\fi}
\def\&CDdownarrow{\if&CDdouble\mathchar"37F\else\mathchar"379\fi}
\let\&CDup\&CDvert
\let\&CDdown\&CDvert
\def\&CD@up{\let\&CDup\&CDuparrow}
\def\&CD@down{\let\&CDdown\&CDdownarrow}
\def\&CD@Up{\def\&CDup{\setboxz@h{$\m@th\&CDuparrow$}%
    \vbox{\copy\z@\kern-3\ex@\box\z@}}}
\def\&CD@Down{\def\&CDdown{\setboxz@h{$\m@th\&CDdownarrow$}%
    \vbox{\copy\z@\kern-3\ex@\box\z@}}}

\def\&CDvline{\ifx\next<\let\next@\CDvline\else\DN@{\CDvline<||>}\fi\next@}
\def\CDvline<#1|#2|#3>{\setboxz@h{$\&CDpref\uparrow$}\&CDuspan\&CDdspan
  \advance\dimen@ii\dp\z@\advance\dimen@\ht\z@\advance\dimen@\dimen@ii
  \llap{$\m@th\vcenter{\hbox{$\ssize#1$}}$}%
  \vrule height\ht\z@ depth\dp\z@ width\z@
  \setboxz@h{\vbox to\dimen@{\offinterlineskip
    \let\up\&CD@up\let\down\&CD@down\let\double\&CDdoubletrue
    \let\Up\&CD@Up\let\Down\&CD@Down\let\dash\&CDdashtrue
    #2\setboxz@h{$\m@th\&CDuparrow$}\&CDmbox{\smash\&CDover}%
    \if&CDdash\&CDvdash\else
      \&CDmbox\&CDup\kern-3\ex@\cleaders\&CDmbox\&CDvert\vfil
      \kern-3\ex@\&CDmbox\&CDdown\fi
    \&CDmbox{\smash\&CDunder}}}%
  \dp\z@-\dimen@ii\ht\z@\dimen@ii\lower\dimen@ii\boxz@
  \rlap{$\m@th\vcenter{\hbox{$\ssize#3$}}$}\FN@\&CD&&}
\def\&CDvdash{\def\k@rn{.5\ex@}%
  \def\b@x##1{\&CDmbox{\vbox to##1{\ifdim##1<2\k@rn\vss
      \else\kern\k@rn\leaders\vrule\vfil\kern\k@rn\fi}}}%
  \ifx\&CDup\&CDvert\else\setbox\tw@\&CDmbox\&CDup\advance\dimen@-\dp\tw@\fi
  \ifx\&CDdown\&CDvert\else\setbox\tw@\&CDmbox\&CDdown\advance\dimen@-\dp\tw@\fi
  \ifdim\dimen@>\z@\loop\advance\dimen@-\dp\z@\ifdim\dimen@>\z@\repeat\fi
  \ifdim\dimen@<\z@\advance\dimen@\dp\z@\fi
  \ifx\&CDup\&CDvert\kern-\k@rn\dimen@ii\z@
    \else\&CDmbox\&CDup\kern\k@rn\dimen@ii\dimen@\dimen@\z@\fi
  \b@x\dimen@\cleaders\b@x{\dp\z@}\vfil\b@x\dimen@ii
  \ifx\&CDdown\&CDvert\kern-\k@rn\else\kern\k@rn\&CDmbox\&CDdown
\fi}

\atdef@+{\FN@\&CDvline} \atdef@ V#1V#2V{\CDvline<#1|\&CD@down|#2>}
\atdef@ A#1A#2A{\CDvline<#1|\&CD@up|#2>}
\atdef@|{\&CDdoubletrue\FN@\&CDvline}
\atdef@\vert{\&CDdoubletrue\FN@\&CDvline}

\def\upmapsto{\def\&CDover{\vrule height.4\ex@ width2.4\ex@}}
\def\downmapsto{\def\&CDunder{\vrule depth.4\ex@ width2.4\ex@}}

\def\&CDhook#1#2#3#4{\def#1{\raise#2\ex@\hbox{%
    $\m@th\scriptscriptstyle\mkern#3mu#4$}}}

\def\&CDuspan{\dimen@\z@}
\def\&CDdspan{\dimen@ii\z@}
\def\vexpand(#1,#2){\ifCD@\dimen@#1\dimen@ii#2%
  \let\&CDuspan\relax\let\&CDdspan\relax\fi}
\def\expand{\ifCD@\dimen@ii2.5\ex@\let\&CDdspan\relax\fi}
\def\Expand{\ifCD@\dimen@ii21.5\ex@\let\&CDdspan\relax\fi}
\def\expandup{\ifCD@\dimen@2.5\ex@\let\&CDuspan\relax\fi}
\def\Expandup{\ifCD@\dimen@21.5\ex@\let\&CDuspan\relax\fi}

%
\expandafter\ifx\csname minCDaw@\endcsname\relax
 \let\minCDaw@\minCDarrowwidth\fi
\expandafter\ifx\csname ssize\endcsname\relax\let\ssize\scriptstyle\fi

\csname cd.def\endcsname

\def\stydate{September 26, 1998}
\immediate\write16{This is DD.DEF by A.Degtyarev <\stydate>}
{\edef\temp{\the\everyjob\immediate\write16{DD.DEF: <\stydate>}}
\global\everyjob\expandafter{\temp}}

\expandafter\edef\csname dd.def\endcsname{%
\catcode`\noexpand\@\the\catcode`\@\catcode`\noexpand\&\the\catcode`\&
  \edef\expandafter\noexpand\csname dd.def\endcsname{\stydate}}
\catcode`\@=11 \catcode`\&=11

\font\lcircle=lcircle10 \chardef\&DDcirc'143 \chardef\&DDbullet'163
\chardef\&DDCirc'150

\newcount\DDlinelength \DDlinelength=3
\newdimen\&DDu
\newdimen\&DDdp
\newcount\&DDln
\def\toz@{to\z@}

\def\&DDhbox#1#2{\setboxz@h{$#1$}\ifdim\wdz@>\z@\lower.5\dimen@ii
  \vbox{\hbox\toz@{$\m@th\ssize#2$}}\fi}
\def\&DDboxz#1{\hbox\toz@{\hss#1\hss}}
\def\&DDmbox#1{\&DDboxz{$\m@th\ssize#1$}}
\def\&DDscript#1#2|#3><#4|#5>{%
  \setboxz@h\toz@{\setboxz@h{$\ssize0$}\dimen@ii\ht\z@
  \def\&k{2.5\p@}\dimen@.5\dimen@i
  \hss\&DDhbox{#2}{\hss#2\kern\&k}\kern\dimen@
  {\baselineskip\&k\lineskip.5\p@\lineskiplimit\lineskip
    \setboxz@h{$#4$}\vbox{\ifdim\wdz@>\z@\&DDmbox{#4}\hbox{}\fi
      \hrule height\dimen@ width\z@}%
    \setboxz@h{$#5$}\vtop{\hrule height\z@ depth\dimen@ width\z@
        \ifdim\wdz@>\z@\advance\baselineskip\dimen@ii\null\&DDmbox{#5}\fi}}%
  \kern\dimen@\&DDhbox{#3}{\kern\&k#3\hss}\hss}%
  \if&DDfl\else\ht\z@\z@\fi\ifdim\dp\z@>\&DDdp\global\&DDdp\dp\z@\fi
  \dp\z@\z@\kern-#1\boxz@\kern#1}

\def\&DDparm#1{\ifx\&DDparm#1\else\setboxz@h{\lcircle#1}%
    \ht\z@\z@\dp\z@\z@\copy\z@\kern-\wdz@\ifdim\wdz@>\dimen@\dimen@\wdz@\fi
    \expandafter\&DDparm\fi}
\def\&DDnode#1){\vrule height\&DDu depth\&DDu width\z@\dimen@\z@
  \hbox to\tw@\&DDu{\let\circ\&DDcirc\let\Circ\&DDCirc\let\bullet\&DDbullet
  \def\.{\&DD&text.}\def\*{\&DD&math*}\hss\&DDparm#1\&DDparm\hss
  \global\dimen@i\dimen@}\atdef@*<{\&DDscript\&DDu}}
\def\&DD&text#1{\vbox\toz@{\vss\&DDboxz{\rm#1}\vss}}
\def\&DD&math#1{\lower\fontdimen22\scriptfont\tw@
    \vbox{\&DDmbox{\mathord#1}}}

\def\&DDskip{\ifnum\&DDln<\z@\&DDln\DDlinelength\fi
  \dimen@\&DDln\&DDu\skip@\tw@\dimen@ plus1fil}
\def\&DDreset{\global\mscount\z@\global\&DDln\m@ne\global\let\if&DDdots\iffalse}

\def\&DDhr#1#2{\cleaders\hrule height#1\p@ depth#2\p@\hskip\skip@}
{\catcode`\^=4 \gdef\&DDhor#1#2#3{^\multispan\mscount\&DDskip
  \if&DDdots{\divide\skip@\tw@\advance\skip@-4\&DDu#1%
    \vbox\toz@{\vss\hbox to8\&DDu{\hss$\,\ldots$\hss}\vss}#1}%
  \else#1\fi
  \hskip-\dimen@ plus-.5fil{\global\dimen@i\p@\&DDscript\z@|><#2|#3>}%
  \hskip\dimen@ plus.5fil%
  \&DDreset^}}

\def\&DDvr#1{\vrule width#1\p@ height\skip@}
{\catcode`\^=4
\gdef\&DDvert#1#2#3{\&DDskip\vbox to\skip@{%
  \ifnum\mscount=\z@\else\dimen@\mscount\&DDu\skip@\tw@\dimen@ plus1fil\fi
  \vskip\dimen@{\global\dimen@i\p@\&DDscript\z@#2|#3><|>}\vskip-\dimen@
  \if&DDdots{\divide\skip@\tw@\advance\skip@-4\&DDu\&DDmbox{#1}%
    \vbox to8\&DDu{\vss\&DDmbox{\vdots}\kern6\p@\vss}\&DDmbox{#1}}%
  \else\&DDmbox{#1}\fi\vss}\&DDreset^^}}
{\def\&DDdef#1{\expandafter\gdef\csname &DD&\string#1\endcsname}
\&DDdef-#1-#2-{\&DDvert{\&DDvr{.4}}{#1}{#2}}
\&DDdef=#1=#2={\&DDvert{\&DDvr{.4}\kern1.2\p@\&DDvr{.4}}{#1}{#2}}
\&DDdef\-#1-#2-{\&DDvert{\&DDvr{1}}{#1}{#2}}
\&DDdef:#1:#2:{\&DDvert{\vbox{\cleaders\vbox{\advance\&DDu-.5\p@\kern.5\p@
    \hrule height\tw@\&DDu width.4\p@\kern.5\p@}\vskip\skip@}}{#1}{#2}}
\&DDdef.{\&DDvert{\&DDvr{0}}\empty\empty}}

{\catcode`\^=4
\gdef\&DDtext#1#2[#3]{^\multispan\mscount\setboxz@h{%
  \lower\fontdimen22\textfont\tw@\vbox{\hbox{#3}}}%
  \if&DDfl\else\ht\z@\z@\fi\ifdim\dp\z@>\&DDdp\global\&DDdp\dp\z@\fi
  \dp\z@\z@#1\boxz@#2\&DDreset}}

{\catcode`\^=4 \gdef\&DD{\bgroup
  \everycr{\noalign{\global\&DDdp\z@\global\let\if&DDfl\iffalse\&DDreset}}%
  \offinterlineskip\tabskip\z@skip\global\let\if&DDfl\iffalse\&DDreset
  \setboxz@h{\lcircle\&DDcirc}\&DDu.5\wd\z@
  \def\CR{^\vrule height\z@ depth\&DDdp width\z@\cr
    \noalign{\global\let\if&DDfl\iftrue}}%
  \atdef@({\&DDnode}
  \atdef@-##1-##2-{\&DDhor{\&DDhr{.2}{.2}}{##1}{##2}}%
  \atdef@=##1=##2={\&DDhor{\&DDhr{1}{-.6}\hskip-\skip@\&DDhr{-.6}{1}}{##1}{##2}}%
  \atdef@\-##1-##2-{\&DDhor{\&DDhr{.5}{.5}}{##1}{##2}}%
  \atdef@:##1:##2:{\&DDhor{\cleaders\hbox{\advance\&DDu-.5\p@\kern.5\p@
    \vrule height.2\p@ depth.2\p@ width\tw@\&DDu\kern.5\p@}\hskip\skip@}{##1}{##2}}%
  \atdef@.{\&DDhor{\&DDhr00}\empty\empty}%
  \atdef@|##1{\csname &DD&\string##1\endcsname}%
  \atdef@[{\&DDtext\hss\hss[}%
  \atdef@\l{\&DDtext\empty\hss}\atdef@\r{\&DDtext\hss\empty}%
  \def\expand{\global\mscount}%
  \def\withdots{\global\let\if&DDdots\iftrue}%
  \halign\bgroup^\hss$\m@th##$\hss\cr}}
\def\endDD{\crcr\egroup\egroup}

\def\DDtop{\vtop\&DD}

\csname dd.def\endcsname

\newdimen\postdemoskip
\postdemoskip\medskipamount


{\catcode`\@=11 \gdef\proclaimfont@{\sl}}

\ifx\=\undefined\let\=\B\fi

\def\Dg#1:{\par\noindent{\bf Dg #1}:\enspace\ignorespaces}
\def\DgIt:{\par\noindent{\bf Dg and It}:\enspace\ignorespaces}

\let\textdot\D

\def\ie{\emph{i.e.}}
\def\eg{\emph{e.g.}}
\def\cf{\emph{cf.}}

\let\Ga\alpha
\let\Gb\beta
\let\Ge\varepsilon
\let\Gf\varphi
\let\Gg\gamma
\let\Go\omega
\let\Gr\rho
\let\Gs\sigma
\let\Gl\lambda
\let\GL\Lambda
\let\GO\Omega
\let\GG\Gamma
\let\GS\Sigma

\let\GD\Delta
\let\Gk\kappa
\let\Gt\theta
\let\Gi\iota

\let\wb\bar
\let\wt\tilde

\let\<\langle
\let\>\rangle

\def\sb#1{_{\scriptscriptstyle#1}}
\def\sp#1{^{\scriptscriptstyle#1}}

\def\sbp{\sb+}
\def\sbm{\sb-}
\def\sppm{\sp\pm}
\def\spp{\sp+}
\def\spm{\sp-}

\def\ZZ{\Z_2}
\def\ls|#1|{\mathopen|#1\mathclose|}
\def\Cp#1{\Bbb P^{#1}}

\def\Per{\let\Gkern\!\operatorname{Per}}
\def\per{\operatorname{per}}
\def\Aut{\operatorname{Aut}}

\def\Gal{\operatorname{Gal}}
\def\Pic{\operatorname{Pic}}
\def\ord{\operatorname{ord}}
\def\discr{\operatorname{discr}}
\def\Gramm{\operatorname{Gr}}

\def\group#1{\text{\sl#1\/}}
\def\grO{\group{O}}
\def\grSO{\group{SO}}
\def\grGL{\group{GL}}
\def\grU{\group{U}}
\def\grSL{\group{SL}}

\def\O{\Cal O}

\def\eL{L^\bullet}

\def\ort{^{\scriptscriptstyle\perp}}

\def\HS{\Cal H}
\def\dHS{\partial\HS}
\def\HC{\Cal C}
\def\hp#1{\frak h_{#1}}

\def\hps#1{\frak s_{#1}}
\def\sh{\operatorname{h}}
\def\psh{\operatorname{h}}
\def\CT{\Cal T}
\def\CJ{\Cal J}
\def\Gsbp{G^0}
\def\MS{\frak M}
\def\KMS{K\MS}
\def\ms{\frak m}
\def\kms{\ms^K}
\def\CQ{\Cal Q}
\def\D{\Bbb D}

\def\K{\Bbb K}
\let\k\Bbbk
\def\P{\Bbb P}
\def\gw{\frak w}
\def\gv{\frak v}
\def\gD{\frak O}
\def\gg{\frak g}

\let\Gkern\relax
\def\G{\Gkern\let\Gkern\relax{}^G}
\def\prhull{\sphat\mkern3mu}
\let\Jper\per

\def\anti-{\hbox{(anti-)}\hskip0pt}

{\catcode`\@=11 \gdef\enddemo{\par\revert@envir\enddemo\endremark}}
\let\oldremark\endremark
\def\endremark{\oldremark\vskip\postdemoskip}

\newhead\subsubsection\subsubsection\endsubhead
{\catcode`\@=11 \gdef\subsubheadfont@{\bf}}

\MakeToc{toc}

\topmatter

\title
Finiteness and quasi-simplicity for symmetric $K3$-surfaces
\endtitle

\author
Alex Degtyarev, Ilia Itenberg, and Viatcheslav Kharlamov
\endauthor

\thanks
Second author is a member of Research Training Network RAAG, supported by
the European Human Potential Program.
\endthanks

\thanks
Third author is a member of Research Training Networks EDGE and RAAG,
supported by the European Human Potential Program.
\endthanks

\address
Bilkent University\endgraf\nobreak
06533 Ankara, Turkey
\endaddress

\email
degt\@fen.bilkent.edu.tr
\endemail

\address
CNRS \endgraf
Institut de Recherche Math\'ematique de Rennes
\endgraf
35042 Rennes Cedex, France
\endaddress

\email
itenberg\@maths.univ-rennes1.fr
\endemail

\address
Universit\'e Louis Pasteur et IRMA (CNRS)\endgraf
7 rue Ren\'e Descartes 67048 Strasbourg Cedex, France
\endaddress

\email
kharlam\@math.u-strasbg.fr
\endemail

\abstract We compare the smooth and deformation equivalence of actions of
finite groups on $K3$-surfaces by holomorphic and anti-holomorphic
transformations. We prove that the number of deformation classes is
finite and, in a number of cases, establish the expected coincidence of
the two equivalence relations. More precisely, in these cases we show
that an action is determined by the induced action in the homology. On
the other hand, we construct two examples to show that, first, in general
the homological type of an action does not even determine its topological
type, and second, that $K3$-surfaces $X$ and $\bar X$ with the same Klein
action do not need to be equivariantly deformation equivalent even if the
induced action on $H^{2,0}(X)$ is real, \ie, reduces to multiplication by
$\pm 1$.
\endabstract

\endtopmatter

\document

\section{Introduction}

\subsection{Questions}\label{questions}
In this paper, we study equivariant deformations of complex $K3$-surfaces
with symmetry groups, where by a symmetry we mean an either holomorphic
or anti-holomorphic transformation of the surface. Although the
automorphism group of a particular $K3$-surface may be infinite, we
confine ourselves to finite group actions and address the following two
questions (see~\ref{terminology}--\ref{deformations} for precise
definitions):
\def\iiitem#1{\leavevmode\hbox{\kern-\rosteritemwd\kern-.5em #1\enspace}\ignorespaces}
\roster
\item""
\iiitem{\emph{finiteness}:} whether the number of actions, counted up to
equivariant deformation and isomorphism, is finite, and
\item""
\iiitem{\emph{quasi-simplicity}:} whether the differential topology of an
action determines it up to the above equivalence.
\endroster
The response to the second question, in the way that it is posed, is
obviously in the negative. For example, given an action on a surface~$X$,
the same action on the complex conjugate surface~$\bar X$ is
diffeomorphic to the original one but often not deformation equivalent to
it. Thus, we pose this question in a somewhat weaker form:
\roster
\item""
\iiitem{\emph{weak quasi-simplicity}:} does the differential topology of
an action determine it up to equivariant deformation and
\anti-isomorphism?
\endroster
Up to our knowledge, these questions have never been posed explicitly,
and, moreover, despite numerous related partial results, they
both remained open.

One may notice a certain ambiguity in the statement of the above
questions, especially in what concerns the quasi-simplicity: we do not
specify whether we consider diffeomorphic actions on true $K3$-surfaces
or, more generally, diffeomorphic actions on surfaces diffeomorphic to a
$K3$-surface. Fortunately, a surface diffeomorphic to a $K3$-surface is a
$K3$-surface, see~\cite{FMbook}, and the two versions turn out to be
equivalent. Thus, we confine ourselves to true $K3$-surfaces and respond
to both the finiteness and (to great extend) weak quasi-simplicity
questions (see~\ref{results}).

Following the founding work by I.~Piatetski-Shapiro and I.~Shafarevich
\cite{PSh-Sh}, we base our study on the global Torelli theorem. When
combined with Vik.~Kulikov's theorem on surjectivity of the period map
\cite{Kulikov}, this fundamental result essentially reduces the
finiteness and quasi-simplicity questions to certain arithmetic problems.
It is this approach that was used by V.~Nikulin in~\cite{N1}
and~\cite{N2}, where he established (partially implicitly) the finiteness
and quasi-simplicity results for polarized $K3$-surfaces with symplectic
actions of finite abelian groups and for those with real structures.
In~\cite{DIK}, we extended these results to real Enriques surfaces. (Note
that a real Enriques surface can be regarded as a $K3$-surface with a
certain action of $\ZZ\times\ZZ$. In~\cite{DIK} we give, in fact, the
full deformation classification of such actions.) While studying real
Enriques surfaces, we got interested in the above questions and obtained
our first results in this direction.

\subsection{Related results}\label{related-results}
One can find a certain similarity between our finiteness results and the
finiteness in theory of moduli of complex structures on $4$-manifolds,
which states (see~\cite{FM} and \cite{F}) that the moduli space of
K\"ahlerian complex structures on a given underlying differentiable compact
$4$-manifold has finitely many components. (By K\"ahlerian we mean a
complex structure admitting a K\"ahler metric. In the case of surfaces this
is a purely topological restriction: the complex structures on a given
compact $4$-manifold~$X$ are K\"ahlerian if and only if the first Betti
number $b_1(X;\Q)$ is even.) Moreover, the moduli space is connected as
soon as there is a K\"ahlerian representative of Kodaira dimension $\le0$
(as it is the case for $K3$-surfaces and complex $2$-tori); for Kodaira
dimension one, there are at most two deformation classes, which are
represented by~$X$ and~$\bar X$, see~\cite{FM}. Examples of general type
surfaces $X$ not deformation equivalent to $\bar X$ are found in
\cite{KK} and \cite{Cat}.

The principle results of our paper can be regarded as an equivariant
version of the above statements for $K3$-surfaces. The finiteness
theorem~\ref{Main-Finiteness} is closely related to a series of results
from theory of algebraic groups that go back to C.~Jordan~\cite{Jordan}.
The original Jordan theorem states that $\grSL(n,\Z)$ contains but a
finite number of conjugacy classes of finite subgroups. A.~Borel and
Harish-Chandra, see \cite{BH} and \cite{B}, generalized this statement to
any arithmetic subgroup of an algebraic group; further recent
generalizations are due to V.~Platonov, see~\cite{Pl}. Note that,
together with the global Torelli theorem, these Jordan type theorems
(applied to the $2$-cohomology lattice of a $K3$-surface) imply that the
number of different finite groups acting faithfully on $K3$-surfaces is
finite. A complete classification of finite groups acting symplectically
(\ie, identically on holomorphic forms) on $K3$-surfaces is found in
Sh.~Mukai~\cite{Mukai} (see also Sh.~Kond\=o~\cite{Kondo} and
G.~Xiao~\cite{Xiao}; the abelian groups where first classified by
Nikulin~\cite{N2}). A sharp bound on the order of a group acting
holomorphically on a $K3$-surface is given by Kond\=o~\cite{K2}.

Among other related finiteness results found in the literature, we would
like to mention  a theorem by Piatetski-Shapiro and
Shafarevich~\cite{PSh-Sh} stating that the automorphism group of an
algebraic $K3$-surface is finitely generated, our~\cite{DIK}
generalization of this theorem to all $K3$-surfaces, and
H.~Sterk's~\cite{St} finiteness results on the classes of irreducible
curves on an algebraic $K3$-surface.
Note that all these results deal with individual surfaces
rather than with their deformation classes. They are related to the
finiteness of the number of conjugacy classes of finite subgroups in the
group of Klein automorphisms of a given variety. As a special case, one
can ask whether the number of conjugacy classes of real structures on a
given variety is finite. For the latter question, the key tool is the
Borel-Serre~\cite{BSerre} finiteness theorem for Galois cohomology of
finite groups; as an immediate consequence, it implies finiteness of the
number of conjugacy classes of real structures on an abelian variety.
In~\cite{DIK} we extended this statement to all surfaces of Kodaira
dimension~$\ge1$ and to all minimal K\"ahler surfaces. Remarkably,
finiteness of the number of conjugacy classes of real structures on a
given rational surface is still an open question.

Unlike finiteness, the quasi-simplicity question does not make much sense
for individual varieties. In the past, it was mainly studied for
deformation equivalence of real structures: given a deformation family of
complex varieties, is a real variety within this family determined up to
equivariant deformation by the topology of the real structure? The first
non trivial result in this direction, concerning real cubic surfaces
in~$\Cp3$, was discovered by F.~Klein and L.~Schl\"afli (see, \eg, the
survey~\cite{survey}). At present, the answer is known for curves
(essentially due to F.~Klein and G.~Weichold, see, \eg, the
survey~\cite{Nat}), complex tori (essentially due to
A.~Comessatti~\cite{Co-abel}), rational surfaces (A.~Degtyarev and
V.~Kharlamov~\cite{DK}), ruled surfaces
(J.-Y.~Welschinger~\cite{Welsch}), $K3$-surfaces (essentially due to
Nikulin~\cite{N1}), Enriques surfaces (see~\cite{DIK}), hyperelliptic
surfaces (F.~Catan\-ese and P.~Frediani~\cite{CF}), and some sporadic
surfaces of general type (\eg, so called \emph{Bogomolov-Miayoka-Yau
surfaces}, see Kharlamov and Kulikov~\cite{KK}).

Note that for the above classes of special surfaces topological
invariants that determine the deformation class are known. Together with
the quasi-simplicity, this implies finiteness (as the invariants take
values in finite sets). Finiteness also holds for varieties of general
type (in any dimension), as for such varieties the Hilbert scheme is
quasi-projective.

\subsection{Terminology convention}\label{terminology}
Unless stated otherwise, all complex varieties are supposed to be
nonsingular, and differentiable manifolds are~$C^\infty$. A \emph{real
variety} $(X, \conj)$ is a complex variety~$X$ equipped with an
anti-holomorphic involution~$\conj$. In spite of the fact that we work
with anti-holomorphic transformations as well, we reserve the term
\emph{isomorphism} for bi-holomoprhic maps, whereas using
\emph{anti-isomorphism} for bi-anti-holomorphic ones.

\subsection{Augmented groups and Klein actions}\label{Klein-actions}
An \emph{augmented group} is a finite group~$G$ supplied with a
homomorphism $\Gk\:G\to\{\pm1\}$. (We do not exclude the case when $\Gk$
is trivial.) Denote the kernel of~$\Gk$ by~$\Gsbp$. A \emph{Klein action}
of a group~$G$ on a complex variety~$X$ is a group action of~$G$ on~$X$
by both holomorphic and anti-holomorphic maps. Assigning~$+1$
(respectively,~$-1$) to an element of~$G$ acting holomorphically
(respectively, anti-holomorphically), one obtains a natural augmentation
$\Gk\:G\to\{\pm1\}$. An action is called \emph{holomorphic}
(respectively, \emph{properly Klein}) if $\Gk=1$ (respectively,
$\Gk\ne1$).

Replacing the complex structure~$J$ on a complex variety~$X$ with its
conjugate~$(-J)$, one obtains another complex variety, commonly denoted
by~$\bar X$, with the same underlying differentiable manifold. An
automorphism of~$X$ is as well an automorphism of~$\bar X$; it can also
be regarded as an anti-holomorphic bijection between~$X$ and~$\bar X$.
Thus, a Klein $G$-action on~$X$ can as well be regarded as a Klein action
on~$\bar X$, with the same augmentation $\Gk\:G\to\{\pm1\}$ and the same
subgroup~$\Gsbp$. These two actions are obviously diffeomorphic, but they
do not need to be isomorphic.

A Klein action of a group~$G$ on a complex variety~$X$ gives rise to the
induced action $G\to\Aut H^*(X)$, $g\mapsto g^*$, in the cohomology ring
of~$X$. Since we deal with $K3$-surfaces, which are simply connected, and
since all elements of~$G$ are orientation preserving in this dimension,
the induced action reduces essentially to the action on the group
$H^2(X)$, regarded as a lattice via the intersection index form. For our
purpose, it is more convenient to work with the \emph{twisted induced
action} $\Gt_X\:G\to\Aut H^2(X)$, $g\mapsto\Gk(g)g^*$. The latter,
considered up to conjugation by lattice automorphisms, is called the
\emph{homological type} of the original Klein action on~$X$. Clearly, it
is a topological invariant.

\subsection{Smooth deformations}\label{deformations}
A (\emph{smooth}) \emph{family}, or \emph{deformation}, of complex
varieties is a proper submersion $p\:X\to S$ with differentiable, not
necessarily compact or complex, manifolds~$X$, $S$ supplied with a
fiberwise integrable complex structure on the bundle $\Ker dp$. The
varieties $X_s=p^{-1}(s)$, $s\in S$, are called \emph{members} of the
family. Given a group~$G$, a family $p\:X\to S$ is called
\emph{$G$-equivariant} if it is supplied with a smooth fiberwise
$G$-action that restricts to a Klein action on each fiber.

Two complex varieties~$X$, $Y$ supplied with Klein actions of a group~$G$
are called \emph{equivariantly deformation equivalent} if there is a
chain $X=X_0,X_1,\ldots,X_k$ of complex varieties~$X_i$ with Klein
actions of~$G$ such that for each $i=0,\dots,k-1$ the varieties $X_i$ and
$X_{i+1}$ are $G$-isomorphic to members of a $G$-equivariant smooth
family. (By a $G$-isomorphism we mean a bi-holomoprhic map $\phi$ such
that $\phi g=g\phi$ for any $g\in G$.)

Clearly, the equivariant deformation equivalence is an equivalence
relation, $G$-equiv\-ari\-ant\-ly deformation equivalent varieties are
$G$-diffeomorphic, and the homological type of a $G$-action is a
deformation invariant.

\subsection{The principal results}\label{results}
Let $X$ be a $K3$-surface with a Klein action of a finite group~$G$. Then
$\Gsbp$ acts on the subspace $H^{2,0}(X)\cong\C$, which gives rise to a
natural representation $\Gr\:\Gsbp\to\C^*$. If $G$ is finite, the image
of~$\Gr$ belongs to the unit circle $S^1\subset\C^*$. We will refer
to~$\Gr$ as the \emph{fundamental representation} associated with the
original Klein action. It is a deformation but, in general, not
topological invariant of the action. A typical example is the same Klein
action on~$\bar X$; its associated fundamental representation is the
conjugate $\bar\Gr\:g\mapsto\overline{\Gr g}\in\C^*$.

As shown below (see~\ref{K3->geometric}), in the case of finite group
actions on a $K3$-surface~$X$ the twisted induced action~$\Gt_X$
determines the subgroup~$\Gsbp$ and `almost' determines the fundamental
representation~$\Gr\:\Gsbp\to S^1$: from~$\Gt_X$, one can recover a pair
$\Gr$, $\bar\Gr$ of complex conjugate fundamental representations.

\theorem[Finiteness Theorem]\label{Main-Finiteness}
The number of equivariant deformation classes of $K3$-surfaces with
faithful Klein actions of finite groups is finite.
\endtheorem

\theorem[Quasi-simplicity Theorem]\label{Main-q-simplicity}
Let~$X$ and~$Y$ be two $K3$-surfaces with finite group~$G$ Klein actions
of the same homological type. Assume that either
\roster
\item
the action is holomorphic, or
\item
the associated fundamental representation $\Gr$ is real, \ie, $\Gr =
\bar\Gr$.
\endroster
Then either $X$ or $\wb X$ is $G$-equivariantly deformation equivalent
to~$Y$. If the associate fundamental representation is trivial, then $X$
and~$\wb X$ are $G$-equivariantly deformation equivalent.
\endtheorem

\remark{Remark}
If $\Gr$ is non-real, the deformation classes of~$X$ and~$\wb X$ are
distinguished by the fundamental representation ($\Gr$ and $\bar\Gr$).
In~\ref{non-real} we give an example when $X$ and $\wb X$ are not
deformation equivalent even though $\Gr$ is real.
\endremark

\remark{Remark}
In~\ref{example} we discuss another example, that of a properly Klein
action whose deformation class is not determined by its homological type
and associated fundamental representation. This is a new phenomenon,
somewhat unusual for $K3$-surfaces. Note however, that the actions
constructed differ by their topology. Thus, they do not constitute a
counter-example to quasi-simplicity of $K3$-surfaces (in its weaker
form), and the problem still remains open.
\endremark

A real variety $(X, \conj)$ with a real (\ie, commuting with~$\conj$)
holomorphic $\Gsbp$-action can be regarded as a complex variety with a
Klein action of the extended group $G = \Gsbp \!\times \,\ZZ$, the
$\ZZ$-factor being generated by $\conj$. Note that, if $X$ is a
$K3$-surface with a real holomorphic $G^0$-action, the associated
fundamental representation $\Gr\:\Gsbp\to\C^*$ is real.

\corollary
Let $X$ and $Y$ be two real $K3$-surfaces with real holomorphic
$\Gsbp$-actions, so that the extended Klein actions of $G=\Gsbp\times\ZZ$
have the same homological type. Then $X$ and~$Y$ are $G$-equivariantly
deformation equivalent.
\qed
\endcorollary

The methods used in the paper can as well be applied to the study of
finite group Klein actions on $2$-dimensional complex tori. (The
corresponding version of global Torelli theorem was first discovered by
Piatetski-Shapiro and Shafarevich~\cite{PSh-Sh} and then corrected by
T.~Shioda~\cite{Shi}). The analogs of~\ref{Main-Finiteness}
and~\ref{Main-q-simplicity} for $2$-tori are
Theorems~\ref{Tori-Finiteness} (finiteness) and~\ref{Tori-q-simplicity}
(quasi-simplicity) proved in Appendix~\ref{tori}. For holomorphic actions
preserving a point this is a known result; it is contained in the
classification of finite group actions on $2$-tori by
A.~Fujiki~\cite{Fujiki}, where a complete description of the moduli
spaces is also given. (The results for holomorphic actions on Jacobians
go back to F.~Enriques and F.~Severi~\cite{ES}, and on general abelian
surfaces, back to G.~Bagnera and M.~de Franchis \cite{BF}.) We give a
short proof not using the classification, extend the results to nonlinear
Klein actions, and compare the complex conjugated actions. As a
straightforward consequence, we obtain analogous results for
hyperelliptic surfaces. A number of tools used in Appendix~\ref{tori} are
close to those used by Fujuki in his study of the relation between
symplectic actions and root systems.

Note that Theorem~\ref{Tori-q-simplicity} is stronger than its
counterpart~\ref{Main-q-simplicity} for $K3$-surfaces: one does not need
any additional assumption on the action. On the other hand, we show that,
in quite a number of cases, a $2$-torus~$X$ is not equivariantly
deformation equivalent to~$\bar X$ (see~\ref{compar}).

Together, Theorems~\ref{Main-Finiteness}, \ref{Main-q-simplicity}
and~\ref{Tori-Finiteness}, and~\ref{Tori-q-simplicity} give finiteness
and quasi-simplicity results for $K3$-surfaces, Enriques surfaces,
$2$-tori, and hyperelliptic surfaces, \ie, for all K\"ahler surfaces of
Kodaira dimension~$0$.

Among other results, not directly related to the proofs, worth mentioning
is Theorem~\ref{family}, which compares the homological types of Klein
actions on a singular $K3$-surface and on close nonsingular ones. There
also is a generalization that applies to any surface provided that the
singularities are simple.

\subsection{Idea of the proof}
As it has already been mentioned, our study is based on the global
Torelli theorem. As is known, in order to obtain a good period space, one
should \emph{mark} the $K3$-surfaces, \ie, fix isomorphisms $H^2(X)\to
L=2E_8\oplus3U$ (see~\ref{common.notation} for the notation).
Technically, it is more convenient to deal with the period space $K\GO_0$
of marked polarized $K3$-surfaces, which, in turn, is a sphere bundle
over the period space $\Per_0$ of marked Einstein $K3$-surfaces
(see~\ref{period.spaces} for details). According to
Kulikov~\cite{Kulikov}, one has $\Per_0=\Per\sminus\GD$, where $\Per$ is
a contractible homogeneous space (the space of positive definite
$3$-subspaces in $L\otimes\R$) and $\GD$ is the set of the subspaces
orthogonal to roots of~$L$.

Now, we fix a finite group~$G$ and an action $\Gt\:G\to\Aut L$. This
gives rise to the equivariant period spaces $K\GO\G_0$ and $\Per\G_0 =
\Per\G \sminus \GD$ of marked $K3$-surfaces with the given homological
type of Klein $G$-action. Note that we are only interested in
\emph{geometric} actions, \ie, those for which the spaces $\Per\G_0$ or
$K\GO\G_0$ are non-empty. Given a $K3$-surface, its markings compatible
with~$\Gt$ differ by elements of the group $\Aut_GL$ of the automorphisms
of~$L$ commuting with~$G$. Thus, the finiteness and the (weak)
quasi-simplicity problems reduce essentially to the study of the set of
connected components of the orbit space $\MS\G=\Per\G_0/\!\Aut_GL$. In
fact, the desired result (connectedness or finiteness of the number of
connected components) can be obtained with a smaller group
$A\subset\Aut_GL$, depending on the nature of the action. (A description
of such `underfactorized' moduli spaces is given
in~\ref{case-1h}--\ref{case-cah}.) Furthermore, the quotient space
$\Per\G_0/A$ can be replaced with a subspace $\Int\GG\sminus\GD$, where
$\GG$ is an appropriate convex (hence, connected) fundamental domain of
the action of~$A$ on~$\Per\G$, and it remains to enumerate the
\emph{walls} in~$\Int\GG$, \ie, the strata of $\GD\cap\Int\GG$ of
codimension~$1$.

\subsection{Contents of the paper}
In Section~\ref{actions.on.lattices} we give the basic definitions and
cite some known results on lattices and group actions on them.
In~\ref{rho} we introduce the notion of \emph{almost geometric} actions.
This notion can be regarded as the `$\Z$-independent' (\ie, defined
over~$\R$) part of the necessary condition for an action to be realizable
by a $K3$-surface. We study the invariant subspaces of an almost
geometric action and show, in particular, that such an action determines
the augmentation of the group and, up to complex conjugation, the
associated fundamental representation.

In Section~\ref{folding-walls} we introduce and study \emph{geometric}
actions, which we define in arithmetical terms. The main goal of this
section are Theorems~\ref{connectedness} and~\ref{finiteness}, which
establish certain connectedness and finiteness properties of
appropriate fundamental domains of groups of
automorphisms of the lattice preserving a given geometric action.

In Section~\ref{The.proof} we introduce the equivariant period and moduli
spaces and show that an action on the lattice is geometric (in the sense
of Section~\ref{folding-walls}) if and only if it is realizable by a
$K3$-surface. We give a detailed description of certain `underfactorized'
moduli spaces and use it to prove the main results.

Section~\ref{Degenerations} deals with equivariant degenerations of
$K3$-surfaces: we discuss the behaviour of the twisted induced action
along the walls of the period space.

In Section~\ref{examples} we discuss two examples to show that, in
general, the deformation type of a Klein action is not determined by its
homological type and associated fundamental representation.

In Appendix~\ref{tori} we treat the case of $2$-tori.

\subsection{Common notation}\label{common.notation}
We
freely use the notation~$\Z_n$ and~$\D_n$ for the
cyclic group of order~$n$ and dihedral group of order~$2n$, respectively.
We use~$A_n$, $D_n$, $E_6$, $E_7$, and~$E_8$ for the even {\bf negative}
definite lattices generated by the root systems of the same name,
and~$U$, for the hyperbolic plane (indefinite unimodular even lattice of
rank~$2$). All other non-standard symbols are explained in the text.

\subsection{Acknowledgements}
We are grateful to T.~Delzant for useful discussions on Bieberbach
groups. This research was started within the frame of
\emph{CNRS}-\emph{T\"UB\textdot ITAK} exchange program, continued with the
support of European networks EDGE and RAAG, and finished during the first
author's visits to \emph{Universit\'e de Rennes}~I and \emph{Universit\'e
Louis Pasteur}, Strasbourg, supported by~\emph{CNRS}.

\section{Actions on lattices}\label{actions.on.lattices}

\subsection{Lattices}\label{Lattices}
An \emph{\rom(integral\rom) lattice} is a free abelian group $L$ of
finite rank supplied with a symmetric bilinear form $b\:L\otimes L\to\Z$.
We usually abbreviate $b(v,w)=v\cdot w$ and $b(v,v)=v^2$. For any
ring~$\GL\supset\Z$ we use the same notation $b$ (as well as $v\cdot w$
and~$v^2$) for the linear extension
$(v\otimes\Gl)\otimes(w\otimes\mu)\mapsto (v\cdot w)\Gl\mu$ of~$b$ to
$L\otimes\GL$. A lattice~$L$ is called \emph{even} if $v^2=0\bmod2$ for
all $v\in L$; otherwise, $L$ is called \emph{odd}. Let
$L\spcheck=\Hom(L,\Z)$ be the dual abelian group. The lattice~$L$ is
called \emph{nondegenerate} (\emph{unimodular}) if the \emph{correlation
homomorphism} $L\to L\spcheck$, $v\mapsto b(v,\,\cdot\,)$, is a
monomorphism (respectively, isomorphism). The cokernel of the correlation
homomorphism is called the \emph{discriminant group} of~$L$ and denoted
by $\discr L$. The group $\discr L$ is finite (trivial) if and only if
$L$ is nondegenerate (respectively, unimodular).

The assignment $(x\bmod L,y\bmod L)\mapsto(x\cdot y)\bmod\Z$, $x,y\in
L\spcheck$ is a well defined bilinear form $b\:\discr L\otimes\discr
L\to\Q/\Z$. If $L$ is even, there also is a quadratic extension
$q\:\discr L\to\Q/2\Z$ of~$b$ given by $x\bmod L\mapsto x^2\bmod2\Z$.

Given a lattice~$L$, we denote by $\Gs\sbp L$ and $\Gs\sbm L$ its inertia
indexes and by $\Gs L=\Gs\sbp L-\Gs\sbm L$, its signature. We call a
nondegenerate lattice~$L$ \emph{elliptic} (respectively,
\emph{hyperbolic}) if $\Gs\sbp L=0$ (respectively, $\Gs\sbp L=1$). The
terminology is not quite standard: we change the sign of the forms, and
we treat a positive definite lattice of rank~$1$ as hyperbolic. This is
caused by the fact that our lattices are related (explicitly or
implicitly) to the Neron-Severi groups of complex surfaces.

A sublattice $M\subset L$ is called \emph{primitive} if the quotient
$L/M$ is torsion free. Given a sublattice $M\subset L$, we denote by
$M\prhull$ its \emph{primitive hull} in $L$, \ie, the minimal primitive
sublattice containing~$M$: $M\prhull=\{v\in L\mid\text{$kv\in M$ for some
$k\in\Z$, $k\ne0$}\}$.

An element $v \in L$ of square~$(-2)$ is called a
\emph{root}.\footnote{Traditionally, the roots are the elements of square
$(-2)$ or $(-1)$. We exclude the case of square $(-1)$ as we only
consider even lattices.} A \emph{root system} is a lattice generated
(over $\Z$) by roots. Recall that any elliptic root system decomposes,
uniquely up to order of the summands, into orthogonal sum of irreducible
elliptic root systems, \ie, those of type $A_n$, $D_n$, $E_6$, $E_7$, or
$E_8$.

\subsection{Automorphisms}
An \emph{isometry} (\emph{dilation}) of a lattice~$L$ is an automorphism
$a\:L\to L$ preserving the form (respectively, multiplying the form by a
fixed number $\ne0$.) All isometries of~$L$ constitute a group; we denote
it by $\Aut L$. If $L$ is non-degenerate, there is a natural
representation $\Aut L\to\Aut\discr L$. Denote its kernel $\Aut^0L$. It
is a finite index normal subgroup of $\Aut
L$ consisting of the `universally extensible' automorphisms. More
precisely, an automorphism~$a$ of~$L$ belongs to $\Aut^0L$ if and only if
$a$ extends to any suplattice $L'\supset L$ identically on~$L\ort$.

Given a vector~$v\in L$, $v^2\ne0$, denote by $\hps{v}$ the reflection
against the hyperplane orthogonal to~$v$, \ie, the isometry of
$L\otimes\R$ defined by $x\mapsto x-((x\cdot v)/v^2)v$. If
$\hps{v}(L)\subset L$ (which is always the case when $v^2=\pm1$ or
$\pm2$), we use the same notation for the induced automorphism of~$L$.
The subgroup $W(L)\subset\Aut L$ generated by the reflections against the
hyperplanes orthogonal to roots of~$L$ is called the \emph{Weil group}
of~$L$. Clearly, $W(L)$ is a normal subgroup of $\Aut L$ and
$W(L)\subset\Aut^0L$.

We recall a few facts on automorphisms of root systems; details can be
found, \eg, in~\cite{Bourbaki}. Let $R$ be an elliptic root system. The
hyperplanes orthogonal to roots in~$R$ divide the space $R\otimes\R$ into
several connected components, called \emph{cameras} of~$R$, and the Weil
group $W(R)$ acts transitively on the set of cameras. For each camera~$C$
of~$R$ there is a canonical semi-direct product decomposition $\Aut
R=W(R)\rtimes S_C$, where $S_C\subset\grO(R\otimes\R)$ is the group of
symmetries of~$C$. (As an abstract group, $S_C$ can be identified with
the group of symmetries of the Dynkin diagram of~$R$.) In particular, if
an element $g\in\Aut R$ preserves~$C$, one has $g\in S_C$. More
generally, if $g$ preserves a face $C'\subset C$, then in the
decomposition $g=sw$, $s\in S_C$, the element $w$ belongs to the Weil
group of the root system generated by the roots of~$R$ orthogonal
to~$C'$.

\subsection{Actions}\label{lattice-notation}
Let $G$ be a group. A $G$-action on a lattice~$L$ is a representation
$\Gt\:G\to\Aut L$. In what follows we always assume~$G$ finite. Given a
ring $\GL\supset\Z$, we use the same notation~$\Gt$ for the extension
$g\mapsto\Gt g\otimes\id_\GL$ of the action to $L\otimes\GL$. Denote by
$\Aut_G(L\otimes\GL)$ the group of $G$-equivariant $\GL$-isometries of
$L\otimes\GL$, \ie, the centralizer of $\Gt G$ in $\Aut(L\otimes\GL)$,
and let $W_G(L)=W(L)\cap\Aut_GL$ and $\Aut_G^0L=\Aut^0L\cap\Aut_GL$.

A submodule $M\subset L\otimes\GL$ is called \emph{$G$-invariant} if $\Gt
g(M)\subset M$ for any $g\in G$; it is galled \emph{$G$-characteristic}
if $a(M)\subset M$ for any $a\in\Aut_G(L\otimes\GL)$.

Let $\K\subset\C$ be a field. For an irreducible $\K$-linear
representation~$\xi$ of~$G$, we denote by $L_\xi(\K)$ the $\xi$-isotypic
subspace of~$L\otimes\K$, \ie, the maximal invariant subspace of
$L\otimes\K$ that is a sum of irreducible representations isomorphic
to~$\xi$. Given a subfield $\k\subset\K$, denote by $L_\xi(\k)$ the
minimal $\k$-subspace of $L\otimes\k$ such that
$L_\xi(\k)\otimes_\k\K\supset L_\xi(\K)$, and for a subring
$\gD\subset\k$, $\gD\ni1$, let $L_\xi(\gD)=L_\xi(\k)\cap(L\otimes\gD)$.
Clearly, $L_\xi(\k)$ is the space of an isotypic $\k$-representation
of~$G$, and $L_\xi(\gD)$ is $G$-invariant and $G$-characteristic. If $\k$
is an algebraic number field and $\gD$ is an order in~$\k$, then
$L_\xi(\gD)$ is a finitely generated abelian group and
$L_\xi(\k)=L_\xi(\gD)\otimes_{\gD}\k$.

We use the shortcut $L^G$ for $L_1(\Z)=\{x\in L\mid\text{$gx=x$ for all
$g\in G$}\}$.

\subsection{Extending automorphisms}\label{automorphisms}
Below, we recall a few simple facts on extending automorphisms of
lattices. All the results still hold if the lattices involved are
supplied with an action of a finite group~$G$ and the automorphisms are
$G$-equivariant. One can also consider lattices defined over an order in
an algebraic number field.

\lemma\label{max-extension}
Let~$M$ be a nondegenerate lattice and $M'\subset M$ a sublattice of
finite index. Then the groups $\Aut M$ and $\Aut M'$ have a common finite
index subgroup.
\qed
\endlemma

\lemma\label{extension}
Let~$M$ be a lattice and $M'\subset M$ a nondegenerate sublattice. Then
the group of automorphisms of $M'$ extending to $M$ has finite index in
$\Aut M'$.
\qed
\endlemma

\lemma\label{rational->integer}
Let~$M$ be a nondegenerate lattice and $A$ a group acting by isometries
on $M\otimes\Q$. Assume that there is a finite index sublattice
$M'\subset M$ such that $a(M')\subset M$ for any $a\in A$. Then~$A$ has a
finite index subgroup acting on~$M$.
\endlemma

\proof
It suffices to apply~\ref{max-extension} to the $A$-invariant sublattice
$\sum_{a\in A}a(M')\subset M$.
\endproof

\corollary\label{two-lattices}
Let~$M\spp$ and~$M\spm$ be two nondegenerate lattices and $J\:M\spm\to
M\spp$ a dilation invertible over $\Q$. Then there exists a finite index
subgroup $A\spp \subset \Aut M\spp$ such that the correspondence
$a\mapsto a\oplus J^{-1}aJ$ restricts to a well defined homomorphism
$A\spp\to\Aut(M\spp\oplus M\spm)$.
\qed
\endcorollary

\subsection{Fundamental polyhedra}
Given a real vector space~$V$ with a nondegenerate quadratic form, we
denote by $\HS(V)$ the space of maximal positive definite subspaces
of~$V$. Note that $\HS(V)$ is a contractible
space of non positive curvature. If
$\Gs\sbp V=1$ (\ie, $V$ is hyperbolic), one can define $\HS(V)$ as the
projectivization $\HC(V)/\R^*$ of the positive cone $\HC(V)=\{x\in V\mid
x^2>0\}$.

Fix an algebraic number field $\k\subset\R$ and let $\gD$ be the ring of
integers of~$\k$. Consider a hyperbolic integral lattice~$M$ and a
hyperbolic sublattice $M'\subset M\otimes\k$ defined over~$\gD$, \ie,
such that $\gD M' \subset M'$. Let $\HS'=\HS(M'\otimes_\gD\R)$. Then any
group~$A$ acting by isometries on~$M$ and preserving~$M'$ acts on~$\HS'$.
Since~$M$ is a hyperbolic integral lattice and $(M')\ort \subset M$ is
elliptic, the induced action is discrete, and the Dirichlet domain with
center at a generic $\k$-rational point of $\HS'$ is a $\k$-rational
polyhedral fundamental domain of the action. Any such domain will be
called a \emph{rational Dirichlet polyhedron} of~$A$ (in~$\HS'$).

The following theorem treats the classical case where $M = M'$ is an
integral lattice and $A = \Aut M$. It is due to C.~L.~Siegel~\cite{Sie},
H.~Garland, M.~S.~Raghunathan~\cite{Garland}, and
N.~J.~Wielenberg~\cite{Wielenberg}.

\theorem\label{Siegel}
Let $M$ be a hyperbolic integral lattice. Then the rational Dirichlet
polyhedra of the full automorphism group $\Aut M$ in $\HS(M)$ are finite.
Unless $M$ has rank~$2$ and represents~$0$, the polyhedra have finite
volume.
\qed
\endtheorem

\corollary\label{convex-hull}
Let $M$ be a hyperbolic integral lattice. Then the closure in
$\HS(M)\cup\partial\HS(M)$ of any rational Dirichlet polyhedron of $\Aut
M$ in $\HS(M)$ is the convex hull of a finite collection of rational
points.
\qed
\endcorollary

\subsection{The fundamental representations}\label{rho}
Let $\Gt\:G\to\Aut L$ be a finite group action on a nondegenerate
lattice~$L$ with $\Gs\sbp L=3$. We will say that $\Gt$ is \emph{almost
geometric} if there is a $G$-invariant flag $\ell\subset\gw$, where
$\gw\subset L\otimes\R$ is a positive definite $3$-subspace and~$\ell$ is
a $1$-subspace with trivial $G$-action.

\lemma\label{d-subspace}
Let $\Gt\:G\to\Aut L$ be a finite group action on a lattice~$L$ with
$d=\Gs\sbp L>0$. Then, for any positive definite $G$-invariant
$d$-subspace $\gw\subset L\otimes\R$, the
induced action $\Gt_{\gw}\:G\to\grO(\gw)=\grO(d)$ is determined by~$\Gt$
up to conjugation in~$\grO(d)$. In particular, the augmentation
\smash{$\Gk\:G\to\grO(\gw)@>{\raise-1pt\vbox{\hbox{$\scriptstyle\det$}}}>>\{\pm1\}$}
is uniquely determined by~$\Gt$.
\endlemma

\proof
Given another subspace~$\gw'$ as in the statement, the orthogonal
projection $\gw'\to\gw$ is non-degenerate and $G$-equivariant. Hence, the
induced representations $\Gt_{\gw},\Gt_{\gw'}\:G\to\grO(d)$ are
conjugated by an element of $\grGL(d)$. Since $G$ is finite, they are
also conjugated by an element of~$\grO(d)$. Indeed, it is sufficient to
treat the case of irreducible representation, where the result follows
from the uniqueness of a $G$-invariant
scalar product up to a constant factor.
\endproof

Given an almost geometric action $\Gt\:G\to\Aut L$, we will always
assume~$G$ augmented via~$\Gk$ above, so that an element $c\in G$ does
not belong to~$\Gsbp=\Ker\Gk$ if and only if it reverses the orientation
of~$\gw$. From~\ref{d-subspace} it follows that the existence of a
$1$-subspace~$\ell$ with trivial $G$-action does not depend on the choice
of a $G$-invariant positive definite $3$-subspace~$\gw$. Furthermore, the
induced action on $\gw_0=\ell\ort\subset\gw$ is also independent
of~$\gw$. Choosing an orientation of~$\gw_0$, one obtains a
$2$-dimensional representation $\Gr\:\Gsbp\to\grSO(\gw_0)=S^1$. In what
follows, we identify~$S^1$ with the unit circle in~$\C$ and often regard
representations in~$S^1$ as one-dimensional complex representations. In
particular, we consider the spaces (lattices) $L_\Gr(\GL)$
(see~\ref{lattice-notation}) associated with~$\Gt$. Note that $L_\Gr(\C)$
is the $\Gr$-eigenspace of~$\Gsbp$. Changing the orientation of~$\gw_0$
replaces~$\Gr$ with its conjugate~$\bar\Gr$. In view of~\ref{d-subspace},
the unordered pair $(\Gr,\bar\Gr)$ is determined by~$\Gt$; we will call
$\Gr$ and~$\bar\Gr$ the \emph{fundamental representations} associated
with~$\Gt$. The order of the image $\Gr(\Gsbp)$ is called the
\emph{order} of~$\Gt$ and is denoted $\ord\Gt$.

\lemma\label{Jrho}
Let $\xi\:\Gsbp\to S^1$ be a non-real representation \rom(\ie,
$\bar\xi\ne\xi$\rom). Then the map $L_\xi(\C)\to L_\xi(\R)$,
$\Go\mapsto\frac12(\Go+\bar\Go)$, is an isomorphism of $\R$-vector
spaces. In particular, the space $L_\xi(\R)$ inherits a natural complex
structure $J_\xi$ \rom(induced from the multiplication by~$i$ in
$L_\xi(\C)$\,\rom), which is an anti-selfadjoint isometry. One has
$J_{\bar\xi}=-J_\xi$.
\endlemma

\proof\nofrills{\/ }
is straightforward. The metric properties of $J_\xi$ follow from the fact
that $\Go^2=0$ for any eigenvector~$\Go$ (of any isometry) corresponding
to an eigenvalue~$\Ga$ with $\Ga^2\ne1$.
\endproof

\lemma\label{real.structure}
Let $\Gt$ be an almost geometric action and $\Gr$ an associated
fundamental representation. Assume that $\Gk\ne1$. Then any element $c\in
G\sminus \Gsbp$ restricts to an involution $c_\Gr\:L_\Gr(\R)\to
L_\Gr(\R)$. If $\Gr$ is not real, then $c_\Gr$ is
$J_\Gr$-anti-linear\rom; in particular, the $(\pm1)$-eigenspaces
$V_\Gr\sppm$ of~$c_\Gr$ are interchanged by~$J_\Gr$.
\endlemma

\proof
Clearly, $c$ takes $\Gr$-eigenvectors of~$\Gsbp$ to $\Gr^c$-eigenvectors,
where $\Gr^c$ is the representation $g\mapsto\Gr(c^{-1}gc)$. Since, by
the definition of fundamental representations, there is a
$\Gr$-eigenvector~$\Go$ taken to a $\bar\Gr$-eigenvector, one has
$\Gr^c=\bar\Gr$ and the space $L_\Gr(\R)$ is $c$-invariant. Furthermore,
the vector $\Re\Go$ is invariant under $c_\Gr^2$. Since $c^2\in \Gsbp$,
one has $c_\Gr^2=\id$.

If $\Gr$ is non-real, then $c$ interchanges $L_\Gr(\C)$ and
$L_{\bar\Gr}(\C)$. Since~$c$ commutes with the complex conjugation, the
isomorphism $\Go\mapsto\frac12(\Go+\bar\Go)$ (see~\ref{Jrho})
conjugates~$c_\Gr$ with the anti-linear involution $\Go\mapsto
c(\bar\Go)$ on $L_\Gr(\C)$.
\endproof

\lemma\label{CnDn}
Let~$\Gt$ be an almost geometric action, $\Gr$ an associated fundamental
representation, and $\k\subset\R$ a field. Then the space $L_\Gr(\k)$ is
$G$-invariant and the induced $G$-action on $L_\Gr(\k)$ factors through
an action of the cyclic group~$\Z_n$ \rom(if $\Gk=1$\rom) or the dihedral
group~$\D_n$ \rom(if $\Gk\ne1$\rom), where $n=\ord\Gt$. The induced
$\Z_n$-action is $\k$-isotypic\rom; the $\D_n$-action is $\k$-isotypic
unless $n\le2$.
\endlemma

\proof
All statements are obvious if $\Gk=1$. Assume that $\Gk\ne1$ and pick an
element $c\in G\sminus \Gsbp$. The intersection $Q=L_\Gr(\k)\cap
c(L_\Gr(\k))$ is defined over~$\k$, and $Q\otimes_\k\R$
contains~$L_\Gr(\R)$ (see~\ref{real.structure}). Hence, $Q\supset
L_\Gr(\k)$ and $L_\Gr(\k)$ is $G$-invariant. Further, the endomorphisms
$c^2$ and $g-c^{-1}gc$ of $L_\Gr(\k)\otimes_\k\R$ (where $g\in \Gsbp$)
are defined over~$\k$ and annihilate $L_\Gr(\R)$
(see~\ref{real.structure} again); due to the minimality of~$L_\Gr(\k)$,
they are trivial.
\endproof

\section{Folding the walls}\label{folding-walls}

\subsection{Geometric actions}\label{geometric-actions}
A finite group action $\Gt\:G\to\Aut L$ on an even nondegenerate
lattice~$L$ with $\Gs\sbp L=3$ is called \emph{geometric} if it is almost
geometric and the sublattice $\eL = (L^G + L_\Gr(\Z))\ort$ contains no
roots, where~$\Gr$ is a fundamental representation of~$\Gt$.

Consider a geometric action~$\Gt$ and fix an associated fundamental
representation~$\Gr$. If $\Gk\ne1$, fix an element $c\in G\sminus \Gsbp$
and denote by~$V_\Gr\sppm$ and $V\sppm$ its $(\pm1)$-eigenspaces in
$L_\Gr(\R)$ and $L_\Gr(\Q)$, respectively (see~\ref{real.structure}
and~\ref{CnDn}). Let $M\sppm=V\sppm\cap L$ be the $(\pm1)$-eigenlattices
of~$c$ in $L_\Gr(\Z)$. If $\Gr\ne1$, the spaces $V_\Gr\sppm$ and $V\sppm$
are hyperbolic. The following lemma is a straightforward consequence
of~\ref{real.structure} and~\ref{CnDn}.

\lemma
The subspaces $V_\Gr\sppm$ and $V\sppm$ and the sublattices~$M\sppm$ are
$G$-char\-ac\-ter\-is\-tic\rom; they are $G$-invariant if and only if
$\ord\Gt\le2$. If $\Gr\ne1$, there is a well defined action of $\Aut_GL$
on $\HS(V_\Gr\sppm)$\rom; it is discrete and, up to isomorphism,
independent of the choice of an element $c\in G\sminus \Gsbp$.
\qed
\endlemma

In view of this lemma one can consider corresponding $G$-actions and
introduce the following rational Dirichlet polyhedra.
\roster\widestnumber\item{}
\item"--\hfill"
$\GG_1\subset\HS(L^G\otimes\R)$ is a rational Dirichlet polyhedron of
$W_G((L^G\oplus\eL)\prhull)$; it is defined whenever $\Gr\ne1$, so that
$\Gs\sbp L^G=1$.
\item"--\hfill"
$\GG_\Gr\sppm\subset\HS(V_\Gr\sppm)$ are some rational Dirichlet
polyhedra of $W_G((M\sppm\oplus\eL)\prhull)$; they are defined whenever
$\Gr$ is real and $\Gk\ne1$. (To define $\GG\spp_\Gr$, one needs to
assume, in addition, that $\Gr\ne1$, so that $\Gs\sbp M\spp=1$.)
\item"--\hfill"
$\GS_\Gr\sppm\subset\HS(V_\Gr\sppm)$ are some rational Dirichlet
polyhedra of $\Aut_G^0(L_\Gr(\Z))$; they are defined whenever $\Gr$ is
non-real and $\Gk\ne 1$.
\endroster

Given a vector $v\in L$, put $\sh(v)=\{x\in L\otimes\R\mid x\cdot v=0\}$
and introduce the following notation:
\roster\widestnumber\item{}
\item"--\hfill"
$\sh_1(v)=\sh(v)\cap(L^G\otimes\R)$;
\item"--\hfill"
if $\Gr$ is real and $\Gk\ne1$, then $\sh_\Gr\sppm(v)=\sh(v)\cap
V_\Gr\sppm$;
\item"--\hfill"
if $\Gr$ is non-real, then $\sh_\Gr(v)=\{x\in L_\Gr(\R)\mid x\cdot
v=J_\Gr x\cdot v=0\}$; if, besides, $\Gk\ne1$, then
$\sh_\Gr\sppm(v)=\sh_\Gr(v)\cap V_\Gr\sppm$.
\endroster
We use the same notation $\psh_1(v)$ and $\psh_\Gr\sppm(v)$ for the
projectivizations of the corresponding spaces in $\HS(L^G\otimes\R)$ and
$\HS(V_\Gr\sppm)$, respectively (whenever the space is hyperbolic).

The goal of this section is to prove the following two theorems.

\theorem\label{connectedness}
Let $\Gt\:G\to\Aut L$ be a geometric action and $\Gr$ an associated
fundamental representation. If $\Gr\ne1$, then for any root $v\in
L_\Gr(\Z)\ort$ the intersection $\sh_1(v)\cap\Int\GG_1$ is empty. If
$\Gr$ is real and $\Gk\ne1$, then for any root $v\in(L^G\oplus
M\sp\mp)\ort$ the intersection $\sh_\Gr\sppm(v)\cap\Int\GG_\Gr\sppm$ is
empty. \rom(For $\GG_\Gr\spp$ to be well defined, one needs to assume, in
addition, that $\Gr\ne1$.\rom)
\endtheorem

\theorem\label{finiteness}
Let $\Gt\:G\to\Aut L$ be a geometric action with non-real associated
fundamental representation~$\Gr$ and $\Gk\ne1$. Then $\GS_\Gr\sppm$
intersects finitely many distinct subspaces $\psh_\Gr\sppm(v)$ defined by
roots $v\in(L^G)\ort$.
\endtheorem

Theorem~\ref{connectedness} is proved at the end of~\ref{inv-sublattice}.
Theorem~\ref{finiteness} is proved in~\ref{proof}.

\subsection{Walls in the invariant sublattice}\label{inv-sublattice}

\theorem\label{inv}
Let~$N$ be an even lattice and~$G$ a finite group acting on~$N$ so that
$(N\G)\ort\subset N$ is negative definite. Let $v\in N$ be a root whose
projection to $N\G\otimes\R$ has negative square. Then either
\roster
\item
the orthogonal complement $(N\G)\ort$ contains a root, or
\item
there is an element of $W_G(N)$ whose restriction to $N\G$ is the
reflection against the hyperplane $\sh(v)\cap(N\G\otimes\R)$.
\endroster
\endtheorem

\corollary\label{inv-poly}
In the above notation, assume that $N$ is hyperbolic and $(N\G)\ort$
contains no roots. Then for any root $v\in N$ the intersection of
$\sh(v)$ with the interior of a rational Dirichlet polyhedron of $W_G(N)$
in $\HS(N\G)$ is empty.
\qed
\endcorollary

To prove Theorem~\ref{inv} we need a few facts on automorphisms of root
systems. Let~$R$ be an even root system and~$G$ a finite group acting
on~$R$. The action is called \emph{admissible} if the orthogonal
complement $(R^G)\ort$ contains no roots, and it is called
\emph{$b$-transitive} if there is a root whose orbit generates~$R$.

\lemma\label{admissible}
Given a finite group~$G$ action on an elliptic root system~$R$, the
following statements are equivalent\rom:
\roster
\item
the action is admissible\rom;
\item
the action preserves a camera of~$R$\rom;
\item
the action factors through the action of a subgroup of the symmetry group
of a camera of~$R$.
\endroster
\endlemma

\proof
An action is admissible if and only if $R^G$ does not belong to a wall
$\sh(v)$ defined by a root~$v\in R$. On the other hand, $R^G$ contains an
inner point of a camera if and only if this camera is preserved by the
action.
\endproof

\corollary\label{two-actions}
Up to isomorphism, there are two faithful admissible $b$-transi\-tive
actions on irreducible even root systems\rom: the trivial action on~$A_1$
and a $\ZZ$-action on~$A_2$ interchanging two roots~$u$, $v$ with $u\cdot
v=1$.
\endcorollary

\proof
The statement follows from Lemma~\ref{admissible}, the classification of
irreducible root systems, and the natural bijection between the
symmetries of a camera and the symmetries of its Dynkin diagram.
\endproof

\proof[Proof of Theorem~\ref{inv}]
Pick a vector~$v$ as in the statement, and consider the sublattice $R
\subset N$ generated by the orbit of~$v$. Under the assumptions, $R$ is
an even root system, and the induced $G$-action on~$R$ is $b$-transitive.
Assume that the action on~$R$ is admissible (as otherwise $(R^G)\ort$,
and thus $(N^G)\ort$, would contain a root). Then, in view
of~\ref{two-actions}, the lattice~$R$ splits into orthogonal sum of
several copies of either~$A_1$ or~$A_2$, and the vector $\bar
v=\sum_{g\in G}g(v)$ has the form $\sum m_ia_i$, $m_i\in\Z$, where each
$a_i$ is a generator of $A_1$ or the sum of two generators of $A_2$
interchanged by the action. Since the $a_i$'s are mutually orthogonal
roots, the composition of the reflections $\hps{a_i}$ is the desired
automorphism of~$N$.
\endproof

\proof[Proof of Theorem~\ref{connectedness}]
The statement for $\GG_1$ follows immediately from Theorem \ref{inv}
applied to $N=L_\Gr(\Z)\ort$. To prove the assertion for~$\GG_\Gr\sppm$,
consider the induced $G$-action $\Gt_{\gw}\:G\to\grO(\gw)$, where $\gw$
is as in the definition of an almost geometric action, see~\ref{rho}, and
note that, under the hypotheses ($\Gr\ne1$ is real), $\Gt_{\gw}$ factors
through the abelian subgroup $C\subset\grO(\gw)$ generated by the central
symmetry~$c$ and a reflection~$s$. Thus, the statement for~$\GG_\Gr\spp$
(respectively,~$\GG_\Gr\spm$) follows from~\ref{inv} applied to the
lattice $N=(L^G\oplus M\spm)\ort$ (respectively, $N=(L^G\oplus
M\spp)\ort$) with the twisted action $g\mapsto r(g)\Gt(g)$, where
$r\:G\to\{\pm 1\}$ is the composition of~$\Gt_{\gw}$ and the homomorphism
$c\mapsto -1$, $s\mapsto 1$ (respectively, $c\mapsto -1$, $s\mapsto -1$).
\endproof

\subsection{The group $\Aut_GL$}\label{AutM}
Let, as before, $\Gt\:G\to\Aut L$ be an almost geometric action and $\Gr$
a fundamental representation of~$\Gt$. Recall (see~\ref{CnDn}) that the
induced $G$-action on $L_\Gr(\Z)$ factors through the group $G'=\Z_n$ (if
$\Gk=1$) or~$\D_n$ (if $\Gk\ne1$), where $n=\ord\Gt>2$. Let~$\K$ be the
cyclotomic field $\Q(\exp(2\pi i/n))$ and let $\k\subset\K$ be the real
part of~$\K$, \ie, the extension of~$\Q$ obtained by adjoining the real
parts of the primitive $n$-th roots of unity. Both~$\K$ and~$\k$ are
abelian Galois extensions of~$\Q$. Denote by $\gD_\K$ and $\gD$ the rings
of integers of~$\K$ and~$\k$, respectively. Unless specified otherwise,
we regard~$\k$ and~$\K$ as subfields of~$\C$ via their standard
embeddings. An isotypic $\k$-representation of~$G'$ corresponding to a
pair of conjugate primitive $n$-th roots of unity will be called
\emph{primitive}.

\lemma\label{all->pair}
For any primitive irreducible $\k$-representation~$\xi$ of~$G'$, the
restriction homomorphism $\Aut_GL\to\Aut_GL_\xi(\gD)$ is well defined and
its image has finite index. If $L=L_\xi(\Z)$, the restriction is a
monomorphism.
\endlemma

\proof
In view of~\ref{extension} and~\ref{CnDn}, it suffices to consider the
case when $L=L_\xi(\Z)$ and $G=G'$. The restriction homomorphism is well
defined as any $G$-equivariant isometry of $L_\xi(\Z)$, after extension
to $L_\xi(\Z)\otimes\k$, must preserve the $\k$-isotypic subspaces. It is
a monomorphism, since $L_\xi(\Q)$ is the minimal $\Q$-vector space such
that $L_\xi(\Q)\otimes\k$ contains $L_\xi(\k)$. (If an element
$g\in\Aut_GL_\xi(\Z)$ restricts to the identity of $L_\xi(\gD)$, then
$\Ker(g-\id)$ is a $\Q$-vector space with the above property; hence, it
must contain $L_\xi(\Q)$.)

It remains to prove that, up to finite index, any $G$-equivariant
$\gD$-automorphism~$g$ of $L_\xi(\gD)$ extends to a $G$-equivariant
automorphism of $L_\xi(\Z)\otimes\gD$ defined over~$\Z$. Up to finite
index, one has an orthogonal decomposition
$L_\xi(\Z)\otimes\gD\supset\bigoplus L_{\xi_i}(\gD)$, the summation over
all primitive irreducible representations~$\xi_i$ of~$G$. For each such
representation~$\xi_i$ there is a unique element $\gg_i\in\Gal(\k/\Q)$
such that $\xi_i=\gg_i\xi$, and the automorphism
$\bigoplus\gg_ig\gg_i^{-1}$ of $\bigoplus L_{\xi_i}(\gD)$ is Galois
invariant, \ie, defined over~$\Z$.
\endproof

Let now $\Gk\ne1$, \ie, $G'=\D_n$. Put
$M\sppm_\xi=V\sppm_\xi\cap(L\otimes\gD)$ and denote by $\Aut M\sppm_\xi$
the group of isometries of $M\sppm_\xi$ defined over~$\gD$. (Note that
$V_\Gr\sppm$ are defined over~$\k$ and thus can be regarded as subspaces
of $L_\Gr(\k)$.)

\lemma\label{pair->half}
For any primitive irreducible $\k$-representation~$\xi$ of~$G'=\D_n$, the
restriction homomorphism $\Aut_GL_\xi(\gD)\to\Aut M\sppm_\xi$ is a well
defined monomorphism, and its image has finite index.
\endlemma

\proof
Again, it suffices to consider the case $G=G'$. Obviously, any
$G$-equivariant automorphism of $L_\xi(\gD)$ preserves~$M\sppm_\xi$. To
prove the converse (say for $M_\xi\spp$), note that, up to a factor, the
map $J_\xi$ is defined over~$\k$ (as this is obviously true for an
irreducible representation, where $\dim_\k V\spp_\xi=\dim_\k
V\spm_\xi=1$), \ie, there is a dilation $J=kJ_\xi$ of $L_\xi(\k)$
interchanging~$V\spp_\xi$ and~$V\spm_\xi$. Furthermore, the factor can be
chosen so that $J(M\spm_\xi)\subset M\spp_\xi$. Since any extension of an
isometry $a\in\Aut M\spp_\xi$ to $L_\xi(\gD)$ must commute with~$J$, on
$M\spp_\xi\oplus M\spm_\xi$ it must be given by $a\oplus J^{-1}aJ'$. On
the other hand, due to~\ref{two-lattices}, the latter expression does
define an extension for all~$a$ in a finite index subgroup of $\Aut
M\spp_\xi$.
\endproof

\corollary\label{finite.copy}
The polyhedron~$\GS_\Gr\sppm$ is the union of finitely many copies of a
rational Dirichlet polyhedron of $\Aut M\sppm_\Gr$ in~$\HS_\Gr\sppm$.
\qed
\endcorollary

\subsection{Dirichlet polyhedra:
the case $\Gf(\ord\Gr)=2$}\label{case.gf=2}
Recall that $\Gf$ is the
Euler function, \ie, $\Gf(n)$ is the number of positive integers $<n$
prime to~$n$. Alternatively, $\Gf(n)$ is the degree of the cyclotomic
extension of~$\Q$ of order~$n$. Consider a hyperbolic sublattice
$M\subset L$ and denote by $\HS=\HS(M\otimes\R)$ the corresponding
hyperbolic space. Given a vector $v\in M$, let
$\psh_M(v)=(\sh(v)\cap\HC(M\otimes\R))/\R^*\subset\HS$.

\lemma\label{line}
Let $\ell\subset\HS$ be a line whose closure intersects the absolute
$\dHS$ at rational points. Then for any integer~$a$ there are at most
finitely many vectors $v\in M$ such that $v^2=a$ and the hyperplane
$\psh_M(v)$ intersects~$\ell$.
\endlemma

\proof
Let $u_1,u_2\in M$ be some vectors corresponding to the intersection
points $\ell\cap\dHS$. Then $u_1, u_2$ span a (scaled) hyperbolic plane
$U\subset M$ and the orthogonal complement $U\ort\subset M$ is elliptic.
Therefore, $U\oplus U\ort$ is of finite index~$d$ in $M$.

Let $v$ be a vector as in the statement. Since $\psh_M(v)$
intersects~$\ell$, one has $v=\Gl bu_1+(\Gl-1)bu_2+v'$ for some $v'\in
\frac{1}{d}\,U\ort$ and $\Gl\in(0,1)$. Thus, the equation $v^2=a$ turns
into $-b^2\Gl(1-\Gl)+(v')^2=a$. Since $dv'$ belongs to a negative
definite lattice, $\Gl(1-\Gl)>0$, and both $\Gl bd$ and $(1-\Gl)bd$ are
integers, this equation has finitely many solutions.
\endproof

\corollary\label{polytope}
Let~$Q\subset\HS$ be a polyhedron whose closure in $\HS\cup\partial\HS$
is a convex hull of finitely many rational points. Then for any
integer~$a$ there are at most finitely many vectors $v\in M$ such that
$v^2=a$ and the hyperplane $\psh_M(v)$ intersects~$Q$.
\endcorollary

\proof
Each edge of~$Q$ either is a compact subset of $\HS$ or has a rational
endpoint on the absolute. In the former case, the edge intersects
finitely many hyperplanes $\psh_M(v)$, as they form a discrete set. In
the latter case, both the intersection points of the absolute and the
line containing the edge are rational, and the edge intersects finitely
many hyperplanes $\psh_M(v)$ due to~\ref{line}. Finally, if a hyperplane
does not intersect any edge of~$Q$, it contains at least $\dim\HS$
vertices of~$Q$ at the absolute and is determined by those vertices.
Since $Q$ has finitely many vertices, the number of such hyperplanes is
also finite.
\endproof

\corollary[Corollary \rm(of~\ref{polytope}
and~\ref{convex-hull})]\label{Siegel.M} Assume that $\Gk\ne 1$ and
$\Gf(\ord\Gt) = 2$ \rom(so that $M_\Gr\sppm$ are defined over $\Z$\rom)
and let $\Pi_\Gr\sppm$ be some rational Dirichlet polyhedra of $\Aut
M_\Gr\sppm$ in~$\HS_\Gr\sppm$. Then for any integer~$a$ there are at most
finitely many vectors $v\in M_\Gr\sppm$ such that $v^2=a$ and the
subspace $\psh_\Gr\sppm(v)$ intersects~$\Pi_\Gr\sppm$
or~$J_\Gr(\Pi_\Gr\sp\mp)$.
\qed
\endcorollary

\subsection{Dirichlet polyhedra: the case $\Gf(\ord\Gt)\ge4$}\label{case.gf>2}
Recall that an algebraic number field~$F$ has exactly $\deg(F/\Q)$
distinct embeddings to~$\C$. Denote by~$r(F)$ the number of real
embeddings (\ie, those whose image is contained in~$\R$), and by~$c(F)$,
the number of pairs of conjugate non-real ones. Clearly, $r(F)+2c(F)=\deg
F$. The following theorem is due to Dirichlet (see, \eg,
\cite{Borevich}).

\theorem
The rank of the group of units \rom(\ie, invertible elements of the ring
of integers\rom) of an algebraic number field~$F$ is $r(F)+c(F)-1$.
\qed
\endtheorem

Let $n=\ord\Gt$ and assume that $\Gf(n)\ge4$. Let~$\k$, $\gD$,
and~$M_\Gr\sppm$ be as in~\ref{AutM}. Note that
$r(\k)=\deg\k=\frac12\Gf(n)\ge2$ and $c(\k)=0$.

\lemma\label{2-compact}
If $\Gk\ne 1$, $\Gf(n)\ge4$, and $\dim_\k V\sppm_\Gr=2$, then the
rational Dirichlet polyhedra of $\Aut M\sppm_\Gr$ in $\HS\sppm_\Gr$ are
compact.
\endlemma

\proof
Since $\HS\sppm_\Gr$ are hyperbolic lines, it suffices to show that the
groups $\Aut M\sppm_\Gr$ are infinite. Consider one of them, say, $\Aut
M\spp_\Gr$. The lattice $M\spp_\Gr$ contains a finite index
sublattice~$M'$ whose Gramm matrix (after, possibly, dividing the form by
an element of~$\gD$) is of the form
$$
\bmatrix0&1\\1&0\endbmatrix\qquad\text{or}\qquad
\bmatrix1&0\\0&-d\endbmatrix\quad\text{with $d>0$ and
$\sqrt{d}\notin\k$}.
$$
In the former case (which occurs if the form represents~$0$ over~$\k$),
the automorphisms of~$M'$ are of the form
$$
A_\Gl=\bmatrix\pm\Gl&0\\0&\pm1/\Gl\endbmatrix,
$$
where $\Gl\in\gD^*$ is a unit of~$\k$. Thus, in this case $\Aut
M\spp_\Gr$ contains a free abelian group of rank $r(\k)-1>0$.

In the latter case, the automorphisms of~$M'$ are of the form
$$
\bmatrix\pm1&0\\0&\pm1\endbmatrix\qquad\text{or}\qquad
B_\Gl=\bmatrix\Ga&d\Gb\\\Gb&\Ga\endbmatrix,
$$
where $\Ga,\Gb\in\gD$ and $\Gl=\Ga+\Gb\sqrt{d}$ is a unit of
$F=\k(\sqrt{d})$ such that $\Ga^2-\Gb^2d=1$. We will show that the group
of such units is at least~$\Z$.

The map $\mu\:\Ga+\Gb\sqrt{d}\mapsto\Ga^2-\Gb^2d$ is a homomorphism from
the group of units of~$F$ to the group of units of~$\k$, and its cokernel
is finite. As $d>0$, the quadratic extension~$F$ of~$\k$ has at least two
real embeddings to~$\C$, \ie, $r(F)\ge2$. Since
$r(F)+2c(F)=2\deg\k=2r(\k)$, one has\footnote{In fact, under the
assumption on the signature of the form, $F$ has exactly two real
embeddings to~$\C$, namely, $\k(\sqrt{d})$ and $\k(-\sqrt{d})$. In
particular, modulo torsion one has $\Ker\mu=\Z$. Indeed, the other
embeddings are $\k(\pm\sqrt{\gg(d)})$, $\gg\in\Gal(\k/\Q)$, $\gg\ne1$,
and since all spaces $L_{\gg\Gr}(\k)$ are negative definite, one has
$\gg(d)<0$.} $\rank\Ker\mu=\frac12r(F)\ge1$.

The coefficients~$\Ga$, $\Gb$ of all integers $\Ga+\Gb\sqrt{d}$ of~$F$
have `bounded denominators', \ie, $\Ga,\Gb\in\frac1N\gD$ for some
$N\in\N$ (since the abelian group generated by~$\Ga$'s and~$\Gb$'s has
finite rank and $\gD$ has maximal rank). Hence, for any $\Gl\in\Ker\mu$,
the map $B_\Gl$ defines an isometry of~$V\spp$ taking $N\cdot M'$
into~$M'$, and Lemma~\ref{rational->integer} applies.
\endproof

\remark{Remark}
Note that, if $\Gf(n)>4$, the form cannot represent~$0$ over~$\k$.
Indeed, otherwise $\Aut M_\Gr\sppm$ would contain a free abelian group of
rank~$\ge2$, which would contradict to the discreteness of the action.
\endremark

Next theorem (as well as Lemma~\ref{2-compact}) can probably be deduced
from the Godement criterion. We chose to give here an alternative
self-contained proof.

\theorem\label{compact}
If $\Gk\ne 1$ and $\Gf(n)\ge4$, then the rational Dirichlet polyhedra of
$\Aut M\sppm_\Gr$ in $\HS\sppm_\Gr$ are compact.
\endtheorem

\proof
Let $m=\dim_\k V\sppm_\Gr$. The assertion is obvious if $m=1$, and it is
the statement of~\ref{2-compact} if $m=2$. If $m>2$ and a rational
Dirichlet polyhedron $\Pi\subset\HS\spp_\Gr$ is not compact, one can find
a line $\HS'=\HS(V'\otimes_\k\R)$, $V'\subset V\spp_\xi$, such that
$\Pi\cap\HS'$ is not compact. (If $\Pi=\HS\spp_\Gr$, one can take
for~$V'$ any hyperbolic $2$-subspace. Otherwise, one can replace~$\Pi$
with one of its non-compact facets and proceed by induction.)
Applying~\ref{2-compact} to $M'=V'\cap L_\rho(\gD)$, one concludes that
the polyhedron $\Pi'\subset\HS'$ of~$\Aut M'$ is compact. On the other
hand, in view of~\ref{extension}, $\Pi\cap\HS'$ must be a finite union of
copies of~$\Pi'$.
\endproof

\corollary\label{finite.copy>2}
Assume that $\Gk\ne 1$ and $\Gf(\ord\Gt)\ge 4$, and let $\Pi_\Gr\sppm$ be
some rational Dirichlet polyhedra of $\Aut M_\Gr\sppm$ in~$\HS_\Gr\sppm$.
Then for any integer~$a$ there are at most finitely many vectors $v\in
M\sppm$ such that $v^2=a$ and the subspace $\psh_\Gr\sppm(v)$
intersects~$\Pi_\Gr\sppm$ or~$J_\Gr(\Pi_\Gr\sp\mp)$.
\qed
\endcorollary

\subsection{Proof of Theorem~\ref{finiteness}}\label{proof}
In view of~\ref{finite.copy}, one can replace~$\GS_\Gr\sppm$ in the
statement with the rational Dirichlet polyhedra~$\Pi_\Gr\sppm$ of $\Aut
M\sppm$ in~$\HS_\Gr\sppm$.

For a root $v\in(L^G)\ort$ denote by~$v\sppm$ its projections to~$V\sppm$
(the $(\pm1)$-eigenspaces of~$c$ on $L_\Gr(\Q)\otimes\R$\,) and
by~$v_\Gr\sppm$, its projections to~$V_\Gr\sppm$. The projections
$v\sppm$ are rational vectors with uniformly bounded denominators, \ie,
there is an integer~$N$, depending only on~$\Gt$, such that $N v\sppm\in
M\sppm$. Under the assumption ($\Gr$ is non-real and $\Gk\ne1$), the set
$\psh_\Gr\spp(v)$ is not empty if and only if each of~$v_\Gr\sppm$ either
is trivial or has negative square. In any case, $(v_\Gr\sppm)^2\le0$ and,
hence, $(v\sppm)^2\le0$. Thus, the squares $(N v\sppm)^2$ take finitely
many distinct integral values, and the statement of the theorem follows
from~\ref{Siegel.M} and~\ref{finite.copy>2}.
\qed

\section{The proof}\label{The.proof}

\subsection{Period spaces related to $K3$-surfaces}\label{period.spaces}
Let $L=2E_8\oplus3U$. Consider the variety~$\Per$ of positive definite
$3$-subspaces in $L\otimes\R$. It is a homogeneous symmetric space (of
noncompact type):
$$
\Per=\grSO\spp(3,19)/\grSO(3)\times\grSO(19).
$$
The orthogonal projection of a positive definite $3$-subspace to another
one is nondegenerate. Hence, one can orient all the subspaces in a
coherent way; this gives an orientation of the canonical $3$-dimensional
vector bundle over~$\Per$. In what follows we assume such an orientation
fixed; the corresponding orientation of a space $\gw\in\Per$ is referred
to as its \emph{prescribed orientation}.

Given a vector $v\in L$  with $v^2=-2$, let $\hp{v}\subset\Per$ be the
set of the $3$-subspaces orthogonal to~$v$. Put
$$
 \tsize
 \Per_0=\Per\sminus\bigcup_{v\in L,\,v^2=-2}\hp{v}.
$$
The space $\Per_0$ is called the \emph{period space of marked Einstein
$K3$-surfaces}.

There is a natural $S^2$-bundle $K\GO\to\Per$, where
$$
K\GO=\{(\gw,\Gg)\mid\gw\in\Per,\;\Gg\in\gw,\;\Gg^2=1\}.
$$
The pull-back $K\GO_0$ of $\Per_0$ is called the \emph{period space of
marked K\"ahler $K3$-surfaces}. Finally, let $\GO$ be the variety of
oriented positive definite $2$-subspaces of~$L\otimes\R$; it is called
the \emph{period space of marked $K3$-surfaces}. One can identify~$\GO$
with the projectivization
$$
 \{\Go\in L\otimes\C\mid\Go^2=0,\;\Go\cdot\bar{\Go}>0\}/\C^*,
 \eqtag\label{GO-C}
$$
associating to a complex line generated by~$\Go$ the plane
$\{\Re(\Gl\Go)\mid\Gl\in\C\}$ with the orientation given by a basis
$\Re\Go$, $\Re i\Go$. Thus, $\GO$ is a $20$-dimensional complex variety,
which is an open subset of the quadric defined in the projectivization of
$L\otimes\C$ by $\Go^2=0$. The spaces~$K\GO_0$ and~$\Per_0$ are
(noncompact) real analytic varieties of dimensions~$59$ and~$57$,
respectively.

\subsection{Period maps}\label{period.maps}
A \emph{marking} of a $K3$-surface~$X$ is an isometry
$\Gf\:H^2(X)\to L$. It is called \emph{admissible} if the orientation of
the space $\gw=\<\Re\Gf(\Go),\Im\Gf(\Go),\Gf(\Gg)\>$, where $\Go\in
H^{2,0}(X)$ and $\Gg$ is the fundamental class of a K\"ahler structure
on~$X$, coincides with its prescribed orientation. A \emph{marked
$K3$-surface} is a $K3$-surface~$X$ equipped with an admissible marking.
Two marked $K3$-surfaces $(X,\Gf)$ and $(Y,\psi)$ are \emph{isomorphic}
if there exists a biholomorphism $f\:X\to Y$ such that $\psi=\Gf\circ
f^*$. Denote by $\CT$ the set of isomorphism classes of marked
$K3$-surfaces.

The \emph{period map} $\per\:\CT\to\GO$ sends a marked $K3$-surface
$(X,\Gf)$ to the $2$-subspace $\{\Re\Gf(\Go)\mid\Go\in H^{2,0}(X)\}$, the
orientation given by $(\Re\Gf(\Go),\Re\Gf(i\Go))$. (We will always use
the same notation~$\Gf$ for various extensions of the marking to other
coefficient groups.) Alternatively, $\per(X,\Gf)$ is the line
$\Gf(H^{2,0}(X))$ in the complex model~\eqref{GO-C} of $\GO$.

A \emph{marked polarized $K3$-surface} is a $K3$-surface~$X$ equipped
with the fundamental class $\Gg_X$ of a K\"ahler structure and an
admissible marking $\Gf\:H^2(X)\to L$. Two marked polarized $K3$-surfaces
$(X,\Gf,\Gg_X)$ and $(Y,\psi,\Gg_Y)$ are \emph{isomorphic} if there
exists a biholomorphism $f\:X\to Y$ such that $\psi=\Gf\circ f^*$ and
$f^*(\Gg_Y)=\Gg_X$. Denote by $K\CT$ the set of isomorphism classes of
marked polarized $K3$-surfaces.

The \emph{period map} $\per^K\:K\CT\to K\GO$ sends a triple
$(X,\Gf,\Gg_X) \in K\CT$ to the point $(\gw,\Gf(\Gg_X))\in K\GO$, where
$\gw=\per(X,\Gf)\oplus\Gf(\Gg_X)\in\Per$ is as above. When this does not
lead to a confusion, we abbreviate $\per^K(X,\Gf,\Gg_X)$ to $\per^K(X)$.

As is known (see~\cite{PSh-Sh} and~\cite{Kulikov}, or~\cite{Siu}), the
period map $\per^K$ is a bijection to~$K\GO_0$, and the image of~$\per$
is $\GO_0$. Moreover, $K\GO_0$ is a fine period space of marked polarized
$K3$-surfaces, \ie, the following statement holds (see~\cite{Beauville}).

\theorem\label{Univer-family}
The space $K\GO_0$ is the base of a universal smooth family of marked
polarized $K3$-sur\-faces, \ie, a family $p\:\Phi\to K\GO_0$ such that
any other smooth family $p'\:X\to S$ of marked polarized $K3$-surfaces is
induced from~$p$ by a unique smooth map $S\to K\GO_0$. The latter is
given by $s\mapsto\per^K(X_s)$, where $X_s$ is the fiber over $s\in S$.
\endtheorem

Since the only automorphism of a $K3$-surface identical on the homology
is the identity (see~\cite{PSh-Sh}), Theorem~\ref{Univer-family} can be
rewritten in a slightly stronger form.

\theorem\label{Univer-family.1}
For any smooth family $p'\:X\to S$ of marked polarized $K3$-sur\-faces
there is a unique smooth fiberwise map $X\to\Phi$
\rom(see~\ref{Univer-family}\rom) that covers the map $S\to K\GO_0$,
$s\mapsto\per^K(X_s)$ of the bases and is an isomorphism of marked
polarized $K3$-surfaces in each fiber.
\endtheorem

\corollary\label{CoHoMap}
Let $(X,\Gg_X)$ and $(Y,\Gg_Y)$ be two polarized $K3$-surfaces and let
$g\:H^2(Y)\to H^2(X)$ be an isometry such that $g(\Gg_Y)=\Gg_X$.
Then\rom:
\roster
\item
if $g(H^{2,0}(Y))=H^{2,0}(X)$, then $g$ is induced by a unique
holomorphic map $X\to Y$, which is a biholomorphism\rom;
\item
if $g(H^{2,0}(Y))= H^{0,2}(X)$, then $-g$ is induced by a unique
anti-holomorphic map $X\to Y$, which is an anti-biholomorphism.
\qed
\endroster
\endcorollary

\subsection{Equivariant period spaces}
In this section we construct the period space of marked polarized
$K3$-surfaces with a $G$-action of a given homological type. Recall that
we define the homological type as the class of the twisted induced action
$\Gt_X\:G\to\Aut H^2(X)$ modulo conjugation by elements of $\Aut H^2(X)$.
A marking takes~$\Gt_X$ to an action $\Gt\:G\to\Aut L$. Note in this
respect that, since we work with admissible markings only, it would be
more natural to consider~$\Gt_X$ up to conjugation by elements of the
subgroup $\Aut L\cap\grO\spp(L\otimes\R)$. However, this stricter
definition would be equivalent to the original one, as the central
element $-\id\in\Aut L$ belongs to $\grO\spm(L\otimes\R)$.

\proposition\label{K3->geometric}
Let~$X$ be a $K3$-surface supplied with a Klein action of a finite
group~$G$. Then the twisted induced action $\Gt_X\:G\to\Aut H^2(X)$ is
geometric, and the augmentation $\Gk\:G\to\{\pm1\}$ and the pair $\Gr,
\bar\Gr \:\Gsbp\to S^1$ of complex conjugated fundamental representations
introduced in~\ref{results} coincide with those determined by~$\Gt_X$
\rom(see~\ref{rho}\rom).
\endproposition

\proof
Since $G$ is finite, $X$ admits a K\"ahler metric preserved by the
holomorphic elements of~$G$ and conjugated by the anti-holomorphic
elements. Take for~$\Gg_X$ the fundamental class of such a metric. Pick
also a holomorphic form on~$X$ and denote by~$\Go$ its cohomology class.
Let $\gw$ be the space spanned by $\Gg_X$, $\Re\Go$, and $\Im\Go$, and
let $\ell \subset \gw$ be the subspace generated by~$\Gg_X$. Then the
flag $\ell \subset \gw$ attests the fact that $\Gt_X$ is almost
geometric, and this flag can be used to define~$\Gk$ and~$\Gr$. As
$\Gg_X$ and~$\Go$ cannot be simultaneously orthogonal to an integral
vector $v\in H^2(X)$ of square~$(-2)$, the action is geometric.
\endproof

Let $\Gt\:G\to\Aut L$ be an almost geometric action on~$L$. The
assignment $g\:\gw\mapsto \Gk(g)g(\gw)$, where $g\in G$ and $-\gw$ stands
for $\gw$ with the opposite orientation, defines a $G$-action on the
space $\Per$. Denote by $\Per\G$ the subspace of the $G$-fixed points and
let $\Per\G_0=\Per\G\cap\Per_0$. There is a natural map
$K\GO\G\to\Per\G$, where
$$
K\GO\G=\{(\gw,\Gg)\mid\gw\in\Per\G,\;\Gg\in\gw\G,\;\Gg^2=1\},
$$
$\gw\G$ standing for the $G$-invariant part of~$\gw$. Put
$K\GO\G_0=\{(\gw,\Gg)\in K\GO\G\mid\gw\in\Per\G_0\}$ and denote by
$\GO\G$ (respectively, $\GO\G_0$) the image of $K\GO\G$ (respectively,
$K\GO\G_0$) under the projection $K\GO\to\GO$. The following statement is
a paraphrase of the definitions.

\proposition\label{non-empty}
An almost geometric action $\Gt: G \to \Aut L$ is geometric if and only
if the space $K\GO\G_0$ \rom(as well as $\Per\G_0$ and $\GO\G_0$\rom) is
non-empty.
\qed
\endproposition

Let $(X, \Gf)$ be a marked $K3$-surface. We will say that a Klein
$G$-action on $X$ and an action $\Gt: G \to \Aut L$ are \emph{compatible}
if for any $g\in G$ one has $\Gt_Xg=\Gf^{-1}\circ\Gt g\circ\Gf$, where
$\Gt_X\:G\to\Aut H^2(X)$ is the twisted induced action. If a marking is
not fixed, we say that a Klein $G$-action on~$X$ is compatible with~$\Gt$
if $X$ admits a compatible admissible marking, \ie, if $\Gt_X$ is
isomorphic to~$\Gt$.

\proposition\label{G-periods}
An action $\Gt\:G\to L$ is compatible with a Klein $G$-action on a marked
$K3$-surface if and only if $\Gt$~is geometric. Furthermore, $K\GO\G_0$
is a fine period space of marked polarized $K3$-surfaces with a Klein
$G$-action compatible with~$\Gt$, \ie, it is the base of a universal
smooth family of marked polarized $K3$-surfaces with a Klein $G$-action
compatible with~$\Gt$.
\endproposition

\proof
The `only if' part follows from~\ref{K3->geometric}, and the `if',
from~\ref{CoHoMap} and~\ref{non-empty}. The fact that $K\GO\G_0$ is a
fine period space is an immediate consequence of \ref{Univer-family.1}.
\endproof

\proposition
Let $\Gk\:G\to\{\pm1\}$ be the augmentation and $\Gr\:\Gsbp\to S^1$ a
fundamental representation associated with~$\Gt$. If $\Gr=1$, then the
spaces $K\GO\G$ and $\GO\G$ are connected. If $\Gr\ne1$, then the space
$K\GO\G$ \rom(respectively, $\GO\G$\rom) consists of two components,
which are transposed by the involution $(\gw,\Gg)\mapsto(\gw,-\Gg)$
\rom(respectively, the involution reversing the orientation of
$2$-subspaces\rom). If, besides, $\Gr\ne\bar\Gr$, the two components of
$K\GO\G$ \rom(or $\GO\G$\rom) are in a one-to-one correspondence with the
two fundamental representations~$\Gr$, $\bar\Gr$.
\endproposition

\proof
Since $\Per$ is a hyperbolic space and $G$ acts on~$\Per$ by isometries,
the space $\Per^G$ is contractible. The projections $K\GO\G\to\Per\G$ and
$K\GO\G_0\to\Per\G_0$ are (trivial) $S^p$-bundles, where $p=0$ if
$\Gr\ne1$, $p=1$ if $\Gr=1$ and $\Gk\ne1$, and $p=2$ if $\Gr=1$ and
$\Gk=1$. Finally, since each space $\gw\in\Per$ has its prescribed
orientation, a choice of a $G$-invariant vector $\Gg\in\gw$ determines an
orientation of $\Gg\ort\subset\gw$ and, hence, a fundamental
representation.
\endproof

\subsection{The moduli spaces}\label{moduli-space}
Fix a geometric action $\Gt: G \to \Aut L$ and consider the space
$\KMS\G=K\GO\G_0/\!\Aut_GL$. In view of~\ref{G-periods}, it is the
`moduli space' of polarized $K3$-surfaces with Klein $G$-actions
compatible with~$\Gt$. Given such a surface~$(X,\Gg_X)$, pick a marking
$\Gf\:H^2(X)\to L$ compatible with~$\Gt$ and denote by
$\kms(X,\Gg_X)=\kms(X)$ the image of $\per^K(X,\Gf,\Gg_X)$ in $\KMS\G$.
Since any two compatible markings differ by an element of $\Aut_GL$, the
point $\kms(X,\Gg_X)$ is well defined. The following statement is an
immediate consequence of~\ref{G-periods} and the local connectedness
of~$K\GO\G_0$.

\proposition
Let~$(X,\Gg_X)$ and~$(Y,\Gg_Y)$ be two polarized $K3$-surfaces with Klein
$G$-actions compatible with~$\Gt$. Then $X$ and $Y$ are $G$-equivariantly
deformation equivalent if and only if $\kms(X)$ and $\kms(Y)$ belong to
the same connected component of $\KMS\G$.
\qed
\endproposition

In \ref{case-1h}--\ref{case-cah} below we give a more detailed
description of period and moduli spaces. We use the notations
of~\ref{geometric-actions}.

\theorem[The case $\Gr=1$, $\Gk=1$]\label{case-1h}
If $\Gr=1$ and $\Gk=1$, then $K\GO\G_0\cong(\HS(L^G)\sminus\GD)\times
S^2$, where $\codim\GD\ge3$. In particular, $K\GO\G_0$ and, hence,
$K\MS\G$ are connected.
\endtheorem

\theorem[The case $\Gr=1$, $\Gk\ne1$]\label{case-1ah}
If $\Gr=1$ and $\Gk\ne1$, then $K\MS\G$ is a quotient of the connected
space $((\HS(L^G)\times\Int\GG_\Gr\spm)\sminus\GD)\times S^1$, where
$\codim\GD\ge2$. In particular, $K\MS\G$ is connected.
\endtheorem

\theorem[The case $\Gr\ne1$ real, $\Gk=1$]\label{case-rh}
If $\Gr\ne1$ is real and $\Gk=1$, then $K\MS\G$ is a quotient of the
two-component space $((\Int\GG_1\times\HS(L_\Gr(\R)))\sminus\GD)\times
S^0$, where $\codim\GD\ge2$. In particular, $K\MS\G$ has at most two
connected components, which are interchanged by the complex conjugation
$X\mapsto\wb X$.
\endtheorem

\theorem[The case $\Gr\ne1$ real, $\Gk\ne1$]\label{case-rah}
If $\Gr\ne1$ is real and $\Gk\ne1$, then $K\MS\G$ is a quotient of the
two-component space
$((\Int\GG_1\times\Int\GG_\Gr\spp\times\Int\GG_\Gr\spm)\sminus\GD)\times
S^0$, where $\codim\GD\ge2$. In particular, $K\MS\G$ has at most two
connected components, which are interchanged by the complex conjugation
$X\mapsto\wb X$.
\endtheorem

\theorem[The case $\Gr$ non-real, $\Gk=1$]\label{case-ch}
If $\Gr$ is non-real and $\Gk=1$, then $K\MS\G$ is a quotient of the
two-component space $((\Int\GG_1\times\P\!_J\Cal C_\Gr)\sminus\GD)\times
S^0$, where $\P\!_J\Cal C_\Gr$ is the space of positive definite
\rom(over~$\R$\rom) $J_\Gr$-complex lines in $L_\Gr(\R)$ and
$\codim\GD\ge2$. In particular, $K\MS\G$ has at most two connected
components, which are interchanged by the complex conjugation
$X\mapsto\wb X$.
\endtheorem

\theorem[The case $\Gr$ non-real, $\Gk\ne1$]\label{case-cah}
If $\Gr$ is non-real and $\Gk\ne1$, then $K\MS\G$ is a quotient of the
space $((\Int\GG_1\times\GS_\Gr\spp)\sminus\GD)\times S^0$, where $\GD$
is the union of a subset of codimension~$\ge2$ and finitely many
hyperplanes of the form
$\Int\GG_1\times(\psh_\Gr\sppm(v)\cap\GS_\Gr\spp)$ defined by roots
\smash{$v\in(L^G)\ort$}. This space has finitely many connected
components\rom; hence, so does $K\MS\G$.
\endtheorem

\proof[Proof of~\ref{case-1h}--\ref{case-rah}]
One has
\roster\widestnumber\item{}
\item"--"
$\Per\G=\HS(L^G\otimes\R)$ in case~\ref{case-1h},
\item"--"
$\Per\G=\HS(L^G\otimes\R)\times\HS(V_\Gr\spm)$ in case~\ref{case-1ah},
\item"--"
$\Per\G=\HS(L^G\otimes\R)\times\HS(L_\Gr(\R))$ in case~\ref{case-rh}, and
\item"--"
$\Per\G=\HS(L^G\otimes\R)\times\HS(V_\Gr\spp)\times\HS(V_\Gr\spm)$ in
case~\ref{case-rah}.
\endroster
Thus, in each case, $\Per\G$ is a product $\prod\HS(L_i\otimes\R)$ of the
hyperbolic spaces of orthogonal indefinite sublattices $L_i\subset L$
such that $\bigoplus_i L_i \oplus \eL$ is a finite index sublattice
in~$L$. Consider the quotient $\CQ_0=\Per\G_0/W$, where $W=\prod W_i$
(the product in $W_G(L)$) and $W_i=1$ if $\Gs\sbp L_i>1$ or
$W_i=W_G((L_i\oplus\eL)\prhull)$ if $\Gs\sbp L_i=1$. The quotient~$\CQ_0$
can be identified with a subspace of $\CQ=\prod\Int\GG_i$, where $\GG_i$
is a fundamental Dirichlet polyhedron of~$W_i$ in $\HS(L_i\otimes\R)$.
(Note that $\GG_i=\HS(L_i\otimes\R)$ unless $\Gs\sbp L_1=1$.) Put
$\GD=\CQ\sminus\CQ_0$; it is the union of the walls $\hp{v}\cap\CQ$ over
all roots $v\in L$.

For a root $v\in L$ one has $\codim(\hp{v}\cap\CQ)\ge\sum\Gs\sbp L_i$,
the summation over all~$i$ such that the projection of~$v$ to~$L_i$ is
nontrivial. Thus, a wall $\hp{v}\cap\CQ$ may have codimension~$1$ only if
$v\in(L_i\oplus\eL)\prhull$ and $\Gs\sbp L_i=1$. However, in this case
$\hp{v}\cap\CQ=\varnothing$ due to~\ref{connectedness}. Hence,
$\codim\GD\ge2$ and the space $\CQ_0$ is connected.
\endproof

\proof[Proof of~\ref{case-ch}]
In this case, $\Per\G_0/W_G((L^G\oplus\eL)\prhull)$ can be identified
with a subset of $\Int\GG_1\times\P\!_J\Cal C_\Gr$, and the proof follows
the lines of the previous one.
\endproof

\proof[Proof of~\ref{case-cah}]
One has $\Per\G = \HS(L^G\otimes\R)\times\HS(V_\Gr\spp)$, and the
quotient space
$\CQ_0=\Per\G_0/(W_G(L_\Gr(\Z)\ort)\cdot\Aut_G^0(L_\Gr(\Z))$ can be
identified with a subset of $\Int\GG_1\times\GS_\Gr\spp$. Now, the
statement follows from~\ref{connectedness} and~\ref{finiteness}.
\endproof

\subsection{Proof of Theorems~\ref{Main-Finiteness} and~\ref{Main-q-simplicity}}\label{main.proof}
Theorem\ref{Main-q-simplicity} follows from~\ref{case-1h}--\ref{case-ch}.
Theorem~\ref{Main-Finiteness} consists, in fact, of two statements:
finiteness of the number of equivariant deformation classes within a
given homological type of $G$-actions (of a given group~$G$), and
finiteness of the number of homological types of faithful actions. The
former is a direct consequence of~\ref{case-1h}--\ref{case-cah}. The
latter is a special case of the finiteness of the number of conjugacy
classes of finite subgroups in an arithmetic group, see~\cite{BH}
and~\cite{B}.
\qed

\section{Degenerations}\label{Degenerations}

\subsection{Passing through the walls}\label{walls}
Let $L=2E_8\oplus3U$. Consider a geometric $G$-action $\Gt\:G\to\Aut L$.
Pick a $G$-invariant elliptic root system~$R\subset L$. Denote by $\bar
R$ the sublattice of~$L$ generated by all roots in $(R+\eL)\prhull$.
Clearly, $\bar R$ is a $G$-invariant root system; it is called the
\emph{$\Gt$-saturation} of~$R$. We say that~$R$ is \emph{$\Gt$-saturated}
if $R=\bar R$. Any $\Gt$-saturated root system~$R$ is \emph{saturated},
\ie, $R$ contains all roots in~$R\prhull$.

Fix a camera~$C$ of~$\bar R$ and denote by~$S_C$ its group of symmetries.
Then, for any $g\in G$, the restriction of~$\Gt g$ to~$\bar R$ admits a
unique decomposition $s_gw_g$, $s_g\in S_C$, $w_g\in W(\bar R)$. Let
$\Gt_R(g)=(\Gt g)w_g^{-1}\in\Aut L$, $w_g$ being extended to~$L$
identically on~\smash{$\bar R\ort$}. We will call the map
$\Gt_R\:G\to\Aut L$ the \emph{degeneration} of~$\Gt$ at~$R$.

\proposition\label{degeneration}
The map\/ $\Gt_R$ is a geometric $G$-action. Up to conjugation by an
element of~$W(\bar R)$, it does not depend on the choice of a camera~$C$
of~${\bar R}$ and is the only action with the following properties\rom:
\roster
\item\local1
the action induced by~$\Gt_R$ on~${\bar R}$ is admissible\rom;
\item\local2
$\Gt$ and $\Gt_R$ induce the same action on each of the following
sets\rom: $\bar R\ort$, $\discr\bar R$, the set of irreducible components
of~$\bar R$.
\endroster
Conversely, if $\bar R\subset L$ is a saturated root system and
$\Gt_R\:G\to L$ is an action satisfying~\loccit1--\loccit2 above, then
$\bar R$ is $\Gt$-saturated and $\Gt_R$ is a degeneration of~$\Gt$
at~$\bar R$.
\endproposition

\proof
Clearly, both $\Gt$ and $\Gt_R$ factor through a subgroup of $\Aut\bar
R\times\Aut\bar R\ort$. The composition of~$\Gt_R$ with the projection to
$\Aut\bar R\ort$ coincides with that of~$\Gt$; the composition of $\Gt_R$
with the projection to $\Aut \bar R$ is the composition of~$\Gt$, the
projection to $\Aut\bar R$, and the quotient homomorphism $\Aut {\bar
R}\to S_C\subset\Aut \bar R$. Hence, $\Gt_R$ is a homomorphism.
Furthermore, another choice of a camera~$C'$ of~$\bar R$ leads to another
representation $\Aut {\bar R}\to S_{C'}\subset\Aut \bar R$, which is
conjugated to the original one by a unique element $w_0\in W(\bar R)$;
the latter can be regarded as an automorphism of~$L$.

All other statements follow directly from the construction. For the
uniqueness, it suffices to notice that, for any irreducible root
system~$R'$ and a camera~$C'$ of~$R'$, the natural homomorphism
$S_{C'}\to\Aut\discr R'$ is a monomorphism.
\endproof

\proposition\label{system-in-Lg}
Let $R$ be a $\Gt$-saturated root system and $R'\subset R$ the sublattice
generated by all roots in $R\cap(L^G)\ort$. Then, up to conjugation by an
element of $W(R)$, the degenerations~$\Gt_R$ and~$\Gt_{R'}$ coincide. In
particular, $\Gt_R$ can be chosen to coincide with~$\Gt$ on $(R')\ort$.
\endproposition

\proof
Take for~$C$ a camera adjacent to the intersection of the mirrors defined
by the roots of~$R'$. Then $C$ has an invariant face (possibly, empty),
and the decomposition $\Gt g|_R=s_gw_g$ has $w_g\in W(R')$ for any $g\in
G$.
\endproof

If the action is properly Klein, one can take for $R$ the $\Gt$-saturated
root system generated by all roots in~$(L^G)\ort$ orthogonal to a given
wall $\psh_\Gr\spp(v)$. The resulting degeneration is called the
\emph{degeneration} at the wall $\psh_\Gr\spp(v)$.

\remark{Remark}
The degeneration construction gives rise to a partial order on the set of
homological types of geometric actions of a given finite group~$G$.
\endremark

\subsection{Degenerations of $K3$-surfaces}\label{degenerations}
Let $(G,\Gk)$ be an augmented group. Denote by $D_\Ge$ the disk
$\{s\in\C\mid\ls|s|<\Ge\}$. The composition of~$\Gk$ and the
$\{\pm1\}$-action via the complex conjugation $s\mapsto\bar s$ is a Klein
$G$-action on~$D_\Ge$. A \emph{$G$-equivariant degeneration of
$K3$-surfaces} is a nonsingular complex $3$-manifold~$X$ supplied with a
Klein $G$-action and a $G$-equivariant (with respect to the above
$G$-action on~$D_\Ge$) proper analytic map $p\:X\to D_\Ge$ so that the
following holds:
\roster\widestnumber\item{}
\item"--"
the projection~$p$ has no critical values except $s=0$;
\item"--"
the fibers $X_s=p^{-1}(s)$ of~$p$ are normal $K3$-surfaces, nonsingular
unless $s=0$.
\endroster
(By a singular $K3$-surface we mean a surface whose desingularization
is~$K3$.) Given a degeneration~$X$, denote by $\pi_s\:\tilde X_s\to X_s$,
$s\in D_\Ge$, the minimal resolution of singularities of~$X_s$, see,
\eg,~\cite{La}. (Note that $\tilde X_s=X_s$ unless $s=0$.) From the
uniqueness of the minimal resolution it follows that any Klein action
lifts from $X_s$ to $\tilde X_s$. Thus, if either $\Gk=1$ or $s$ is real,
$\tilde X_s$ inherits a natural Klein action of~$G$.

\theorem\label{family}
Let $p\:X\to D_\Ge$ be a $G$-equivariant degeneration of $K3$-surfaces.
Pick a regular value $t\in D_\Ge$, real, if $\Gk\ne1$. Denote by
$R\subset H^2(X_t)$ the subgroup Poincar\'e dual to the kernel of the
inclusion homomorphism $H_2(X_t)\to H_2(X)=H_2(X_0)$. Then $R$ is a
saturated elliptic root system and the twisted induced $G$-action on
\smash{$H^2(\tilde X_0)$} is isomorphic to the degeneration at~$R$ of the
twisted induced $G$-action on $H^2(X_t)$.
\endtheorem

\remark{Remark}
A statement analogous to Theorem~\ref{family} holds in a more general
situation, for a family of complex surfaces whose singular fiber at $s=0$
has at worst simple singularities, \ie, those of type $A_n$, $D_n$,
$E_6$, $E_7$, or~$E_8$. The only difference is the fact that one can no
longer claim that the root system~$R$ is saturated, and one should
consider the degeneration at~$R$ {\bf without} passing to its saturation
first. (In particular, the algebraic definition of degeneration should be
changed. Our choice of the definition, incorporating the saturation
operation, was motivated by our desire to assure that the result should
be a geometric action.) The proof given below applies to the general case
with obvious minor modifications.
\endremark

\proof
It is more convenient to switch to the twisted induced actions~$\Gt_s$ in
the homology groups $H_2(X_s)$, $s\in D_\Ge$; they are Poincar\'e dual to
the twisted induced actions in the cohomology.

Let $\Gi_s\:H_2(\tilde X_s)\to H_2(X)$, $s\in D_\Ge$, be the composition
of~$(\pi_s)_*$ and the inclusion homomorphism. Put $R_s=\Ker\Gi_s$.
Consider sufficiently small $G$-invariant open balls $B_i\subset X$ about
the singular points of~$X_0$ and let $B=\bigcup B_i$. One can assume that
$t$ is real and sufficiently small, so that $M_i=X_t\cap B_i$ are Milnor
fibers of the singular points. Then there is a $G$-equivariant
diffeomorphism $d'\:X_t\sminus B\to X_0\sminus B$.

Recall that all singular points of the $K3$-surface~$X_0$ are simple and
$R_0$ is a saturated elliptic root system (see Lemma~\ref{K3-saturated}
below). In particular, $d'$ extends to a diffeomorphism $d\:X_t\to\tilde
X_0$. Note that neither~$d$ nor the induced isomorphism $d_*\:H_2(X_t)\to
H_2(\tilde X_0)$ is canonical and $d_*$ does not need to be
$G$-equivariant. However, $d_*$ does preserve the $G$-action on the sets
of irreducible components of the root systems~$R_t$, $R_0$ (as it is just
the $G$-action on the set of singular points of~$X_0$), and, in view of
natural identifications $R_s\ort=H_2(X_s\sminus B)/\Tors$ and $\discr
R_s=H_1(\partial(X_s\sminus B))$, $s=t,0$, and the fact that $d'$
commutes with~$G$, the restrictions of~$d_*$ to~$R_t\ort$ and $\discr
R_t$ are $G$-equivariant. Finally, the action induced by~$\Gt_0$ on~$R_0$
is admissible: it preserves the camera defined by the exceptional
divisors in~$\tilde X_0$ (see~\ref{admissible}). Thus, after identifying
$H_2(X_t)$ and $H_2(\tilde X_0)$ via~$d_*$, the actions $\Gt=\Gt_t$ and
$\Gt_R=\Gt_0$ satisfy \iref{degeneration}1--\ditto2, and
\ref{degeneration} implies that $\Gt_0$ is the degeneration of~$\Gt_t$ at
$R_t$.
\endproof

For completeness, we outline the proof of the following lemma, which
refines the well known fact that a $K3$-surface can have at worst simple
singular points.

\lemma\label{K3-saturated}
Let $X$ be a $K3$-surface. Then any negative definite sublattice
$R\subset H^2(X)$ generated by classes of irreducible curves is a
saturated root system.
\endlemma

\proof
As it follows from the adjunction formula, any irreducible curve
$C\subset X$ of negative self-intersection is a $(-2)$-curve, \ie, a
non-singular rational curve of self-intersection~$(-2)$. Thus, any
sublattice~$R$ as in the statement is an elliptic root system generated
by classes of irreducible $(-2)$-curves.

 From the Riemann-Roch theorem it follows that, given a root $r\in\Pic X$,
there is a unique $(-2)$-curve $C\subset X$ whose cohomology class is
$\pm r$. Thus, the set of all roots in $\Pic X$ splits into disjoint
union $\GD\sbp\cup\GD\sbm$, where $\GD\sbp$ is the set of
\emph{effective} roots (those realized by curves) and $\GD\sbm=-\GD\sbp$.
Furthermore, the set $\GD\sbp$ is closed with respect to positive linear
combinations and the function $\#\:\GD\sbp\to\N$ counting the number of
components of the curve representing a root $r\in\GD\spp$ is a well
defined homomorphism, in the sense that, whenever a root~$r$ is
decomposed into $\sum a_ir_i$ for some $r_i\in\GD\sbp$ and $a_i\in\N$,
one has $r\in\GD\sbp$ and $\#r=\sum a_i\,\#r_i$. (Note that, if $X$ is
algebraic, the roots $r\in\GD\sbp$ with $\#r=1$ define the walls of the
rational Dirichlet polyhedron of $\Aut\Pic X$ in $\HS(\Pic X\otimes\R)$
containing the fundamental class of a K\"ahler structure, see,
\eg,~\cite{PSh-Sh} or~\cite{DIK}. If $X$ is non-algebraic, they define
the walls of a distinguished camera of $\Pic X$.)

Let now $R\in\Pic X$ be a root system as in the statement and $\bar
R\supset R$ its saturation in $\Pic X$. Consider the subsets
$\bar\GD\sb\pm=\bar R\cap\GD\sb\pm$. They form a partition of the set of
roots of~$\bar R$, one has $\bar\GD\sbm=-\bar\GD\sbp$, and $\bar\GD\sbp$
is closed
with respect to positive linear combinations. Hence, there is a
unique camera~$C$ of~$\bar R$ such that $\bar\GD\sbp$ is the set of roots
positive with respect to~$C$ (see, \eg,~\cite{Bourbaki}); this means that
the roots $r_1,\dots,r_k\in\bar\GD\sbp$ defining the walls of~$C$ form a
basis of~$\bar R$ and each root $r\in\bar\GD\sbp$ is a positive linear
combination of the $r_i$'s. Hence, any root $r\in\bar\GD\sbp$ with
$\#r=1$ must be one of~$r_i$'s. Since $R$ is generated by such roots, one
has $R=\bar R$.
\endproof

\section{Are $K3$-surfaces quasi-simple?\enspace}\label{examples}

\subsection{$\KMS\G$ with walls}
Here, we construct an example of a geometric action of the group $G=\D_3$
(with $\Gr$ non-real and $\Gk\ne1$) whose associated space $\KMS\G$ has
more than two components, \ie, the action of $\Aut_GL$ on the set of
connected components of $\Per\G_0$ is not transitive. This shows that the
assumptions on the action in Theorem~\ref{Main-q-simplicity} cannot be
removed. However, the resulting Klein actions on $K3$-surfaces are not
diffeomorphic (see~\ref{non-diffeo}), \ie, they do not constitute a
counter-example to quasi-simplicity of $K3$-surfaces.

\proposition\label{example}
There is a homological type of $\D_3$-action on $L\cong3U\oplus2E_8$
realizable by six $\D_3$-equivariant deformation classes of
$K3$-surfaces. More precisely, there is a geometric action of $G=\D_3$
on~$L$ such that the corresponding moduli space $\KMS\G$ consists of
three pairs of complex conjugate connected components.
\endproposition

\proof
Fix a decomposition $L=P\oplus Q$, where $P\cong2U$ and $Q\cong
U\oplus2E_8$. Define a $\D_3$-action on~$L$ as follows. On~$Q$, the
$\Z_3$ part of~$\D_3$ acts trivially, and each nontrivial involution
of~$\D_3$ acts via multiplication by~$-1$. On~$P$, fix a basis $u_1$,
$v_1$, $u_2$ and $v_2$ so that $u_i^2=v_i^2=0$, $u_i\cdot v_i=1$, and
$u_i\cdot u_j=v_i\cdot v_j=u_i\cdot v_j=0$ for $i\ne j$. Choose an
order~$3$ element~$t$ and an order~$2$ element~$s$ in $\D_3$, and define
their action on~$P$ by the matrices
$$
T=\bmatrix\format\r&&\quad\r\\0&0&-1&0\\0&-1&0&-1\\1&0&-1&0\\0&1&0&0\endbmatrix,\qquad
S=\bmatrix\format\r&&\quad\r\\0&0&1&0\\1&0&0&1\\1&0&0&0\\0&1&-1&0\endbmatrix,
$$
respectively. Note that $\eL$ is trivial; hence, according
to~\ref{G-periods}, the constructed $\D_3$-action on~$L$ is realizable by
a Klein $\D_3$-action on a $K3$-surface.

The associated fundamental representation of the constructed action is
non-real. Hence, $\KMS\G \cong(\Per\G_0/\!\Aut_GL)\times S^0$, and it
suffices to show that $\Per\G_0/\!\Aut_GL$ has three connected
components.

One has $L^G=Q$ and $L_\Gr(\Z)=P$. The lattice $M\spp$
(the $(+1)$-eigenlattice of~$s$) is generated by $w_1=u_1+v_1+u_2$ and
$w_2=u_1+u_2-v_2$, and one has $w_1^2=2$, $w_2^2=-2$, and $w_1\cdot
w_2=0$. We assert that the only nontrivial automorphism of $M\spp$ that
extends to an equivariant automorphism of~$P$ is the multiplication
by~$-1$; thus, $\Aut_GP=\{\pm1\}$. Indeed, $\Aut M\spp$ consists of the
four automorphisms $w_1\mapsto\Ge_1w_1$, $w_2\mapsto\Ge_2w_2$, where
$\Ge_1,\Ge_2=\pm1$, and the equivariant extension to $P\otimes\Q$ is
uniquely given by the additional conditions $t(w_i)\mapsto\Ge_it(w_i)$.
If $\Ge_1\ne\Ge_2$, the extension is not integral.

Thus, the action of $\Aut_GL$ on $\HS\spp$ is trivial, the fundamental
domain $\GS_\Gr\spp$ coincides with $\HS\spp$, and, in view
of~\ref{case-cah}, one has
\smash{$\Per\G_0/\!\Aut_GL=(\tilde\GG_1\times\HS\spp)\sminus\GD$}, where
$\tilde\GG_1=\Int\GG_1/\!\Aut Q$ and $\GD$ is the union of a subset of
codimension~$\ge2$ and the hyperplanes
$\tilde\GG_1\times\psh_\Gr\sppm(v)$ defined by roots $v\in P$. (Since
$\dim\HS\spp=1$, each nonempty set $\psh_\Gr\sppm(v)$ is a hyperplane.)
Let $v\in P$ be a root and $v\sppm$ its projections to $V\sppm$. Since
$2v\sppm\in M_\Gr\sppm$ and $M_\Gr\spp$ has no vectors of square~$-4$,
the condition $\psh_\Gr\spp(v)\ne\varnothing$ implies that either
$v\spp=0$ (and then $(v\spm)^2=-2$), or $(v\spp)^2=-2$ (and then
$v\spm=0$), or $(2v\spp)^2=-8-(2v\spm)^2=-2$ or~$-6$. Each $M_\Gr\sppm$
contains, up to sign, one vector of square~$(-2)$ and two vectors of
square~$(-6)$. Comparing their images under~$J_\Gr$, one concludes that
the space $\HS\spp$ is divided into three components by the two walls
$\psh_\Gr\spp(w_2)$ and $\psh_\Gr\spp(2w_2-w_1)$.
\endproof

\subsubsection{}\label{new-action}
Before discussing this example in more details, introduce another
geometric $\D_3$-action on~$L$ with the same sublattice
$L^G=Q=U\oplus2E_8$. In the above notation, replace $S$ with the matrix
$$
S'=\bmatrix\format\r&&\quad\r\\1&0&-1&0\\0&1&0&0\\0&0&-1&0\\0&-1&0&-1\endbmatrix,
$$
and keep the rest unchanged. For the new action, one has $M_\Gr\sppm\cong
U(2)$. The only possible wall in~$\HS\spp$ is~$\psh_\Gr\spp(w\spp)$,
where $w\spp\in M_\Gr\spp$ is the only vector of square~$-4$. However,
$J_\Gr w\spp$ is not proportional to the vector $w\spm\in M_\Gr\spm$ of
square~$-4$; hence, the action is realized by a single $\D_3$-equivariant
deformation class of $K3$-surfaces.

In view of the following lemma, there are exactly two (up to isomorphism)
geometric $\D_3$-actions on~$L$ with $L^G\cong U\oplus2E_8$.

\lemma\label{two-Z3-actions}
Up to automorphism, there are three non-trivial $\Z_3$-actions on the
lattice $P\cong2U$\rom; their invariant sublattices are isomorphic to
either~$A_2$, or $A_2(-1)$, or~$0$. The last action admits two, up to
isomorphism, extension to a $\D_3$-action.
\endlemma

\proof
Let~$t\in\Z_3$ be a generator. Pick a primitive vector~$u_1$ of
square~$0$ and let $u_2=t(u_1)$. If $t(u_1)=u_1$ for any such~$u_1$, the
action is trivial. If $u_1\cdot u_2=a\ne0$, then $u_1$, $u_2$, and
$t^2(u_1)$ span a sublattice~$P'$ of rank three. In this case $a=\pm1$,
and the action is uniquely recovered using the fact that its restriction
to~$(P')\ort$ (a sublattice of rank one) is trivial. Finally, if
$u_1\cdot u_2=0$ and $u_1$, $u_2$ are linearly independent, then one must
have $t(u_2)=-u_1-u_2$. Completing $u_1$, $u_2$ to a basis $u_1$, $v_1$,
$u_2$, $v_2$ as in the proof of~\ref{example}, one can see that the
system $T^3=\id$, $\Gramm=T^*\Gramm T$ (where $T$ is the matrix of~$t$
and $\Gramm$ is the Gramm matrix) has a unique solution (the one
indicated in the proof of~\ref{example}).

Consider the last action and an involution $s\:P\to P$, $ts=st^{-1}$. The
invariant space~$M\spp$ of~$s$ is either~$U$, or~$U(2)$, or
$\<2\>\oplus\<-2\>$. The consideration above shows that the $\Z_3$-orbit
of any primitive vector~$u_1$ of square~$0$ is standard and spans a
sublattice of rank~$2$. Start from $u_1\in M\spp$ and complete it to a
basis $u_1$, $v_1$, $u_2$, $v_2$ as above. The set of solutions to the
system $TS=ST^{-1}$, $S^2=\id$, $\Gramm=S^*\Gramm S$ for the matrix~$S$
of~$s$ depends on one parameter~$a$, $s(v_2)=au_1-v_2$, and a change of
variables shows that only the values $a=0$ or~$1$ produce essentially
different actions (with $M\spp\cong U(2)$ or $\<2\>\oplus\<-2\>$,
respectively).
\endproof

\subsection{Geometric models}\label{example-discussion}
In this section, we give a geometric description
(via elliptic pencils) of the six families constructed in~\ref{example}.
At a result, at the end of the section we prove the following statement.

\proposition\label{non-diffeo}
All three pairs of complex conjugate deformation families constructed
in~\ref{example} differ by the topological type of the $\D_3$-action.
\endproposition

Fix a decomposition $Q=\Pic X\cong2E_8\oplus U$. Let $e_1',\ldots,e_8'$,
$e_1'',\ldots,e_8''$ be some standard bases for the $E_8$-components and
$u$, $v$ a basis for the $U$~component, so that $u^2=v^2=0$ and $u\cdot
v=1$. Under an appropriate choice of~$\Gg$ (a small perturbation of
$u+v$) the graph of $(-2)$-curves on~$X$ is the following:
$$
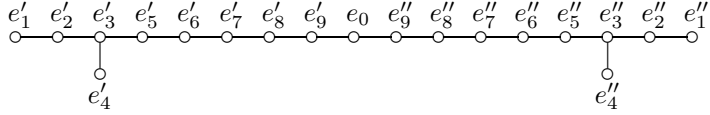
\def\.#1{@(\circ)@*<|><\textstyle\strut\ #1|>}
\DDtop
\.{e_1'}@---\.{e_2'}@---\.{e_3'}@---\.{e_5'}@---\.{e_6'}@---
 \.{e_7'}@---\.{e_8'}@---\.{e_9'}@---\.{e_0}@---
 \.{e_9''}@---\.{e_8''}@---\.{e_7''}@---\.{e_6''}@---\.{e_5''}@---
 \.{e_3''}@---\.{e_2''}@---\.{e_1''}\cr
@.@.@|||@.@.@.@.@.@.@.@.@.@.@.@.@|||\cr
@.@.@(\circ)@*<|><|\textstyle\strut{e_4'}>@.@.@.@.@.@.@.@.@.@.@.@.
 @(\circ)@*<|><|\textstyle\strut{e_4''}>\CR
\endDD
$$
Here $e_0=u-v$, $e_9'=v-2e_1'-4e_2'-6e_3'-3e_4'-5e_5'-4e_6'-3e_7'-2e_8'$,
and $e_9''=v-2e_1''-4e_2''-6e_3''-3e_4''-5e_5''-4e_6''-3e_7''-2e_8''$.

Consider the equivariant elliptic pencil $\pi\:X\to\Cp1$ defined by the
effective class~$v$. From the diagram above it is clear that the pencil
has a section~$e_0$ and two singular fibers of type~$\tilde E_8$, whose
components are $e_1',\ldots,e_9'$ and $e_1'',\ldots,e_9''$, respectively,
and has no other reducible singular fibers. (We use the same notation for
a $(-2)$-curve and for its class in~$L$.) Counting the Euler
characteristic shows that the remaining singular fibers are either
$4\tilde A_0^*$, or $2\tilde A_0^*+\tilde A_0^{**}$, or $2\tilde
A_0^{**}$. (Here $\tilde A_0^*$ and $\tilde A_0^{**}$ stand for a
rational curve with a node or a cusp, respectively.) In any case, at
least one of these singular fibers must also remain fixed under the
$\Z_3$-action; hence, the $\Z_3$-action on the base of the pencil has
three fixed points and thus is trivial. This implies, in particular, that
the pencil has no fibers of type $\tilde A_0^*$: the normalization of
such a fiber would have three fixed points (the two branches at the node
and the point of intersection with~$e_0$) and the $\Z_3$-action on it
and, hence, on the whole surface would have to be trivial. Thus, the
types of the singular fibers of the pencil are $2\tilde A_0^{**}+2\tilde
E_8$.

Let us study the action of~$\Z_3$ on the fibers of the pencil. Each fiber
has at least one fixed point: the point of intersection with~$e_0$. For
nonsingular fibers this implies that
\roster
\item
they all have $j$-invariant $j=0$ (as there is only one elliptic curve
admitting a $\Z_3$-action with a fixed point), and
\item
each nonsingular fiber has two fixed points more.
\endroster
Denote the closure of the union of these additional fixed points by~$C$.
This is a curve fixed under the $\Z_3$-action. In particular, it must
intersect the cuspidal fibers at the cusps. The action on the $\tilde
E_8$~singular fibers can easily be recovered starting from the points of
intersection with~$e_0$ and using the following simple observation: in
appropriate coordinates $(x,y)$ a generator $g\in\Z_3$ acts via
$(x,y)\mapsto(x,\Ge y)$ in a neighborhood of a point of a fixed curve
$y=0$, and via $(x,y)\mapsto(\Ge^2x,\Ge^2y)$ in a neighborhood of an
isolated fixed point $(0,0)$. (Here $\Ge$ is the eigenvalue of~$\Go$:
$g(\Go)=\Ge\Go$.) One concludes that the components~$e_3'$, $e_7'$,
$e_3''$, and~$e_7''$ are fixed, the intersection points of pairs of other
components are isolated fixed points, and $C$ intersects the $\tilde
E_8$~fibers at some points of~$e_1'$ and~$e_1''$. In particular, the
restriction $\pi\:C\to\Cp1$ is a double covering with four branch points;
hence, $C$ is a nonsingular elliptic curve.

Let $\tilde X$ be $X$ with isolated fixed points blown up and $\tilde
Y=\tilde X/\Z_3$. This is a rational ruled surface with two singular
fibers~$\tilde F'$, $\tilde F''$ (the images of the $\tilde E_8$ fibers
of~$X$), whose adjacency graphs are as follows:
$$
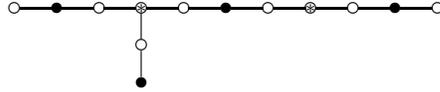
\def\.{@(\circ)}
\def\b{@(\bullet)}

\DDtop
\.@---\b@---\.@---\s@---\.@---\b@---\.@---\s@---\.@---\b@---\.\cr
 @.@.@.@|||\cr
 @.@.@.\.\cr
 @.@.@.@|||\cr
 @.@.@.\b\CR
\endDD
$$
(Here $\DDcenter@(\circ)\endDD$\,, $\DDcenter@(\bullet)\endDD$\,, and
$\DDcenter@(\circ\*)\endDD$ stand for a nonsingular rational curve of
self-intersection $-1$, $-3$, and~$-6$, respectively; an edge corresponds
to a simple intersection point of the curves.) The image~$\tilde R$ of
the section~$e_0$ has self-intersection~$(-6)$ and intersects the
rightmost curve in the graph; the image~$\tilde D$ of the section~$C$ has
self-intersection~$0$ and intersects the leftmost curve in the graph. The
branch divisor of the covering $\tilde X\to\tilde Y$ is $\tilde R+\tilde
D+\text{(the $(-6)$-components)}-\text{(the $(-3)$-components)}$.

Contract the singular fibers of~$\tilde Y$ to obtain a geometrically
ruled surface~$Y$. Denote by~$R$, $D$, $F'$, and~$F''$ the images
of~$\tilde R$, $\tilde D$, $\tilde F'$, and~$\tilde F''$, respectively.
The contraction can be chosen so that $R^2=0$, \ie,
$Y\cong\Cp1\times\Cp1$. Then $D^2=8$ and $D$ is a curve of bi-degree
$(2,2)$. It is tangent to~$F'$, $F''$, and $R$ passes through the
tangency points.

The above construction respects the $\D_3$-action on~$X$, and $Y$
inherits a canonical real structure in respect to which $D$, $R$, $F'$,
and $F''$, as well as the base of the pencil, are real; one has
$Y_{\R}=S^1\times S^1$.

Recall that, up to rigid isotopy, the embedding $D_{\R}\subset Y_{\R}$ is
one of the following:
\roster
\item\label{empty}
$D_\R$ is empty;
\item\label{oval}
$D_\R$ consists of one oval (a component contractible in~$Y_{\R}$);
\item\label{ovals}
$D_\R$ consists of two ovals;
\item\label{(0,1)}
$D_\R$ consists of two components, each realizing the class $(0,1)$ in
$H_1(Y_{\R})$;
\item\label{(1,0)}
$D_\R$ consists of two components, each realizing the class $(1,0)$ in
$H_1(Y_{\R})$;
\item\label{(1,1)}
$D_\R$ consists of two components, each realizing the class $(1,1)$ in
$H_1(Y_{\R})$.
\endroster
(The basis in $H_1(Y_{\R})$ is chosen so that $R_{\R}$ realizes $(1,0)$
and $F'_\R$ realizes $(0,1)$.) Now, one can easily indicate four
topologically distinct types of the action. Since $p'$ and~$p''$ are on
the same generatrix~$R$, the embedding $D_{\R}\subset Y_{\R}$ is either
\def\ibox#1{\hbox to\dimen8{#1\hss}}
\setbox0\hbox{(0)}\dimen8\wd0
\roster\widestnumber\item{b}
\item"\ibox{(a)}"
as in~\therosteritem{\ref{oval}}, or
\item"\ibox{(b)}"
as in~\therosteritem{\ref{ovals}} (the points $p'$, $p''$ are in the same
component of~$D_\R$), or
\item"\ibox{(c)}"
as in~\therosteritem{\ref{(0,1)}} (the points $p'$, $p''$ are in the
different components of~$D_\R$).
\endroster
In the latter case, there are two possibilities:
\setbox0\hbox{(c$''$)}\dimen8\wd0
\roster\widestnumber\item{c$''$}
\item"\ibox{(c$'$)}"
$F'_\R$ and $F''_\R$ belong to (the closure of) the same component of
$Y_\R\sminus D_\R$,
\item"\ibox{(c$''$)}"
$F'_\R$ and $F''_\R$ belong to (the closure of) distinct components of
$Y_\R\sminus D_\R$.
\endroster

Note that, according to Lemma~\ref{two-actions}, any model constructed
does necessarily realize either the action of~\ref{example}, or the
action of~\ref{new-action}.

The models of types (a) and (b) (resp., (a) and (c$'$)) can be joined
through a singular elliptic $K3$-surface whose desingularization has a
fiber of type $\tilde A_2$. In view of Proposition~\ref{family}, these
types realize the action of~\ref{example}. Hence, the remaining type
(c$''$) realizes the action of~\ref{new-action}.

\proof[Proof of \ref{non-diffeo}]
The surfaces in question are represented by the above models
of
types~(a), (b) and~(c$'$), which differ topologically: by the number of
components of $C_\R \cong D_\R$ and by whether $C_\R$ has a component
bounding a disk in $X_\R$.
\endproof

\subsection{The four families
in their Weierstra{\ss}{} form}\label{Weirstrass} Since the four families
constructed
above are
Jacobian fibrations (\ie, have sections), are isotrivial, and have
singular fibers of type $2\tilde A_0^{**}+2\tilde E_8$, their
Weierstra{\ss}{} equations are of the form
$$
y^2z=x^3+(u^2-v^2)^5p_2(u,v)z^3,
$$
where $(u:v)$ are homogeneous real coordinates in~$\Cp1$, $p_2$ is a
degree~$2$ homogeneous real polynomial with simple roots other than
$u=\pm v$, and $(x,y,z)$ are regarded as coordinate in the bundle
$\P(\O(6)\oplus\O(4)\oplus\O)$ over~$\Cp1(u:v)$. Isomorphisms between
such elliptic fibrations are given by projective transformations in
$\Cp1(u:v)$ and coordinates changes of the form $x\mapsto k^4x$,
$y\mapsto k^6y$, $z\mapsto z$, $u\mapsto ku$, $v\mapsto kv$, $k\in\R^*$.
By means of such isomorphisms the equation can be reduced to one of the
following four families:
$$
\alignat2
&y^2 z= x^3 + (u^2-v^2)^5 (u-cv)(u-\bar cv) z^3,&&\quad{c\ne\bar c},\\
&y^2 z= x^3 \pm (u^2-v^2)^5 (u-av)(u-bv) z^3,&&\quad{-1<a<b<1},\rlap{\quad and}\\
&y^2 z= x^3 + (u^2-v^2)^5 (u-av)(u-bv) z^3,&&\quad{-1<a<1<b}.
\endalignat
$$
The $\tilde E_8$ singular fibers are those with
$u^2=v^2$.
Each of the surfaces can be equipped with any of the two
$\D_3$-actions generated by the complex conjugation and the
multiplication of~$x$ by either $\exp(2\pi i/3)$ or $\exp(-2\pi i/3)$.

The exceptional family, \ie, that with the action  of~\ref{new-action},
is the one with the last equation. To see this, one can explicitly
construct two cycles in $M\sbm$ with square~$0$ and intersection~$2$. For
one of them, we pick a skew-invariant under the complex conjugation
circle $\xi$ in an elliptic fiber between $u=av$ and $u=v$ and drag it
along a loop in $\Cp1(u:v)$ around $u=-v$ and $u=av$. The other
(singular) cycle is constructed from a circle $\eta$ in the same fiber
with $T\eta=\bar\eta$, where $T$ is the monodromy operator about the
fiber $u=v$. We drag it along a loop around $u=v$ and pull its ends
together into the cusp of the fiber $u=av$.

Note that the real part of the double section of the surfaces in the
first family has only one connected component, so it correspond to the
series~(a). One component of the double section of the surfaces given by
the second equation with the sign~$-$ bounds a disc in the real part of
the surface, so it corresponds to series~(b). The same equation with the
sign~$+$ gives series~(c$'$).

Thus, one obtains another description of the six disjoint families
constructed in~\ref{example}. The bijection between the set of
isomorphism classes of $K3$-surfaces with a $\D_3$-action such that
$L^G=U\oplus2E_8$ and the set of surfaces given by the above four
equations (considered up to projective transformations of the base and
rescalings) can be used for an alternative proof of~\ref{example}.

\subsection{Distinct conjugate components with the same real~$\Gr$}
In this section, we construct an example of a geometric action~$\Gt$ of a
certain group $G=\wt T_{192}$ (with $\Gr\ne1$ real and $\Gk=1$) whose
moduli space has two distinct components interchanged by the conjugation
$X\mapsto\bar X$. Note that, since $\Gr$ is real, the components are {\bf
not} distinguished by the associated fundamental representations.

Recall that the group~$T_{192}$ can be described as follows. Consider the
form $\Phi(u,v)=u^4+v^4-2\sqrt{-3}\,u^2v^2$. Its group of unitary
isometries is the so called \emph{binary tetrahedral group}
$T_{24}\subset\grU(2)$; it can be regarded as a $\Z_3$-extension of the
Klein group~$Q_8=\{\pm1,\pm i,\pm j,\pm k\}\subset\Bbb H$. (Note that the
double projective line ramified at the roots of~$\Phi$ is a hexagonal
elliptic curve. An order three element of~$T_{24}$ can be given, \eg, by
the matrix
$$
q=\frac{1}{-1+i\sqrt3}\bmatrix-1-i&\phantom{-}1-i\\-1-i&-1+i\endbmatrix,
$$
whose determinant is $(-1+i\sqrt3)/2$.) The center of~$T_{24}$ is
$\{\pm1\}\subset Q_8$. Identify two copies of $T_{24}/Q_8\cong\Z_3$ via
$[q]\mapsto[q]^{-1}$ and let $T'$ be the fibered central product
$(T_{24}\times_{\Z_3}T_{24})/\{c_1=c_2\}$, where $c_1$ and $c_2$ are the
central elements in the two factors. Then $T_{192}$ is the semi-direct
product $T'\rtimes\Z_2$, the generator~$t$ of~$\Z_2$ acting via
transposing the factors.

Denote by~$\wt T_{192}$ the extension of $T_{192}$ by an element~$c$
subject to the relations $c^2=c_1=c_2$, $c^{-1}tc=c_1t$, and $ac=ca$ for
any~$a$ in either of the two copies of~$T_{24}\subset T_{192}$. Augment
this group via $\Gk\:\wt T_{192}\to\wt T_{192}/T_{192}=\ZZ$.

\proposition\label{non-real}
There is a geometric action of $G= \wt T_{192}$ on $L = 3U \oplus 2E_8$
such that the associated fundamental representation~$\Gr$ is real and the
corresponding moduli space $\KMS\G$ consists of a pair of conjugate
points~$X$, $\bar X$.
\endproposition

\proof
Consider the quartic $X\subset\Cp3$ given by the polynomial
$\Phi(x_0,x_1)+\Phi(x_2,x_3)$. According to Mukai~\cite{Mukai}, there is
a $T_{192}$-action on~$X$ with $\Gr = 1$. It can be described as follows.
The central product $(T_{24}\times T_{24})/\{c_1=c_2\}$ acts via block
diagonal linear automorphisms of $\Phi\oplus\Phi$, the two factors acting
separately in $(x_0,x_1)$ and $(x_2,x_3)$. The fundamental representation
of the induced action on~$X$ has order~$3$, and its kernel extends to a
symplectic $T_{192}$-action via the involution
$(x_0,x_1)@:<--->(x_2,x_3)$.

The described $T_{192}$-action on~$X$ extends to a $\wt T_{192}$-action,
the element $c\in\wt T_{192}$ acting via
$(x_0:x_1:x_2:x_3)\mapsto(ix_0:ix_1:x_2:x_3)$, so that $\Gr(c)=-1$.
Choosing an isometry $H^2(X) \to L$, one obtains a geometric \smash{$\wt
T_{192}$}-action on~$L$.

Fix a marking $H^2(X)=L$ and consider the twisted induced action on~$L$.
We assert that the corresponding period space consists of two points~$X$
and~$\bar X$, both admitting a unique embedding into~$\Cp3$ compatible
with a projective representation of~$\wt T_{192}$. Indeed, as it follows,
\eg, from the results of Xiao~\cite{Xiao}, for any action of the group
$G'=T_{192}$ with $\Gr=1$ one has $\rank\eL=19$; hence,
$\Per^{\smash{G'}}$ is a single point $\gw\subset L\otimes\R$ and
$K\GO^{\smash{G'}}=S^2$. Passing to $G=\wt T_{192}$ decomposes~$\gw$ into
$\ell=\gw\G$ and $\ell\ort$ and reduces $K\GO\G$ to a pair of points.
Since the action is induced from~$\Cp3$, the line~$\ell$ is generated by
an integral vector of square~$4$, and this is the only (primitive)
polarization of the surface compatible with the action.

It remains to show that $X$ does not admit an anti-holomorphic
automorphism commuting with~$\wt T_{192}$. Any such automorphism would
preserve~$\ell$ and, hence, would be induced from an anti-holomorphic
automorphism~$a$ of $\Cp3$. Since~$a$ commutes with $\wt T_{192}$, it
must fix the four intersection points of~$X$ with the line $C$ given by
$\{x_0=x_1=0\}$. In particular, $a$ must preserve $C$. On the other hand,
the roots of~$\Phi$ do not lie on a circle and, thus, cannot be fixed by
an anti-homography.
\endproof

\appendix{Finiteness and quasi-simplicity for $2$-tori}\label{tori}

\subsection{Klein actions on $2$-tori}\label{tori-intro}
In this section we prove analogs of Theorems~\ref{Main-Finiteness}
and~\ref{Main-q-simplicity} for complex $2$-tori (or just $2$-tori, for
brevity). The \emph{homological type} of a finite group~$G$ Klein action
on a $2$-torus~$X$ is the twisted induced action $\Gt_X\:G\to\Aut H^2(X)$
on the lattice $H^2(X)\cong3U$, considered this time up to conjugation by
{\bf orientation preserving} lattice automorphisms. As in the case of
$K3$-surfaces, one has $H^{2,0}(X)\cong\C$, and the action of~$\Gsbp$ on
$H^{2,0}(X)$ gives rise to a natural representation $\Gr\:\Gsbp\to\C^*$,
called the \emph{associated fundamental representation}. Both~$\Gt_X$
and~$\Gr$ are deformation invariants of the action; $\Gt_X$ is also a
topological invariant.

Our principal results for $2$-tori are the following two theorems.

\theorem[Finiteness Theorem]\label{Tori-Finiteness}
The number of equivariant deformation classes of complex $2$-tori with
faithful Klein actions of finite groups of uniformly bounded order
\rom(for any given bound\rom) is finite.
\endtheorem

\remark{Remark}
Note that the order of groups acting on $2$-tori and not containing pure
translations is bounded (\cf~\ref{Old-Finiteness} below).
In particular, there are finitely many deformation classes of such
actions.
\endremark

\theorem[Quasi-simplicity Theorem]\label{Tori-q-simplicity}
Let~$X$ and~$Y$ be two complex $2$-tori with diffeomorphic finite
group~$G$ Klein actions. Then either $X$ or $\wb X$ is $G$-equivariantly
deformation equivalent to~$Y$. If the associate fundamental
representation is trivial, then $X$ and~$\wb X$ are $G$-equivariantly
deformation equivalent.
\endtheorem

\theorem[Corollary]\label{Hyperelliptic}
The number of equivariant deformation classes of hyperelliptic surfaces
with faithful Klein actions of finite groups is finite. If~$X$ and~$Y$
are two hyperelliptic surfaces with diffeomorphic finite group~$G$ Klein
actions, then either $X$ or $\wb X$ is $G$-equivariantly deformation
equivalent to~$Y$.
\qed
\endtheorem

Recall that, after fixing a point~$0$ on a $2$-torus~$X$, one can
identify~$X$ with the quotient space $T_0(X)/H_1(X;\Z)$ and thus regard
it as a group. Then with each \anti-automorphism~$t$ of~$X$ one can
associate its linearization~$dt$ preserving~$0$, and, hence, any Klein
action~$\Gt$ on~$X$ gives rise to its linearization~$d\Gt$ consisting of
\anti-holomorphic autohomomorphisms of~$X$. As is known (see,
\ie,~\cite{Vinbergsurvey} or~\cite{Charlap}), the original action~$\Gt$
is uniquely determined by~$d\Gt$ and a certain element $a(\Gt)\in
H^2(G;H_1(X))=H^1(G;T_0(X)/H_1(X;\Z))$, the latter depending only on the
equivalence class of the extension $1\to H_1(X)\to\Cal G\to G\to 1$,
where $\Cal G$ is the lift of~$G$ to the group of \anti-holomorphic
transformations of the universal covering $T_0X$ of~$X$. In particular,
$a(\Gt)$ is a topological invariant.

Clearly, both the homological type of a Klein action~$\Gt$ and its
fundamental representation~$\Gr$ depend only on the linearization~$d\Gt$.
Since the group $H^2(G;H_1(X))$ is finite and $a(\Gt)$ is a topological
invariant, the general case of~\ref{Tori-Finiteness}
and~\ref{Tori-q-simplicity} reduces to the case of linear actions. Thus,
from now on, {\bf we consider only actions preserving~$0$.} All
\anti-automorphisms preserving~$0$ are group homomorphisms, and they all
commute with the automorphism $-\id\:X\to X$. For simplicity, we {\bf
always assume that $-\id\in G$}. For such actions, we prove
theorems~\ref{Old-Finiteness} and~\ref{Old-q-simplicity} below, which
imply~\ref{Tori-Finiteness} and~\ref{Tori-q-simplicity}.

\theorem\label{Old-Finiteness}
The number of equivariant deformation classes of complex $2$-tori with
faithful linear Klein actions of finite groups preserving $0$ is finite.
\endtheorem

\theorem\label{Old-q-simplicity}
Let~$X$ and~$Y$ be two complex $2$-tori with linear finite group~$G$
Klein actions of the same homological type.
Then either $X$ or $\wb X$ is $G$-equivariantly deformation
equivalent to~$Y$. If the associate fundamental representation is
trivial, then $X$ and~$\wb X$ are $G$-equivariantly deformation
equivalent.
\endtheorem

These theorems are proved at the end of Section~\ref{tori.spaces}.

\remark{Remark}
Note that, speaking about linear actions, Theorem~\ref{Old-q-simplicity}
is somewhat stronger than~\ref{Tori-q-simplicity}, as it also asserts
that the diffeomorphism type of a linear action is determined by its
homological type.
\endremark

\remark{Remark}
In the case of real actions (see \ref{results}), the surfaces $X$ and
$\bar X$ are obviously equivariantly isomorphic. The same remark applies
to~\ref{Hyperelliptic}, which gives us \emph{gratis} the following
generalization of the corresponding result by F.~Catanese and
P.~Frediani~\cite{CF} for real structures on hyperelliptic surfaces: Let
$X$ and $Y$ be two complex $2$-tori with real structures and with real
holomorphic $\Gsbp$-actions, so that the extended Klein actions of
$G=\Gsbp\times\ZZ$ have the same homological type and the same value of
$a(\theta)$. Then $X$ and~$Y$ are $G$-equivariantly deformation
equivalent.
\endremark

\subsection{Periods of marked $2$-tori}
Let $\GL$ be an oriented free abelian group of rank~$4$. Put
$L=\bigwedge^2\GL\spcheck$. The orientation of~$\GL$ defines an
identification $\bigwedge^4\GL\spcheck=\Z$ and turns $L$ into a lattice
via $\Jper\:L\otimes L\to\bigwedge^4\GL\spcheck=\Z$. It is isomorphic
to~$3U$. Denote $\Aut\spp L =\Aut L\cap\grSO\spp(L\otimes\R)$.

Let~$\CJ$ be the set of complex structures on~$\GL\otimes\R$ compatible
with the orientation of~$\GL$. Let, further, $\GO$ be the set of oriented
positive definite $2$-subspaces in $L\otimes\R$. As in~\eqref{GO-C}, one
can identify $\GO$ with the space $\{\Go\in
L\otimes\C\mid\Go^2=0,\;\Go\cdot\bar{\Go}>0\}/\C^*$. Both~$\CJ$ and~$\GO$
have natural structures of smooth manifolds. Let $\Jper\:\CJ\to\GO$ be
the map defined via $J\mapsto(x^1+iJ^*x^1)\wedge(x^2+iJ^*x^2)$, where
$J\in\CJ$, $J^*$ is the adjoint operator on~$L\spcheck$, and $x^1,x^2\in
L\spcheck\otimes\R$ are any two vectors generating $L\spcheck\otimes\R$
over~$\C$ (with respect to the complex structure~$J^*$).

The following statement is essentially contained in~\cite{PSh-Sh}
and~\cite{Shi}.

\proposition\label{JtoGO}
The map $\Jper\:\CJ\to\GO$ is a well defined diffeomorphism. The map
$\grSL(\GL)\to\Aut\spp L$, $\Gf\mapsto\wedge^2\Gf^*$, is an
epimorphism\rom; its kernel is the center $\{\pm1\}\subset\grSL(\GL)$. An
element $\Gf\in\grSL(\GL)$ commutes with a complex structure $J\in\CJ$ if
and only if its image $\wedge^2\Gf^*$ preserves $\Jper J$.
\endproposition

\proof
We will briefly indicate the proof. A simple calculation in coordinates
shows that the map $\Jper\:\CJ\to\GO$ is an immersion and generically
one-to-one. (Remarkably, the equations involved are partially linear.)
Since $\CJ$ and~$\GO$ are connected homogeneous spaces of the same
dimension, $\Jper$ is a diffeomorphism.

The map $\grSL(\GL\otimes\R)\to\grO(L\otimes\R)$,
$\Gf\mapsto\wedge^2\Gf^*$, is a homomorphism of Lie groups of the same
dimension. Hence, it takes the connected group $\grSL(\GL\otimes\R)$ to
the component of unity $\grSO\spp(L\otimes\R)$. The pull-back of
$\Aut\spp L\subset\grSO\spp(L\otimes\R)$ is a discrete subgroup of
$\grSL(\GL\otimes\R)$ containing $\grSL(\GL)$; on the other hand, the
latter is a maximal discrete subgroup (see~\cite{Ra}); hence, it
coincides with the pull-back.

The last statement follows from the naturallity of the construction: one
has $\Jper(\Gf J\Gf^{-1})=\wedge^2\Gf^*(\Jper J)$.
\endproof

A \emph{$1$-marking} of a $2$-torus~$X$ is a group isomorphism
$\Gf_1\:\GL\to H_1(X)$. We call a $1$-marking \emph{admissible} if it
takes the orientation of~$\GL$ to the canonical orientation of $H_1(X)$
(induced from the complex orientation of~$X$). A \emph{$2$-marking}
of~$X$ is a lattice isomorphism $\Gf\:H^2(X)\to L$. Since
$H^2(X)=\bigwedge^2H^1(X)$, every $1$-marking~$\Gf_1$ defines a
$2$-marking $\Gf=\wedge^2\Gf_1^*$. A $2$-marking is called
\emph{admissible} if it has the form $\wedge^2\Gf_1^*$ for some
admissible $1$-marking~$\Gf_1$. Any two admissible $1$-marking differ by
an element of $\grSL(\GL)$; in view of~\ref{JtoGO}, any two admissible
$2$-markings differ by an element of $\Aut\spp L$ and any admissible
$2$-marking has the form $\wedge^2\Gf_1^*$ for exactly two admissible
$1$-markings~$\Gf_1$.

 From now on by a $1$- (respectively, $2$-) marked torus we mean a
$2$-torus with a fixed admissible $1$- (respectively, $2$-) marking.
Isomorphisms of marked tori are defined in the obvious way
(\cf~\ref{period.maps}). Clearly, $1$-marked tori have no automorphisms;
the group of (marked) automorphisms of a $2$-marked torus is
$\{\pm\id\}$.

Consider the space $\Phi=\CJ\times(\GL\otimes\R)/\GL$ and the projection
$p\:\Phi\to\CJ$. The bundle $\Ker dp$ has a tautological complex
structure, which converts $p\:\Phi\to\CJ$ to a family of $1$-marked tori.
This family is obviously universal. In view of~\ref{JtoGO}, this implies
the following statement, called the global Torelli theorem for $2$-marked
tori.

\theorem\label{GTT-tori}
The  family $p\:\Phi\to\GO$ is a universal smooth family of $2$-marked
complex $2$-tori, \ie, any other smooth family $p'\:X\to S$ of $2$-marked
complex $2$-tori is induced from~$p$ by a unique smooth map $S\to\GO$.
\endtheorem

\subsection{Equivariant period spaces}\label{tori.spaces}
The following statement is similar to~\ref{K3->geometric}; it relies on
Proposition~\ref{JtoGO} and on the fact that a finite group action admits
an equivariant K\"ahler metric.

\proposition\label{tori->geometric}
Given a Klein action of a finite group~$G$ on a complex $2$-torus~$X$,
the twisted induced action $\Gt_X\:G\to\Aut H^2(X)$ is almost geometric
\rom(see~\ref{rho}\rom)\rom;
its image belongs to $\Aut\spp H^2(X)$.
\qed
\endproposition

Now, we proceed as in the case of $K3$-surfaces. Let $\Gt\:G\to\Aut\spp
L$ be an almost geometric action, and denote by $\GO\G\subset\GO$ the
fixed point set of the induced action $g\:\gv\mapsto\Gk(g)g(\gv)$,
$\gv\in\GO$. (As before, $-\gv$ stands here for~$\gv$ with the opposite
orientation.) Then the following holds.

\proposition\label{tori.G-periods}
The space $\GO\G$ is a fine period space of $2$-marked complex $2$-tori
with a Klein $G$-action compatible with~$\Gt$, \ie, it is the base of a
universal smooth family of $2$-marked complex $2$-tori with a Klein
$G$-action compatible with~$\Gt$.
\qed
\endproposition

\proposition\label{tori.G-connectedness}
Let $\Gk\:G\to\{\pm1\}$ be the augmentation and $\Gr\:\Gsbp\to S^1$ a
fundamental representation associated with~$\Gt$. If $\Gr=1$, then the
space $\GO\G$ is connected. If $\Gr\ne1$, then the space $\GO\G$ consists
of two components, which are transposed by the involution
$\gv\mapsto-\gv$.
\endproposition

\proof
As in the case of $K3$-surfaces, one can consider the contractible space
$\Per\G$ and sphere bundle $K\GO\G\to\Per\G$ and use the fibration
$K\GO\G\to\GO\G$ with contractible fibers.
\endproof

\proof[Proof of Theorems~\ref{Old-Finiteness} and~\ref{Old-q-simplicity}]
Theorem~\ref{Old-q-simplicity} follows from~\ref{tori.G-periods}
and~\ref{tori.G-connectedness}. In view of~\ref{Old-q-simplicity},
Theorem~\ref{Old-Finiteness} follows from the finiteness of the number of
homological types of faithful actions, \cf~\ref{main.proof}.
\endproof

\subsection{Comparing $X$ and $\bar X$}\label{compar}
As a refinement of Theorem~\ref{Tori-q-simplicity}, we show that in most
cases the $2$-tori~$X$ and~$\bar X$ are not equivariantly deformation
equivalent.

\proposition\label{order>2}
Consider a faithful finite group~$G$ Klein action on a complex
$2$-torus~$X$. Assume that $\Gsbp$ has an element of order $>2$ acting
non-trivially on holomorphic $2$-forms. Then $X$ is not $G$-equivariantly
deformation equivalent to~$\bar X$.
\endproposition

\proof
Let $g\in G$ be an element as in the statement. The assertion is obvious
if the associated fundamental representation~$\Gr$ is non-real. Thus, one
can assume that $\Gr$~is real and $\Gr(g)=-1$. A simple calculation
(using the fact that $g$ is orientation preserving, $\ord g>2$, and
$\wedge^2g^*$ has eigenvalue~$(-1)$ of multiplicity~$\ge2$) shows that in
this case the eigenvalues of the action of~$g$ on~$\GL$ are of the form
$\xi$, $\bar\xi$, $-\bar\xi$, $-\xi$ for some $\xi\notin\R$. Hence, there
is a distinguished square root $\sqrt
g\in\grSL(\GL\otimes\nomathbreak\R)$. (One chooses the arguments of the
eigenvalues in the interval $(-\pi,\pi)$ and divides them by~$2$.) The
automorphism $\wedge^2(\sqrt g)^*$ has order~$4$ on the (only)
$g$-skew-invariant $2$-subspace $\gv$; hence, it defines a distinguished
orientation on~$\gv$.
\endproof

\remark{Remark}
As a comment to the proof of Proposition~\ref{order>2}, we would like to
emphasize a difference between $K3$-surfaces and $2$-tori. Under the
assumptions of~\ref{order>2}, if~$\Gr$ is real, it is still possible that
there is an element~$a \in \Aut_G\spp L$ interchanging the two
points~$\gv$ and~$-\gv$ of $\GO\G$ (representing~$X$ and~$\bar X$).
However, unlike the case of $K3$-surfaces, this does not imply that $X$
and $\bar X$ are $G$-isomorphic; an additional requirement is that a lift
of $a$ to $\grSL(\GL)$ should commute with~$G$.
\endremark

\subsection{Remarks on symplectic actions}\label{symplectic-actions}
We would like to outline here a simple way to enumerate all symplectic
(\ie, identical on the holomorphic $2$-forms) finite group actions on
$2$-tori. (This result is contained in the classification by
Fujiki~\cite{Fujiki}, who calls symplectic actions special.) Our approach
follows that of Kond\=o~\cite{Kondo} to the similar problem for
$K3$-surfaces.

In view of~\ref{tori->geometric} and~\ref{tori.G-periods}, it suffices to
consider finite group actions on $L \cong 3U$ identical on a positive
definite $3$-subspace in $L \otimes \R$. Let $\Gt\:G\to\Aut\spp L$ be
such an action and $\eL = (L^G)\ort$. Then, $\eL$ is a negative definite
lattice of rank $\le 3$, and  the induced $G$-action on~$\eL$ is
orientation preserving and trivial on $\discr \eL$ (as so it is on
$\discr L^G$). Standard calculations with discriminant forms
(\cf~\cite{Kondo}) show that~$\eL$ can be embedded to~$E_8$ (the only
negative definite unimodular even lattice of rank~$8$), and the
$G$-action on~$\eL$ extends to~$E_8$ identically on $E\G_8 = (\eL)\ort
\subset E_8$. Since $\Aut E_8=W(E_8)$, the lattice~$\eL$ is the
orthogonal complement of a face of a camera of~$E_8$. Hence, $\eL$ is a
root system contained in~$A_3$, $A_2\oplus A_1$, or~$3A_1$, and
$G/\Ker\Gt$ is a subgroup of $W(\eL)\cap\grSO(\eL\otimes\R)$. It remains
to observe that any such lattice admits a unique (up to isomorphism)
embedding to~$L$ and, hence, the pair $\eL$, $G/\Ker\Gt\subset W(\eL)$
determines a $G$-action on~$L$ up to automorphism.

In particular, one obtains a complete list of finite groups~$G$ acting
faithfully and symplectically on $2$-tori. One has $\Ker\Gt=\{\pm\id\}$
and the group~$G/\Ker\Gt$ is a subgroup of
$W(\eL)\cap\grSO(\eL\otimes\nomathbreak\R)$ for $\eL=A_3$, $A_2\oplus
A_1$, or~$3A_1$, \ie, of $\frak A_4$ (alternating group on $4$~elements),
$\frak S_3$ (symmetric group on $3$~elements), or $\ZZ\oplus\ZZ$. Lifting
the action from~$L$ to~$\GL$, see~\ref{JtoGO}, one finds that $G$ is a
subgroup of the binary tetrahedral group~$T_{24}$, binary dihedral
group~$Q_{12}$, or Klein (quaternion) group~$Q_8$.

\enddocument